\documentclass[review,3p]{elsarticle}

\usepackage{lineno,hyperref}
\modulolinenumbers[1]

\usepackage{graphicx}
\usepackage{amsmath,amssymb,empheq,amsfonts,amsthm,bm}
\usepackage{booktabs}
\DeclareMathOperator{\E}{\mathbb{E}}
\DeclareMathOperator{\V}{\mathbb{V}}
\newtheorem{theorem}{Theorem}
\newcommand\norm[1]{\left\lVert#1\right\rVert}

\usepackage{floatrow}
\usepackage{caption}
\usepackage{subcaption}
\usepackage[usenames,dvipsnames]{pstricks}
\usepackage{booktabs}
\usepackage{caption}
\usepackage{subcaption}
\usepackage{multirow}
\usepackage[justification=centering]{caption}
\usepackage{physics}
\usepackage{tikz}
\usepackage{pgfplots}
\usepackage{pgf}
\usepackage{soul}

    \pgfplotsset{compat=newest}
    \newlength\figureheight
    \newlength\figurewidth
\usepackage{tikz}
\newcommand\fracfun[1]{\textrm{frac}\left(#1\right)}

\DeclarePairedDelimiter\floor{\lfloor}{\rfloor}

\usetikzlibrary{patterns}
\usetikzlibrary{calc}

\journal{Probabilistic Engineering Mechanics}
\begin{document}

\begin{frontmatter}

\title{h- and p-refined Multilevel Monte Carlo Methods for \\ Uncertainty Quantification in Structural Engineering}

\author[CW]{P. Blondeel\corref{cor1}} 
\ead{philippe.blondeel@kuleuven.be}

\author[CW]{P. Robbe} 

\author[BOKU]{C. Van hoorickx} 

\author[BOKU]{G. Lombaert} 

\author[CW]{S. Vandewalle} 

\address[CW]{KU Leuven, Department of Computer Science, NUMA Section,\\ Celestijnenlaan 200A, 3001 Leuven, Belgium}
\address[BOKU]{KU Leuven, Department of Civil Engineering, Structural Mechanics Section,\\ Kasteelpark Arenberg 40, 3001 Leuven, Belgium}

\cortext[cor1]{Corresponding author}

\begin{abstract}
Practical structural engineering problems  are often characterized by significant uncertainties. Historically, one of the prevalent methods to account for this uncertainty  has been the standard Monte Carlo (MC) method. Recently, improved sampling methods have been proposed, based on the idea of variance reduction by employing a hierarchy of mesh refinements. We combine an h- and p-refinement hierarchy with the Multilevel Monte Carlo (MLMC) and Multilevel Quasi-Monte Carlo (MLQMC) method. We investigate the applicability of these novel combination methods on three  structural engineering problems, for which the uncertainty resides in the Young's modulus: the static response of a cantilever beam with elastic material behavior, its static response with elastoplastic behavior, and its dynamic response with elastic behavior. The uncertainty is either modeled by means of  one random variable sampled from a univariate Gamma distribution or with multiple random variables sampled from a gamma random field. This random field  results from  a truncated Karhunen--Lo\`eve (KL) expansion.  In this paper, we compare the computational costs of these Monte Carlo methods. We demonstrate that MQLMC and MLMC have a significant speedup with respect to MC, regardless of the mesh refinement hierarchy used. We empirically demonstrate that the MLQMC cost is optimally proportional to $\epsilon^{-1}$ under certain conditions, where $\epsilon$ is the tolerance on the root-mean-square error (RMSE). In addition, we show that, when the uncertainty is modeled as a random field, the multilevel methods combined with p-refinement have a significant lower computation cost than their counterparts based on h-refinement. We also illustrate the effect the uncertainty models have on the uncertainty bounds in the solutions. An uncertain Young's modulus modeled as a single random variable has much larger uncertainty bounds on its solution than an uncertain Young's modulus modeled as a random field. \end{abstract}



\end{frontmatter}

\nolinenumbers

\section{Introduction}

There is an increasing need to accurately simulate and compute solutions to engineering problems whilst taking into account model uncertainties. Methods for  uncertainty quantification and propagation in structural engineering can be categorized into two groups: non-sampling methods and sampling methods. Examples of non-sampling methods are the perturbation method  and the Stochastic Galerkin Finite Element method. The perturbation method is based on a Taylor series expansion approximating the mean and variance of the solution \cite{Kleiber}.
The method is quite effective, but its use is restricted to models with a limited number of relatively small uncertainties.
The Stochastic Galerkin method, first proposed by Ghanem and Spanos \cite{GhanemAndSpanos}, is based on a spectral representation in the stochastic space. It transforms the uncertain coefficient partial differential equation (PDE) problem by means of a Galerkin projection technique into a large coupled system of deterministic PDEs. This method allows for somewhat larger numbers of uncertainties and is quite accurate. However, it is highly intrusive and memory demanding, making its implementation cumbersome and restricting its use to rather low stochastic dimensions.

Sampling methods, on the other hand, are typically non-intrusive. Each sample corresponds to a deterministic solve for a set of specified parameter values. Two particularly popular examples are the Stochastic Collocation method \cite{Bab} and the Monte Carlo (MC) method \cite{Fishman}.
The former samples a stochastic PDE at a carefully selected multidimensional set of collocation points. After this sampling, a Lagrange interpolation is performed leading to a polynomial response surface. From this, the relevant stochastic characteristics can easily be computed in a post-processing step. However, as is also the case for Stochastic Galerkin, the Stochastic Collocation method suffers from the curse of dimensionality: the computational cost grows exponentially with the number of random variables considered in the problem.
The MC method on the other hand, selects its samples randomly and does not suffer from the curse of dimensionality. A drawback is its slow convergence as a function of the number of samples taken.
The convergence of Monte Carlo can be accelerated in a variety of ways. For example, alternative non-random selections of sampling points can be used, as in Quasi-Monte Carlo \cite{Caflish,Niederreiter} and Latin Hypercube \cite{Loh} sampling methods. Also, variance reduction techniques, such as Multilevel Monte Carlo (MLMC) \cite{Giles}, Multilevel Quasi-Monte Carlo (MLQMC) \cite{Giles3} and its generalizations, see, e.g., \cite{PJ, PJ2}, can speed up the method. These improved Monte Carlo methods are based on a hierarchy of increasing resolution meshes where samples on coarser meshes are computationally less expensive than on finer meshes.
As a side note, we  mention that there also exist hybrid variants which exhibit both a sampling and non-sampling character. This type of methods combine, for example, the Stochastic Finite Element methodology with Monte Carlo sampling or a multi-dimensional cubature method, see, e.g., \cite{Ghanem_Hyrbid,Acharjee}.

Monte Carlo methods have since long been used in the field of structural engineering, for example in problems of structural dynamics \cite{Masanobu} or in elastoplastic problems where the structure's reliability is assessed \cite{Pulido}.
In this work we combine the MLMC and MLQMC method with an h- and p-refinement mesh hierarchy. These combinations are then applied to a structural engineering problem discretized by means of the finite element method. The problems are defined as the static response of a cantilever beam with elastic material behavior, its static response with elastoplastic behavior, and its dynamic response with elastic behavior. The uncertainty resides in the  Young's modulus.  We consider two different representations of the uncertainty: a homogeneous and a heterogeneous one. The homogeneous representation consists of a random variable sampled from a univariate gamma distribution. For the heterogeneous representation,  we do not use a Gaussian or lognormal random field, as is often the case in other works in the literature, but we use a gamma random field to model the uncertainty. This field is obtained by combining a Karhunen--Lo\`eve expansion with a memoryless transformation. We illustrate the effect of the uncertainty model on the uncertainty bounds of the solution. These uncertainty bounds are computed from the resulting probability density function of the solution.  The obtained MLMC and MLQMC results combined with h- and p-refinement will be compared in terms of computational cost  with results from a standard MC simulation.
\\
This paper is structured as follows. In section 2, we formulate the mathematical model, introduce the problem statement and describe how the uncertainty is modeled. Section 3 recalls the MLMC and MLQMC methods, and provides some additional algorithmic implementation details. In Section 4, numerical results are presented. First, we illustrate  the uncertainty propagation towards the solution for the static elastic, static  elastoplastic and the dynamic elastic response. Second, the performance of standard MC is compared with the performances of MLMC and MLQMC for the static  response, both for an elastic and an elastoplastic material model. Both h- and p-refinement schemes are considered.  The fifth and last section offers concluding remarks and details some paths for further research.

\section{The mathematical model}
\subsection{Beam models and material parameters}
\begin{figure}[H]
\begin{subfigure}[b]{0.54\textwidth}
\centering
\includegraphics[height=1.9cm]{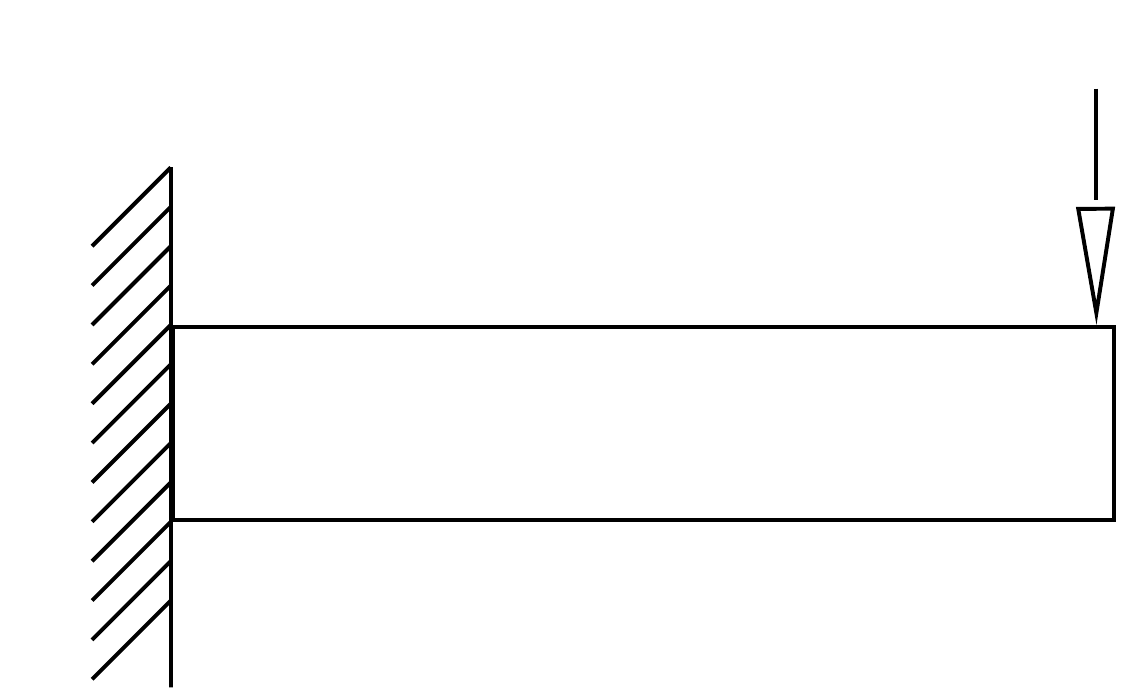}
\end{subfigure}
\begin{subfigure}[b]{0.45\textwidth}
\centering
\includegraphics[height=1.9cm]{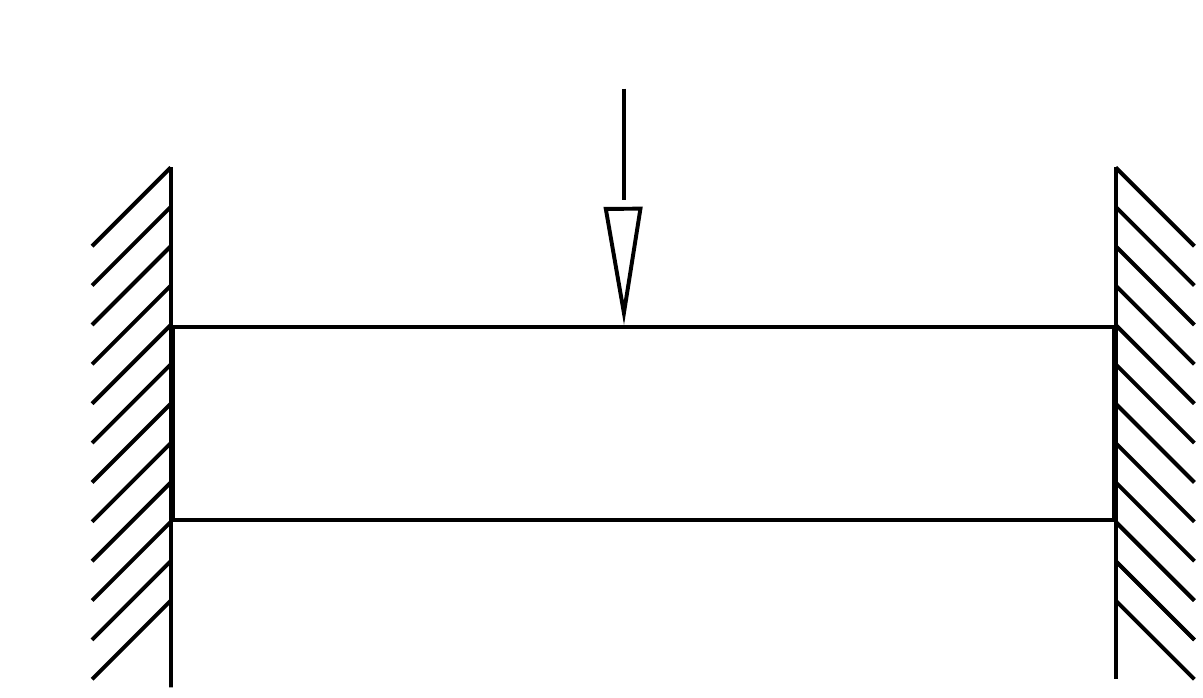}
\end{subfigure}
\caption{Cantilever beam loaded at its right end (left) and beam clamped at both ends loaded in the middle (right).}
\label{fig:bean_configurations}
\end{figure}

The considered engineering problem is the response of a cantilever beam clamped at one side and a beam clamped at both sides  as seen in Fig.\ref{fig:bean_configurations}, assuming plane stress. We consider three different responses. First, we consider the spatial displacement of a concrete beam with an elastic material model, clamped at both ends (static elastic case). Secondly, we consider the spatial displacement of a steel beam with an elastoplastic material model, clamped at both ends (static elastoplastic case). Finally, we consider the frequency response of a concrete beam, clamped at its left end (dynamic elastic case). An overview is given in Tab.\,\ref{Tab:cases}.  The dimensions of the beam are $2.5\,\mathrm{m}$ (length) by $0.25\,\mathrm{m}$ (height) by $1\,\mathrm{m}$ (width)  for the elastic cases and $10^{-3}\,\mathrm{m}$ (width) for the elastoplastic case. The material parameters of the concrete are as follows: a mass density of $2500\,\mathrm{kg/m^3}$, a Poisson ratio of $0.15$ and a Young's modulus subject to some uncertainty, as specified below. The material parameters of the steel are as follows: a yield strength of $240\,\mathrm{MPa}$, a Poisson ratio of $0.25$ and a Young's modulus subjected to uncertainty, as specified below. In order to model the material uncertainty, two uncertainty models will be considered. The first model is a homogeneous Young's modulus characterized by means of a single random variable. The second  model is a heterogeneous Young's modulus represented as a random field. Both uncertainty models will be used to compute the stochastic characteristics in all cases.

\begin{table}[H]
\centering
\scalebox{0.7}{
\begin{tabular}{cccccc}
\toprule
 Case &  Material &  Configuration &  Response &   Domain & Uncertainty \\
    \cmidrule(rl{4pt}){1-6}  
\multirow{2}{*}{Static Elastic} & \multirow{2}{*}{Concrete} & \multirow{2}{*}{Fig.\,\ref{fig:bean_configurations} (right)} & \multirow{2}{*}{Spatial displacement}& \multirow{2}{*}{Elastic}&Homogeneous\\
& && & &Heterogeneous \\
    \cmidrule(rl{4pt}){1-6}  
 \multirow{2}{*}{Static Elastoplastic} & \multirow{2}{*}{Steel} & \multirow{2}{*}{Fig.\,\ref{fig:bean_configurations} (right)} & \multirow{2}{*}{Spatial displacement}& \multirow{2}{*}{Elastoplastic}&Homogeneous\\
& && & &Heterogeneous \\
    \cmidrule(rl{4pt}){1-6}  
 \multirow{2}{*}{Dynamic Elastic} & \multirow{2}{*}{Concrete} & \multirow{2}{*}{Fig.\,\ref{fig:bean_configurations} (left)} & \multirow{2}{*}{Frequency response}& \multirow{2}{*}{Elastic}&Homogeneous\\
& && & &Heterogeneous \\
\bottomrule
\end{tabular}}
\caption{Overview of the different considered cases.}
\label{Tab:cases}
\end{table}

\subsubsection{The homogeneous model}
Following \cite{YongLiu}, we opt to describe the Young's modulus in the homogeneous model by means of a univariate gamma distribution. This distribution is characterized by a shape parameter $\alpha$ and a scale parameter $\beta$, and its probability density function given by
\begin{linenomath*}
\begin{equation}
f(x|\alpha,\beta) = \dfrac{1}{\beta^{\alpha} \Gamma(\alpha)} x^{\alpha-1} e^{\left(-\dfrac{x}{\beta}\right)}.
\end{equation}
\end{linenomath*}

\begin{figure}[H]
\centering
\begin{subfigure}[b]{0.54\textwidth}
\scalebox{0.45}{
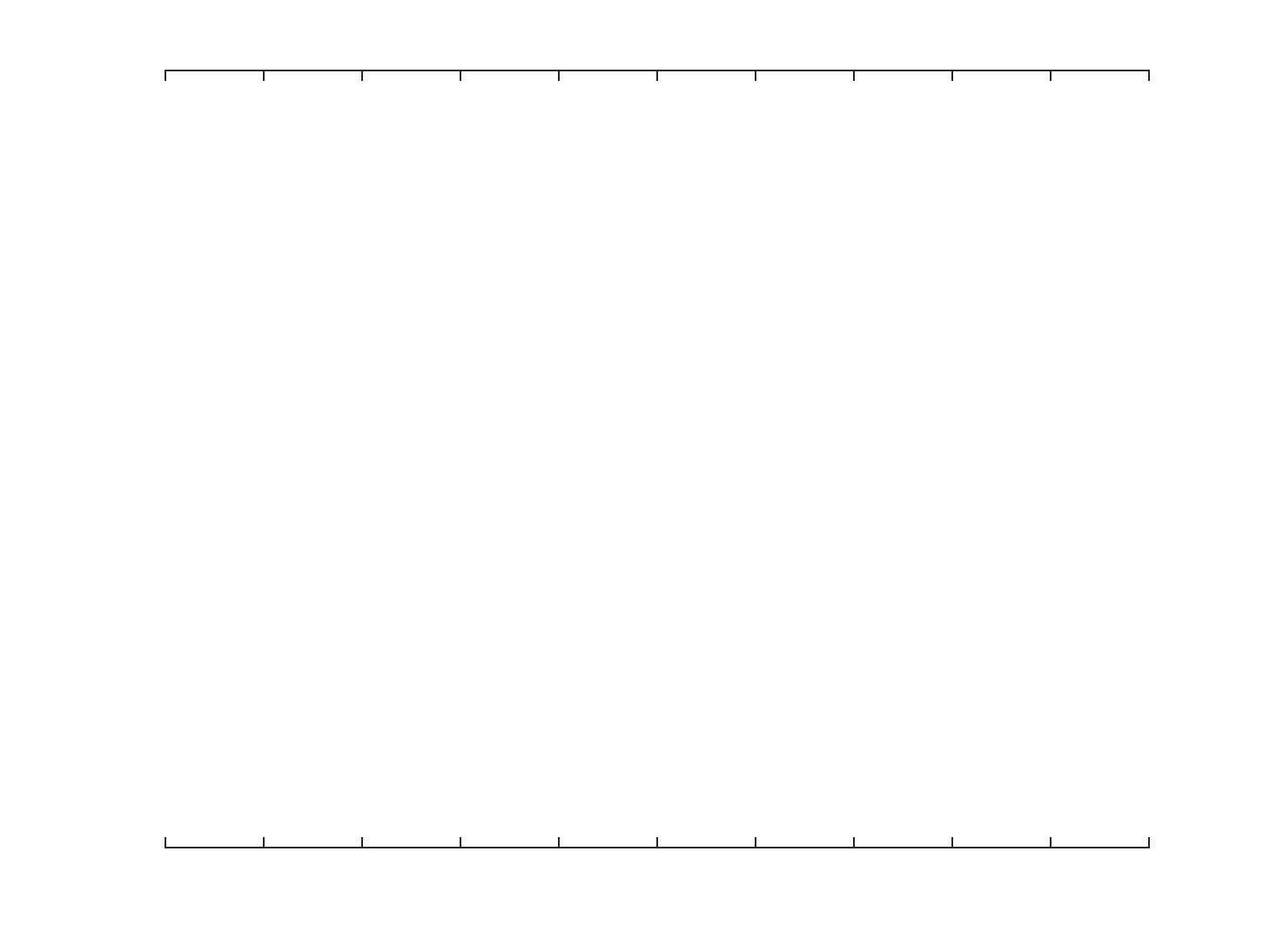}
\end{subfigure}
\begin{subfigure}[b]{0.45\textwidth}
\scalebox{0.45}{
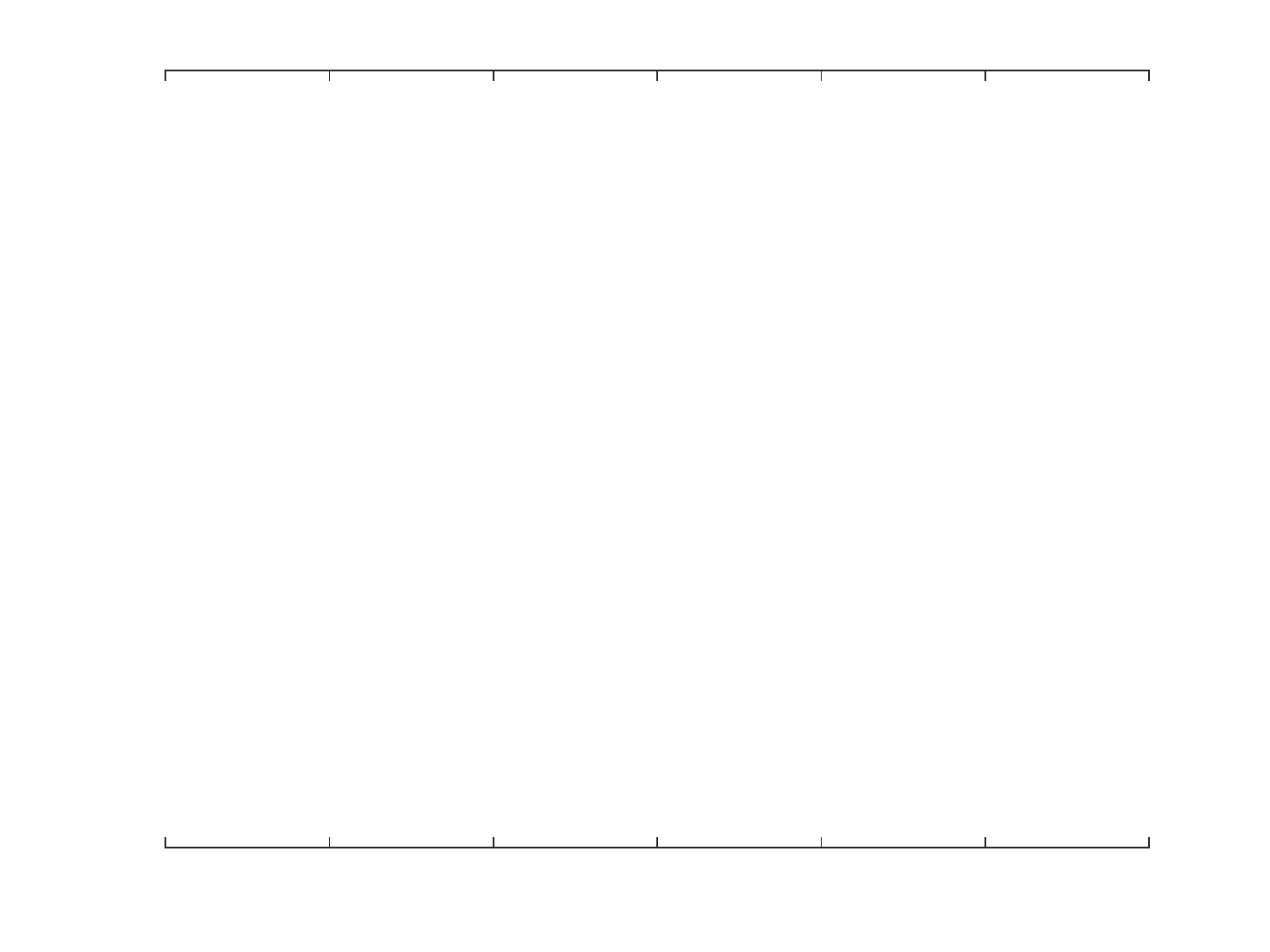}
\end{subfigure}
\centering
\caption{Probability density function for Young's modulus as a univariate distribution for concrete (left) and steel (right). Shown also is the mean $\mu$  and the standard deviation $\sigma$.}
\label{fig:pdf_gam_prior}
\end{figure}

The corresponding mean value and variance can be computed  as $\mu=\alpha \beta$ and
$\sigma^2=\alpha \beta^2$ respectively. In this paper, we  select $\alpha\!=\!7.1633$ and $\beta\!=\!4.1880 \times 10^9$ in order to model the material uncertainty for the concrete beam, which are based on values coming from \cite{Ellen}. This leads to a mean of $30\,\mathrm{GPa}$ and a standard deviation of $11.2\,\mathrm{GPa}$. For modeling the material uncertainty in the steel beam, we select $\alpha\!=\!934.2$ and $\beta\!=\!0.214\times 10^9$, see \cite{Hess}. This gives a mean of $200\,\mathrm{GPa}$ and a standard deviation of $6.543\,\mathrm{GPa}$. The gamma distribution for both materials is plotted in Fig.\,\ref{fig:pdf_gam_prior}.

\subsubsection{The heterogeneous model}
The Young's modulus with spatially varying uncertainty will be represented by means of a (truncated) gamma random field. The construction of this random field is done by means of a classic, two-step process. First, a (truncated) Gaussian random field is generated, using a Karhunen--Lo\`eve (KL) expansion \cite{Loeve}. Next, this Gaussian random field is transformed into a gamma random field with a memoryless transformation \cite{Grigoriu}.

Consider a Gaussian random field $Z(\mathbf{x},\omega)$, where $\omega$ is a random variable, with exponential covariance kernel,
\begin{linenomath*}
\begin{equation}\label{eq:covariance_kernel}
C(\mathbf{x},\mathbf{y}):= \sigma^2 \exp\left(-\dfrac{\norm{\mathbf{x}-\mathbf{y}}_p}{\lambda}\right)\,.
\end{equation}
\end{linenomath*}
We select the 2-norm ($p\!=\!2$), a correlation length $\lambda\!=\!0.3$ and a standard deviation $\sigma \!=\!1.0$. The corresponding KL expansion can then be formulated as follows:
\begin{linenomath*}
\begin{equation}
Z(\mathbf{x},\omega)=\overline{Z}(\mathbf{x},.)+\sum_{n=1}^{\infty}  \sqrt{\theta_n} \xi_n(\omega) b_n(\mathbf{x})\,.
\label{eq:KLExpansion}
\end{equation}
\end{linenomath*}
 $\overline{Z}(\mathbf{x},.)$ denotes the mean of the field, and is set to zero. The $\xi_n(\omega)$ denote i.i.d.\,standard normal random variables.
The symbols $\theta_n$  and $b_n(\mathbf{x})$ respectively denote  the eigenvalues and eigenfunctions of the covariance kernel corresponding to Eq.\,\eqref{eq:covariance_kernel}, which are found by solving the following eigenvalue problem:
\begin{linenomath*}
\begin{equation}
\int_D C(\mathbf{x},\mathbf{y})b_n({\mathbf{y}})\mathrm{d}\mathbf{y} = \theta_n b_n({\mathbf{x}}).
\end{equation}
\end{linenomath*}

These can be approximated  by means of a numerical collocation scheme, i.e., by solving 
\begin{linenomath*}
\begin{equation}
\int_D C(\mathbf{x}_k,\mathbf{y})b_n({\mathbf{y}})\mathrm{d}\mathbf{y} = \theta_n b_n({\mathbf{x}_k}), \quad k=1,2,\dots, M,
\label{eq:Fredholm}
\end{equation}
\end{linenomath*}
in some well-chosen integration points $\mathbf{x}_k$. Following the Nystr\"om method \cite{Atkinson}, the integral in Eq.\,\eqref{eq:Fredholm}, is approximated by a numerical integration scheme which uses the collocation points as quadrature nodes:
\begin{linenomath*}
\begin{equation}
\sum_{q=1}^M w_q C(\mathbf{x}_k,\mathbf{y}_q)\widetilde{b_n}({\mathbf{y}_q}) = \widetilde{\theta}_n \widetilde{b_n}({\mathbf{x}_k}), \quad k=1,2,\dots, M.
\label{eq:colloc}
\end{equation}
\end{linenomath*}
In matrix notation, this becomes
\begin{linenomath*}
\begin{equation}
\Sigma W \widetilde{B}_n = \widetilde{\theta}_n \widetilde{B}_n,
\label{eq:eigvalprob}
\end{equation}
\end{linenomath*}
where $\Sigma$ is a symmetric positive semi-definite matrix with entries $\Sigma_{k,q}=C(\mathbf{x}_k,\mathbf{y}_q)$, $W$ is a diagonal matrix containing the weights $w_q$ on its diagonal and $B$ is a vector with entries $B_{n,q} = b_n(\mathbf{x}_{q})$. The matrix eigenvalue problem, Eq.\,\eqref{eq:eigvalprob}, can be reformulated in an equivalent matrix eigenvalue problem
\begin{linenomath*}
\begin{equation}
\Psi \widetilde{B}^*_n = \widetilde{\theta}_n \widetilde{B}^*_n,
\end{equation}
\end{linenomath*}
where $\widetilde{B}^*_n = \sqrt{W} \widetilde{B}_n$ and $\Psi = \sqrt{W} \Sigma  \sqrt{W}$. $\Psi$ is symmetric positive semi-definite. This implies that the eigenvalues $\widetilde{\theta}_n$ are nonnegative real values and the eigenvectors $\widetilde{B}_n^*$ are orthogonal to each other.
Using Eq.\,\eqref{eq:colloc}, the Nystr\"om interpolation value for the eigenfunctions $b_n(\mathbf{x})$ is obtained:
\begin{linenomath*}
\begin{equation}
\widetilde{b}_n(\mathbf{x}) = \frac{1}{\widetilde{\theta}_n}\sum_{q=1}^M \sqrt{w_q} \widetilde{B}_{n,q}^* C(\mathbf{x},\mathbf{y}_q),
\end{equation}
\end{linenomath*}
where $\widetilde{B}_{n,q}^*$ stands for the q-th element of eigenvector $\widetilde{B}_n^*$. These eigenvalues and eigenfunctions, can after a suitable normalization, be used as an approximate eigenpair in the KL expansion.

In an actual implementation, the number of KL-terms in Eq.\,\eqref{eq:KLExpansion} is truncated to a finite value $s$, i.e.,
\begin{linenomath*}
\begin{equation}
Z(\mathbf{x},\omega)=\overline{Z}(\mathbf{x},.)+\sum_{n=1}^{s}  \sqrt{\theta_n} \xi_n(\omega) b_n(\mathbf{x})\,.
\end{equation}
\end{linenomath*}
 This  number of uncertain parameters depends on the magnitude and on the decay rate of the successive eigenvalues.

\begin{figure}[H]
\centering
\begin{subfigure}[b]{0.54\textwidth}
\scalebox{0.45}{
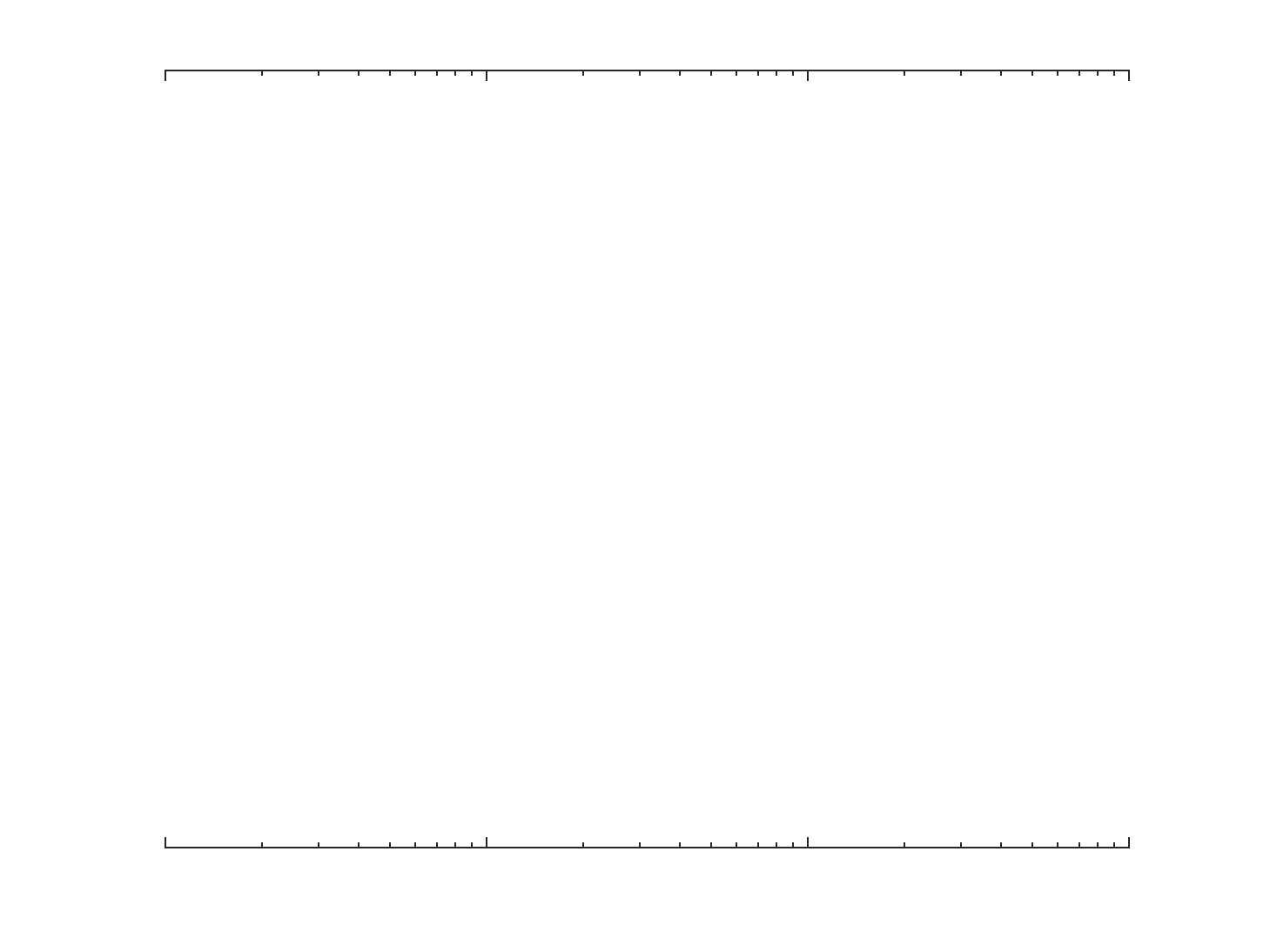}
\end{subfigure}
\begin{subfigure}[b]{0.45\linewidth}
\scalebox{0.45}{
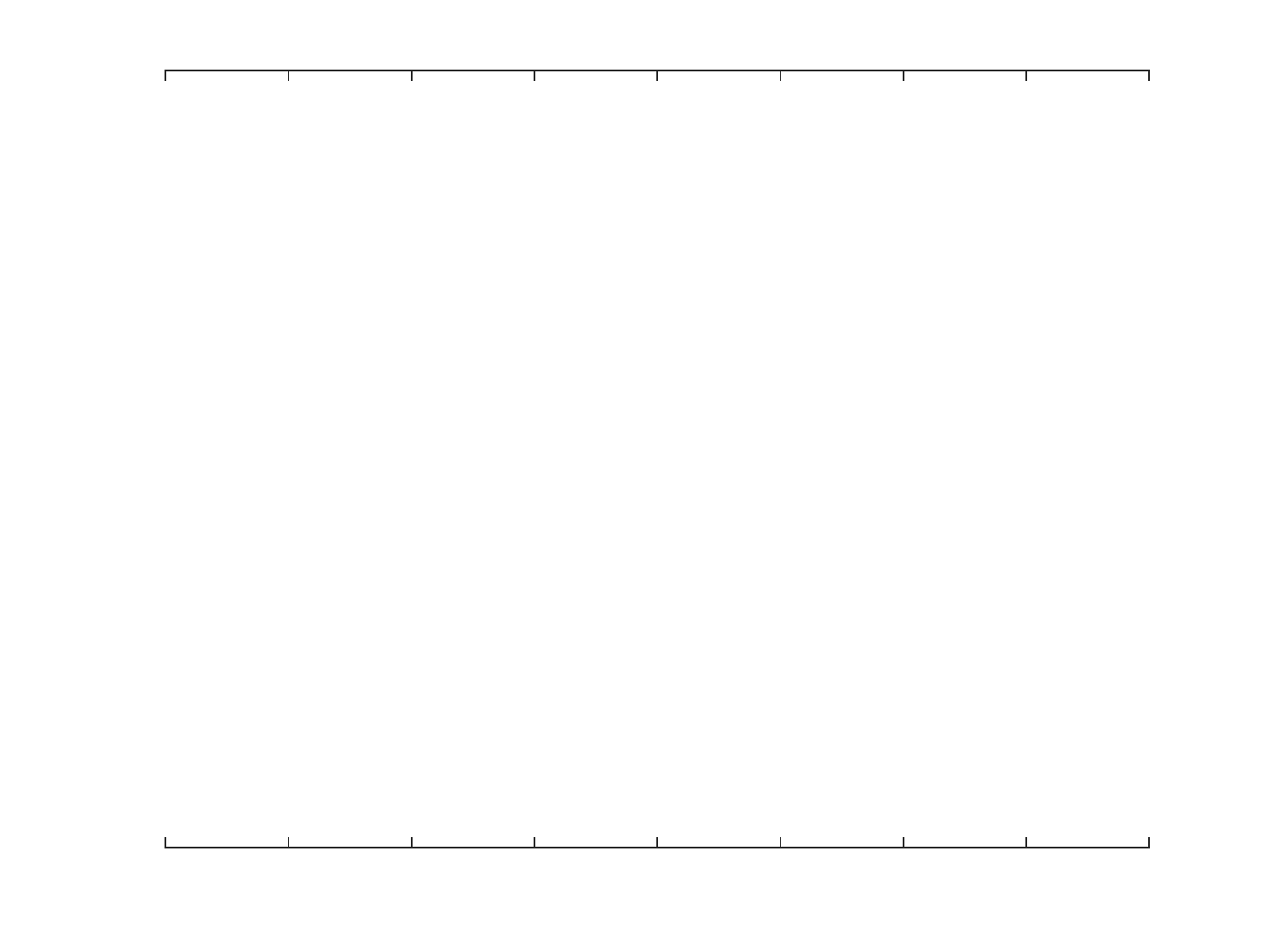}
\end{subfigure}
\caption{Magnitude of the eigenvalues and their cumulative sum (left) and memoryless transformation used to generate the gamma random field (right).}
\label{fig:Eigenvalf}
\end{figure}

The eigenvalues for the exponential covariance function, Eq.\,\eqref{eq:covariance_kernel}, are plotted in Fig.\,\ref{fig:Eigenvalf} (left). The percentage of the variance that is accounted for as a function of the number of included eigenvalues corresponds to the cumulative sum of the eigenvalues, and is also illustrated in the figure. A cumulative sum of 1.0 corresponds to 100\% of the variance of the field being accounted for. Inclusion of the first 101 KL-terms is  sufficient to represent 90\% of the variance of the random field.

\begin{figure}[H]
\begin{subfigure}[b]{0.54\textwidth}
\centering
\scalebox{0.45}{
	\makebox[\textwidth]{
	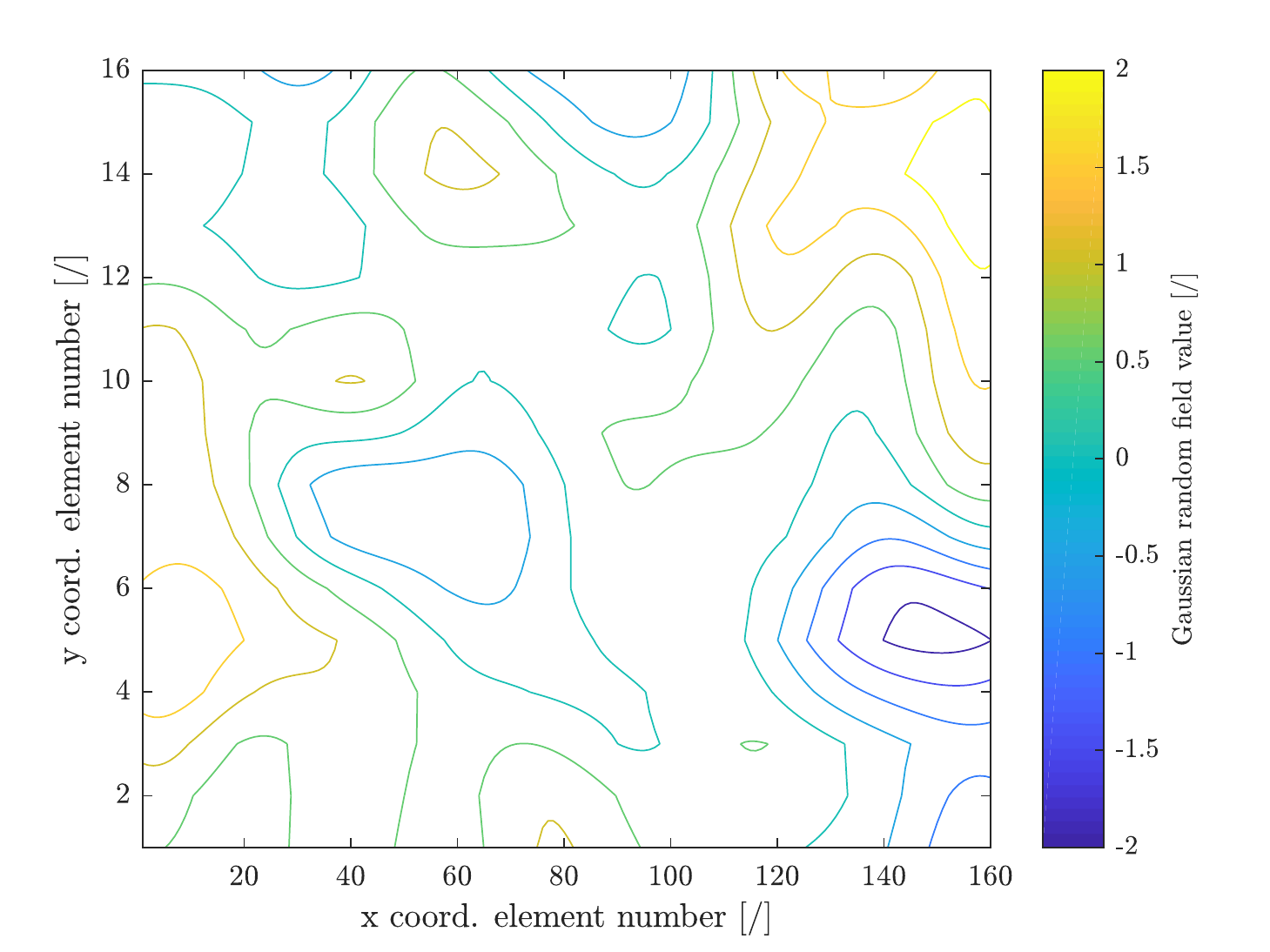}}
\label{fig:Gaussian Random field}
\end{subfigure}
\begin{subfigure}[b]{0.45\linewidth}
\centering
	\scalebox{0.45}{
		\makebox[\textwidth]{
	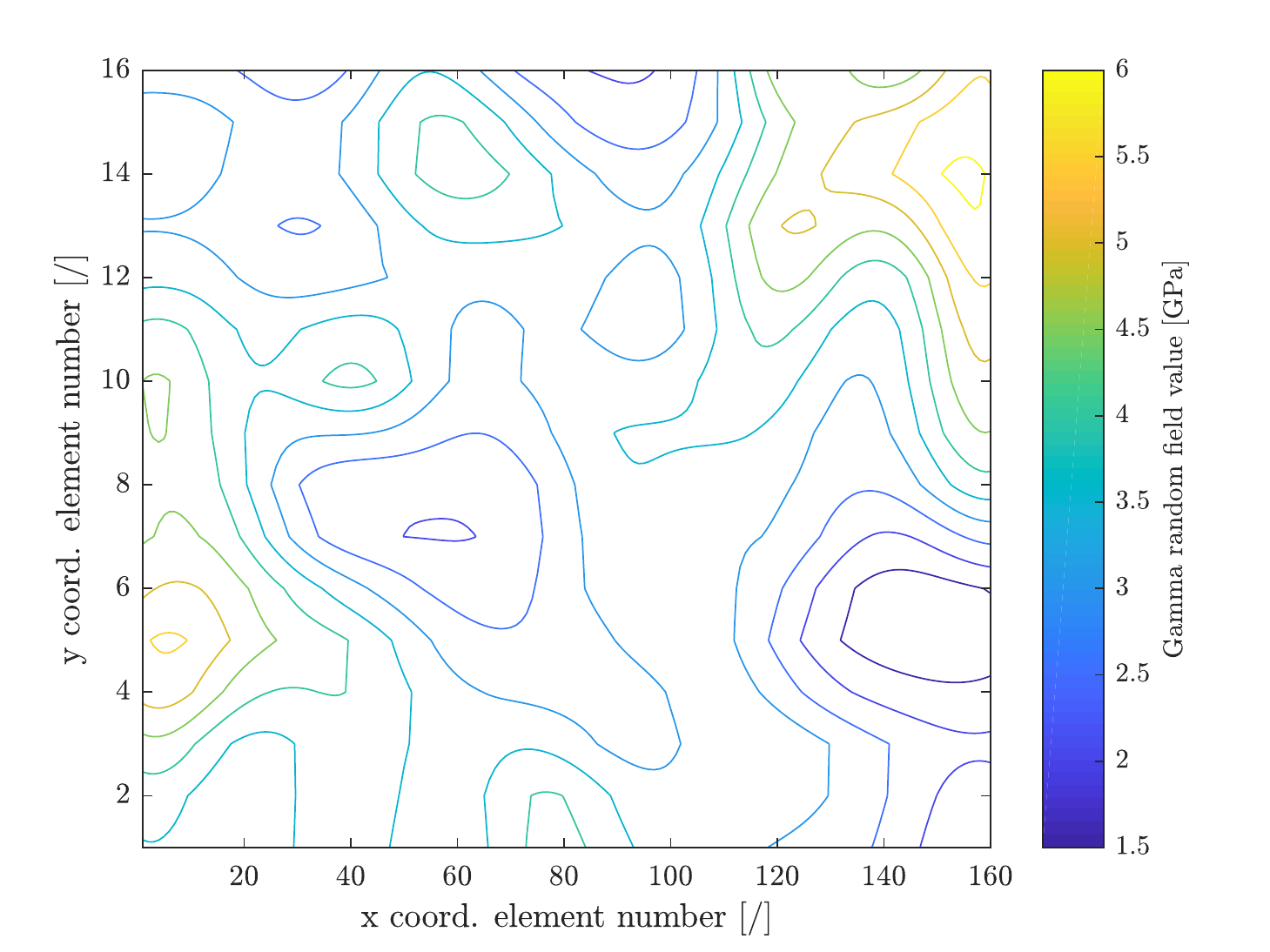}}
\label{fig:Gamma random field}
\end{subfigure}
\caption{Gaussian random field (left) and the corresponding gamma random field (right).}
\label{fig:contour}
\end{figure}

Once the Gaussian field has been generated, a memoryless transformation is applied pointwise,
\begin{linenomath*}
\begin{equation}\label{eq:MemTrans}
g(y) = F^{-1}\left[\Phi(y)\right],
\end{equation}
\end{linenomath*}
in order to obtain the gamma random field \cite{Grigoriu}.
Here, $F$ denotes the marginal cumulative density function (CDF) of the target distribution and $\Phi$ the marginal CDF of the standard normal distribution. This transformation is depicted in Fig.\,\ref{fig:Eigenvalf} (right) with the solid line representing a realization of $F$, and the dashed line representing a realization of $\Phi$.
Contour plots of a realization of a Gaussian random field and the corresponding gamma random field, with variance and mean of the concrete material, are presented for illustration purposes in Fig.\,\ref{fig:contour}.

\subsection{Problem discretization and solution}
The Finite Element method will be used to compute the responses of the beam assuming plane stress. An equidistant, regular rectangular mesh is applied consisting of Lagrange  quadrilateral elements. The underlying equations and solution methods are reviewed hereunder.
\\
\\
For the static elastic case, the system equation is of the form \label{Displacement_eq}
\begin{linenomath*}
\begin{equation}
\mathbf{K} \mathbf{\underline{u}} = \mathbf{\underline{f}},
\end{equation}
\end{linenomath*}
with $\mathbf{K}$ the global stiffness matrix, $\mathbf{\underline{f}}$ the global nodal force vector and $\mathbf{\underline{u}}$ the displacement. The global stiffness matrix and nodal force vector are obtained from the element stiffness matrices $\mathbf{K^e}$ and the element force vectors $\mathbf{\underline{f}^e}$. These are computed numerically by evaluation of the following integrals by means of Gauss quadrature:
\begin{linenomath*}
\begin{equation}\label{eq:K}
\mathbf{K^e} = \int_{\Omega} \mathbf{B}^T \mathbf{D} \mathbf{B} d \Omega
~~~~\mbox{and}~~~
\mathbf{\underline{f}^e} = \int_{\Gamma_t} \mathbf{N}^T \mathbf{\overline{t}}_n d \Gamma_t.
\end{equation}
\end{linenomath*}

The element nodal force vector $\mathbf{\underline{f}^e}$ is modeled as a Neumann boundary condition, where  $\mathbf{\overline{t}}_n$ stands for the surface traction specified as a force per unit area  and $\mathbf{N}$ is the element shape function matrix, integrated over the free element boundary $\Gamma_t$. The element stiffness matrix $\mathbf{K^e}$ is obtained by integrating the matrix  $\mathbf{B^TDB}$ over the element's surface $\Omega$. Matrix $\mathbf{B}$ is defined as $\mathbf{LN}$ with $\mathbf{L}$ the derivative matrix specified below, and  $\mathbf{D}$  is the elastic constitutive matrix for plane stress, containing the element-wise material parameters,
\begin{linenomath*}
\begin{equation}
\mathbf{L} = 
\begin{bmatrix}
\pdv{}{x} & 0  \\
0 & \pdv{}{y}  \\
\pdv{}{y} & \pdv{}{x}  \\
\end{bmatrix}\,
~~~~\mbox{and}~~~
\mathbf{D} = \dfrac{E}{1-\nu^2}
\begin{bmatrix}
1 & \nu & 0 \\
\nu & 1 & 0 \\
0 & 0 & \dfrac{1-\nu}{2} \\
\end{bmatrix}\,.
\label{eq:BDmatrix}
\end{equation}
\end{linenomath*}

For the dynamic  case, the following equation is obtained:
\begin{linenomath*}
\begin{equation}\label{eq:dyn}
\left(\mathbf{K}(1+\imath \,\eta) - (2 \pi f)^2 \mathbf{M}\right)\mathbf{\underline{u}} = \mathbf{\underline{f}}
~~~\mbox{with}~~\mathbf{M}^e= \int_{\Omega} \mathbf{N}^T \rho \mathbf{N} d \Omega\,.
\end{equation}
\end{linenomath*}
Matrix $\mathbf{M}$ denotes the system mass matrix obtained from the assembly of the element mass matrices $\mathbf{M^e}$. $f$ denotes the frequency, $\rho$ the volumetric mass density of the material and $\imath$ the imaginary unit. Hysteretic damping is applied, with $\eta$ the damping loss factor.
\\
\\
The approach for solving the static elastoplastic case differs due to  the nonlinear stress-strain relation in the plastic domain. The plastic region is governed by the von Mises yield criterion with isotropic linear hardening. An incremental load approach is used starting with a force of $0\,\mathrm{N}$. The methods used to solve the elastoplastic problem are based on Chapter 2 $\S$4 and Chapter 7 $\S$3 and $\S$4 of \cite{Borst}. For this case, the system equation takes the following form:
\begin{linenomath*}
\begin{equation}\label{Displacement_eq_plast}
\mathbf{K} \Delta\mathbf{\underline{u}} = \mathbf{\underline{r}},
\end{equation}
\end{linenomath*}
where  $\Delta\mathbf{\underline{u}}$ stands for the resulting displacement increment. The vector $\mathbf{\underline{r}}$ is the residual, 
\begin{linenomath*}
\begin{equation}
\mathbf{\underline{r}}=\mathbf{\underline{f}}+\Delta\mathbf{\underline{f}}-\mathbf{\underline{q}},
\end{equation}
\end{linenomath*}
where  $\mathbf{\underline{f}}$ stands for the sum of the external force increments applied in the previous steps, $\Delta\mathbf{\underline{f}}$  for the applied load increment of the current step and   $\mathbf{\underline{q}}$  for the internal force resulting from the stresses
\begin{linenomath*}
\begin{equation}
\mathbf{\underline{q}}=\int_{\Omega} \mathbf{B}^T \bm{\sigma} d\Omega.
\end{equation}
\end{linenomath*}

First the displacement increment of all the nodes is computed according to Eq.\,\eqref{Displacement_eq_plast}, with an initial system stiffness matrix $\mathbf{K}$ resulting from the assembly of the element stiffness matrix $\mathbf{K^e}$,  computed by means of a Gauss quadrature
\begin{linenomath*}
\begin{equation} \label{K_el_plast}
\mathbf{K^e} = \int_{\Omega} \mathbf{B}^T \mathbf{D}^{ep} \mathbf{B} d \Omega,
\end{equation}
\end{linenomath*}
where $\mathbf{D}^{ep}$ denotes the elastoplastic constitutive matrix. The initial state of $\mathbf{D}^{ep}$ is the elastic constitutive matrix from Eq.\,\eqref{eq:BDmatrix}. Secondly, the strain increment $\Delta \varepsilon$ is computed,
\begin{linenomath*}
\begin{equation}
\Delta \varepsilon = \mathbf{B}\Delta\mathbf{\underline{u}}.
\end{equation}
\end{linenomath*}
Thirdly, the nonlinear stress-strain relationship,
\begin{linenomath*}
\begin{equation}
d\bm{\sigma} = \mathbf{D}^{ep} d\varepsilon,
\end{equation}
\end{linenomath*}
is integrated by means of a backward Euler method. The backward Euler method essentially acts as an elastic predictor-plastic corrector; an initial stress state that is purely elastic is computed and then projected in the direction of the yield surface so as to obtain the plastic stress state. Due to the implicit nature of the integrated stress-strain relation, this equation must be supplemented with the integrated form of the hardening rule and the yield condition. This system of nonlinear equations is then solved with an iterative Newton-Raphson method.
Afterwards, the consistent tangent stiffness matrix is computed \cite{Agusti}. This matrix is then used to compute the updated element stiffness matrix, Eq.\,\eqref{K_el_plast}, resulting in an updated system stiffness matrix $\mathbf{K}$. 
The inner iteration step of solving the stress-strain relation and the updated system stiffness matrix is repeated for each outer iteration step which solves Eq.\,\eqref{Displacement_eq_plast}. The outer step  consists in balancing the internal forces with the external ones as to satisfy the residual, which in our case equals $10^{-4}$ times the load increment. The procedure used is incremental-iterative, relying on the iterative Newton-Raphson method.  This process is repeated for each load increment.

\section{The Multilevel Monte Carlo and Quasi-Monte Carlo method}
\subsection{Mesh refinement hierarchies}
\label{sec:Mesh refinement hierarchies}
The Multilevel Monte Carlo method (MLMC) and Multilevel Quasi-Monte Carlo method (MLQMC) are extensions of the standard Monte Carlo (MC) method, see, e.g.,~\cite{Giles, Giles3,Giles2}. These methods rely on a clever combination of many computationally cheap low resolution samples and a relatively small number of higher resolution, but computationally more expensive samples. MLMC and MLQMC require a predefined hierarchy of meshes in order to  work properly. We propose two different types of hierarchies: a hierarchy of meshes based on an increasing number of finite elements, i.e. h-refinement, and a hierarchy based on increasing the order of the polynomial shape function of the finite elements while retaining the same number of elements, i.e. p-refinement. These hierarchies will be indexed from $0$ to $L$, with $0$ indicating the coarsest approximation and $L$ the finest approximation.

\begin{figure}[H]
\centering
\includegraphics[height=1.9cm]{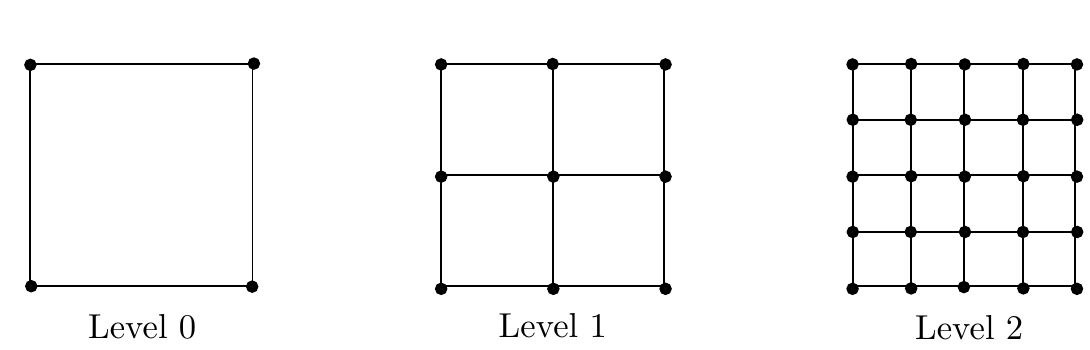}
\caption{Illustrative example of an h-refinement mesh hierarchy used in the MLMC and MLQMC method.}
\label{fig:gridrefinement}
\end{figure}

An example of the first type of mesh hierarchy is shown in Fig.\,\ref{fig:gridrefinement}. As in the Multigrid setting, it is common to use a geometric relation for the number of degrees of freedom between the different levels. We set the number of finite elements for a mesh at level $\ell$  proportional to $2^{d\ell}$, where $d$ is the dimension of the problem $(d=2)$.

\begin{figure}[H]
\centering
\includegraphics[height=1.9cm]{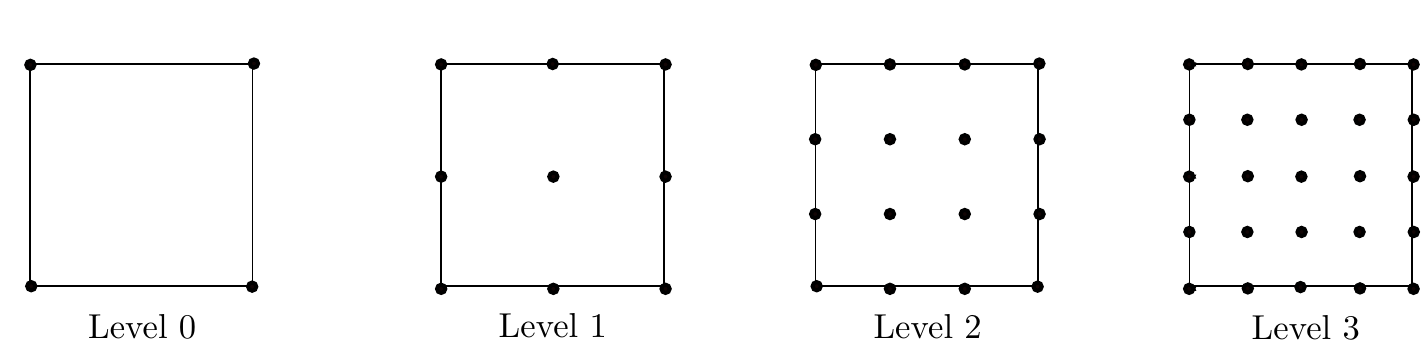}
\caption{Illustrative example of a p-refinement mesh hierarchy used in the MLMC and MLQMC method.}
\label{fig:gridhigherorder}
\end{figure}

An example of a p-refinement hierarchy is shown in Fig.\,\ref{fig:gridhigherorder}. In this work, the higher order elements are defined to be linear, quadratic, cubic, quartic and quintic Lagrangian quadrilateral elements.

We will use both h- and p-hierarchies as levels for MLMC and MLQMC.

\subsection{Multilevel Monte Carlo}
Let $\E[P_L(\omega)]$, or $\E[P_L]$ for short, be the expected value of a particular quantity of interest $P$ depending on a random variable $\omega$, discretized on mesh $L$. The standard MC estimator for $\E[P_L]$ using $N_L$ samples on mesh $L$, denoted as $Q^{\textrm{MC}}_L$, can be written as
\begin{linenomath*}
\begin{equation}
Q^{\textrm{MC}}_L={\frac{1}{N_L}}\sum_{n=1}^{N_L} P_L(\omega^n)\,.
\label{eq:SimpleMC}
\end{equation}
\end{linenomath*}
Multilevel Monte Carlo, on the other hand, starts from a reformulation of $\E[P_L]$ as a telescoping sum. The expected value of the quantity of interest on the finest mesh is expressed as the expected value of the quantity of interest on the coarsest mesh, plus a series of correction terms (or \emph{differences}):
\begin{linenomath*}
\begin{equation}
\E[P_L]=\E[P_0]+\sum_{\ell=1}^L \E[P_\ell -P_{\ell-1}]\,.
\label{Eq:TelescopingSum}
\end{equation}
\end{linenomath*}
Each term in the right-hand side is then estimated separately by a standard Monte Carlo estimator with $N_\ell$ samples, i.e.,
\begin{linenomath*}
\begin{equation}
Q^{\textrm{MLMC}}_L=\frac{1}{N_0}\sum_{n=1}^{N_0} P_0(\omega^{n}) + \sum_{\ell=1}^L \left \{ \frac{1}{N_\ell} \sum_{n=1}^{N_\ell} \left( P_\ell(\omega^{n})-P_{\ell-1}(\omega^{n})\right) \right \},
\label{eq:MLMC}
\end{equation}
\end{linenomath*}
where $Q^{{\textrm{MLMC}}}_L$ is the Multilevel Monte Carlo estimator for the expected value $\E[P_L]$, which is a discrete approximation for the expected value of the quantity of interest, $\E[P]$. The mean square error (MSE) is defined as 
\begin{linenomath*}
\begin{equation}\label{eq:MSE}
\begin{split}
\textrm{MSE}(Q^{\textrm{MLMC}}_L) & := \E\left[\left(Q^{\textrm{MLMC}}_L - \E\left[P\right]\right)^2\right] \\
& := \V\left[Q^{\textrm{MLMC}}_L\right] + \left(\E\left[Q^{\textrm{MLMC}}_L\right] - \E\left[P\right]\right)^2,
\end{split}
\end{equation}
\end{linenomath*}
with $\V\left[\mathrm{\cdot}\right]$ denoting the variance of a random variable $\mathrm{\cdot}$.
The MLMC estimator in Eq.\,\eqref{eq:MLMC} can be written as a sum of $L+1$ estimators for the expected value of the difference on each level, i.e.,
\begin{linenomath*}
\begin{equation}
Q^{\textrm{MLMC}}_L = \sum_{\ell = 0}^L Y_{\ell}, \quad \text{where} \quad Y_{\ell} = \frac{1}{N_\ell} \sum_{n=1}^{N_\ell} \left( P_\ell(\omega^{n})-P_{\ell-1}(\omega^{n})\right).
\end{equation}
\end{linenomath*}
where we defined $P_{-1}\coloneqq0$.

Because of the telescoping property, the MLMC estimator is an unbiased estimator for the quantity of interest on the finest mesh, i.e.,
\begin{linenomath*}
\begin{equation}
\E[Q^{\textrm{MLMC}}_L] = \E[P_L].
\end{equation}
\end{linenomath*}

Denoting by $V_{\ell}$ the variance of the difference, $V_{\ell} = \V(P_{\ell} - P_{\ell-1})$, the variance of the estimator can be written as
\begin{linenomath*}
\begin{equation}
\V[Q^{\textrm{MLMC}}_L] = \sum_{\ell=0}^L \frac{V_{\ell}}{N_{\ell}}.
\label{eq:var}
\end{equation}
\end{linenomath*}

In order to ensure that the MSE in Eq.\,\eqref{eq:MSE} is below a given tolerance $\epsilon^2$, it is sufficient to enforce that the variance $\V[Q^{\textrm{MLMC}}_L]$, given in Eq.\,\eqref{eq:var}, and the squared bias $(\E[P_L-P])^2$ are both less than $\epsilon^2/2$.
The condition on the variance of the estimator can be used to determine the number of samples needed on each level $\ell$. Following the classic argument by Giles in~\cite{Giles}, we minimize the total cost of the MLMC estimator
\begin{linenomath*}
\begin{equation}
\text{cost}(Q^{\textrm{MLMC}}) = \sum_{\ell=0}^{L} N_\ell C_\ell, 
\label{eq:cost_sample}
\end{equation}
\end{linenomath*}
where $C_\ell$ denotes the cost to compute a single realization of the difference $P_\ell-P_{\ell-1}$, subject to the constraint
\begin{linenomath*}
\begin{equation}
\sum_{\ell=0}^L \frac{V_{\ell}}{N_{\ell}} \leq \frac{\epsilon^2}{2}.
\end{equation}
\end{linenomath*}
Treating the $N_\ell$ as continuous variables, we find
\begin{linenomath*}
\begin{equation}\label{eq:nopt}
 N_\ell = \frac{2}{\epsilon^2} \sqrt{\frac{V_\ell}{C_\ell}} \sum_{\ell=0}^L \sqrt{V_\ell C_\ell} .
\end{equation}
\end{linenomath*}
Note that if $\E[P_\ell]\rightarrow\E[P]$, then $V_\ell\rightarrow0$ as $\ell$ increases. Hence, the number of samples $N_\ell$ will be a decreasing function of $\ell$. This means that most samples will be taken on the coarse mesh, where samples are cheap, whereas increasingly fewer samples are required on the finer, but more expensive meshes. In practice, the number of samples must be truncated to $\lceil N_\ell \rceil$, the least integer larger than or equal to $N_\ell$. 

Using Eq.\,\eqref{eq:nopt}, the total cost of the MLMC estimator, from Eq.\,\eqref{eq:cost_sample}, can be written as
\begin{linenomath*}
\begin{equation}\label{eq:cost}
\text{cost}(Q^{\textrm{MLMC}}) = \frac{2}{\epsilon^2}\left(\sum_{\ell=0}^L \sqrt{V_\ell C_\ell}\right)^{2}.
\end{equation}
\end{linenomath*}
This can be interpreted as follows. When the variance $V_\ell$ decreases faster with increasing level $\ell$ than the cost increases, 
the dominant computational cost is located on the coarsest level. The computational cost is then proportional to $V_0 C_0$, which is small because $C_0$ is small. Conversely, if the variance decreases slower with increasing level $\ell$ than the cost increases, 
the dominant computational cost will be located on the finest level $L$, and proportional to $V_L C_L$. This quantity is small because $V_L$ is small. For comparison, the computational cost of a Monte Carlo simulation that reaches the same accuracy is proportional to $V_0 C_L$.

The second term in Eq.\,\eqref{eq:MSE} is used to determine the maximum number of levels $L$. A typical MLMC implementation is level-adaptive, i.e., starting from a coarse finite element mesh, finer meshes are only added if required to reach a certain accuracy. Assume that the convergence $\E[P_\ell]\rightarrow\E[P]$ is bounded as $|\E[P_{\ell} - P]| = \mathcal{O}(2^{-\alpha \ell})$. Then we can use the heuristic
\begin{linenomath*}
\begin{equation}\label{eq:bias_constraint}
\abs{\E[P_L-P]} = \abs{\sum_{\ell=L+1}^{\infty} \E[P_\ell-P_{\ell-1}]} \approx \frac{\abs{\E[P_L-P_{L-1}]}}{2^\alpha-1}
\end{equation}
\end{linenomath*}
and check for convergence using $|\E[P_L-P_{L-1}]|/(2^\alpha-1)\leq\epsilon/\sqrt{2}$, see~\cite{Giles} for details.

\subsection{Multilevel Quasi-Monte Carlo}
\label{Sec:MLQMC}
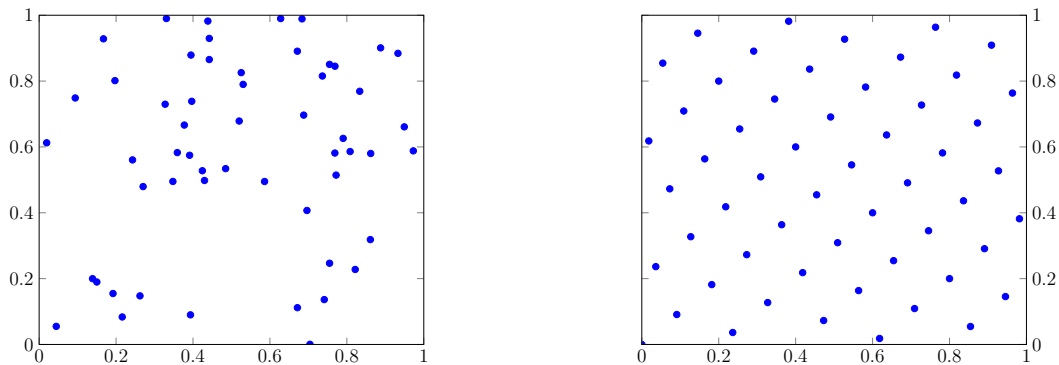
\begin{figure}[H]
\begin{subfigure}[b]{0.54\linewidth}
\centering
\scalebox{0.6}{
%
\begin{tikzpicture}
\begin{axis}[%
scale only axis,
xmin=0,
xmax=1,
ymin=0,
ymax=1,
every x tick label/.append style={font=\large},
every x label/.append style={font=\large},
every y tick label/.append style={font=\large},
every y label/.append style={font=\large},
axis background/.style={fill=white},
clip marker paths=true, axis on top=true]
\addplot [color=blue,only marks,mark=*,mark options={solid,fill=blue},forget plot]
  table[row sep=crcr]
  {
0.627896379614169	0.989872153631504 \\
0.771980385554245	0.514423456505704\\
0.93285357027882	0.884281023126955\\
0.972740854003014	0.588026055308497\\
0.192028349427775	0.154752348656045\\
0.138874202829155	0.199862822857452\\
0.696266337082995	0.406954837138907\\
0.0938200267748656	0.748705718215691\\
0.525404403859336	0.825583815786156\\
0.530344218392863	0.789963029944531\\
0.861139811393332	0.318524245398992\\
0.484853333552102	0.534064127370726\\
0.393456361215266	0.089950678770581\\
0.671431139674026	0.111705744193203\\
0.741257943454206	0.136292548938299\\
0.520052467390387	0.678652304800188\\
0.347712671277525	0.495177019089661\\
0.149997253831683	0.18971040601758\\
0.586092067231462	0.495005824990221\\
0.262145317727807	0.147608221976689\\
0.0444540922782385	0.0549741469061882\\
0.754933267231179	0.850712674289007\\
0.242785357820962	0.560559527354885\\
0.442402313001943	0.929608866756663\\
0.687796085120107	0.696667200555228\\
0.359228210401861	0.58279096517584\\
0.736340074301202	0.815397211477421\\
0.394707475278763	0.879013904597178\\
0.683415866967978	0.988911616079589\\
0.704047430334266	0.000522375356944771\\
0.442305413383371	0.865438591013025\\
0.0195776235533187	0.612566469483999\\
0.330857880214071	0.989950205708831\\
0.424309496833137	0.527680069338442\\
0.270270423432065	0.479523385210219\\
0.197053798095456	0.801347605521952\\
0.82172118496131	0.227842935706042\\
0.429921409383266	0.49809429119639\\
0.887770954256354	0.900852488532005\\
0.391182995461163	0.574661219130188\\
0.769114387388296	0.845178185054037\\
0.396791517013617	0.738640291995402\\
0.808514095887345	0.585987035826476\\
0.755077099007084	0.246734525985975\\
0.377395544835103	0.666416217319468\\
0.216018915961394	0.0834828136026227\\
0.790407217966913	0.625959785171583\\
0.949303911849797	0.660944557947342\\
0.327565434075205	0.729751855317221\\
0.67126437045174	0.890752116325322\\
0.438644982586956	0.982303222883606\\
0.833500595588975	0.769029085335896\\
0.768854252429615	0.581446487875398\\
0.167253545494722	0.928313062314188\\
0.861980478702072	0.580090365758442\\
};
\end{axis}
\end{tikzpicture}%
}
\end{subfigure}
\begin{subfigure}[b]{0.45\linewidth}
\centering
\scalebox{0.6}{
%
\begin{tikzpicture}

\begin{axis}[%
scale only axis,
xmin=0,
xmax=1,
ymin=0,
ymax=1,
every x tick label/.append style={font=\large},
every x label/.append style={font=\large},
every y tick label/.append style={font=\large},
every y label/.append style={font=\large},
axis background/.style={fill=white},
clip marker paths=true, axis on top=true,
 ylabel near ticks,yticklabel pos=right,
]
\addplot [color=blue,only marks,mark=*,mark options={solid,fill=blue},forget plot]
  table[row sep=crcr]
  {%
0	0\\
0.0181818181818182	0.618181818181818\\
0.0363636363636364	0.236363636363636\\
0.0545454545454545	0.854545454545454\\
0.0727272727272727	0.472727272727273\\
0.0909090909090909	0.0909090909090908\\
0.109090909090909	0.709090909090909\\
0.127272727272727	0.327272727272727\\
0.145454545454545	0.945454545454545\\
0.163636363636364	0.563636363636363\\
0.181818181818182	0.181818181818182\\
0.2	0.8\\
0.218181818181818	0.418181818181818\\
0.236363636363636	0.036363636363637\\
0.254545454545455	0.654545454545454\\
0.272727272727273	0.272727272727273\\
0.290909090909091	0.890909090909091\\
0.309090909090909	0.50909090909091\\
0.327272727272727	0.127272727272727\\
0.345454545454545	0.745454545454546\\
0.363636363636364	0.363636363636363\\
0.381818181818182	0.981818181818182\\
0.4	0.6 \\
0.418181818181818	0.218181818181819\\
0.436363636363636	0.836363636363636\\
0.454545454545455	0.454545454545455\\
0.472727272727273	0.0727272727272741\\
0.490909090909091	0.690909090909091\\
0.509090909090909	0.309090909090909\\
0.527272727272727	0.927272727272726\\
0.545454545454545	0.545454545454547\\
0.563636363636364	0.163636363636364\\
0.581818181818182	0.781818181818181\\
0.6	0.399999999999999\\
0.618181818181818	0.0181818181818194\\
0.636363636363636	0.636363636363637\\
0.654545454545455	0.254545454545454\\
0.672727272727273	0.872727272727271\\
0.690909090909091	0.490909090909092\\
0.709090909090909	0.109090909090909\\
0.727272727272727	0.727272727272727\\
0.745454545454545	0.345454545454544\\
0.763636363636364	0.963636363636365\\
0.781818181818182	0.581818181818182\\
0.8	0.199999999999999\\
0.818181818181818	0.818181818181817\\
0.836363636363636	0.436363636363637\\
0.854545454545454	0.0545454545454547\\
0.872727272727273	0.672727272727272\\
0.890909090909091	0.290909090909089\\
0.909090909090909	0.90909090909091\\
0.927272727272727	0.527272727272727\\
0.945454545454545	0.145454545454548\\
0.963636363636364	0.763636363636365\\
0.981818181818182	0.381818181818183\\
};
\end{axis}
\end{tikzpicture}%
}
\end{subfigure}
\caption{Example of points sampled for MLMC (left) and MLQMC (right).}
\label{fig:points}
\end{figure}

One of the major differences with MLMC is that for MLQMC, the individual sample points are not chosen at random but according to a deterministic rule, see for example Fig.\,\ref{fig:points}.

In this paper, we use rank-1 lattice rules similar to \cite{Giles3}. These points have the following representation: the n-th sample point $\mathbf{x}_n$ is defined as
\begin{linenomath*}
\begin{equation}
\mathbf{x}_n = \fracfun{\frac{n}{N}\mathbf{z}},
\label{eq:QMC_points_no_shift}
\end{equation}
\end{linenomath*}
where $\fracfun{x} = x - \floor{x}, x>0$. Vector $\mathbf{z}$ is an $s$-dimensional vector of positive integers, and $N$ is the number of points in the lattice rule.

Due to the deterministic nature of the MLQMC points, a shift has to be introduced in order to  obtain unbiased estimates of the quantities of interest, as discussed in section 2.9 of \cite{Sloan}. Eq.\,\eqref{eq:QMC_points_no_shift} is rewritten as 
\begin{linenomath*}
\begin{equation}
\mathbf{x}_{i,n} = \fracfun{\frac{n}{N}\mathbf{z} + \Delta_i},
\end{equation}
\end{linenomath*}
where $\Delta$ is a shift or offset, uniformly distributed in $\left[0,1\right]^s$. In practice, multiple random shifts must be chosen, labeled $\Delta_1, \Delta_2, ..., \Delta_R$, in order to allow for the computation of the variance of the estimator, and hence the MSE. The MLQMC estimator is then written as 
\begin{linenomath*}
\begin{equation}
Q^{\textrm{MLQMC}}_L= \frac{1}{R_0}\sum_{i=1}^{R_0}\frac{1}{N_0}\sum_{n=1}^{N_0} P_0(\mathbf{x}_{i,n}) + 
 \sum_{\ell=1}^L \frac{1}{R_\ell}\sum_{i=1}^{R_\ell}\left \{ \frac{1}{N_\ell} \sum_{n=1}^{N_\ell} \left( P_\ell(\mathbf{x}_{i,n})-P_{\ell-1}(\mathbf{x}_{i,n})\right) \right \}.
\end{equation}
\end{linenomath*}
We choose the number of shifts to be constant on each level, i.e., $R_\ell=R$, $\ell=0,1,\ldots,L$. A value $R=10$ will be chosen in our numerical experiments. Contrary to MLMC, the number of samples for MLQMC is not the result of an optimization problem, as in Eq.\,\eqref{eq:nopt}. For MLQMC an adaptive algorithm is used, see \cite{Giles3}. Starting with an initial number of samples, this algorithm multiplies the number of samples on the level with maximum ratio $V_\ell/C_\ell,$ with a constant factor until the variance of the estimator is smaller than $\frac{\epsilon^2}{2}$, where $V_\ell$ is defined as $\V(P_{\ell} - P_{\ell-1})$. In our implementation this multiplication constant is chosen as $1.2$. 

The MLQMC method is expected to work particularly well if the number of subsequent uncertainties decays rapidly, see \cite{KUO,Kuo2016}, as is the case with a smooth random field generated according to a KL expansion, where the magnitude of the successive eigenvalues is decaying rapidly. For more information, see \cite{SLOAN19981}. 
\subsection{Cost Theorem}

Having introduced both methods, we  now present a complexity theorem for MLQMC, which also covers the MLMC method, when $\delta = 1$, see Theorem 1. More details can be found in  \cite{KUO} and on page 76 of \cite{Teckentrup}.
\begin{theorem}
Given the  positive constants $\alpha, \beta, \gamma, c_1, c_2, c_3$ such that $\alpha \geq \dfrac{1}{2} \mathrm{min}\left(\beta,\delta^{-1}\gamma\right)$ with $\delta \in \left(1/2,1\right]$ and assume that the following conditions hold: 
\begin{enumerate}
\item $\lvert \E[P_{\ell} - P] \rvert \leq c_1 2^{-\alpha \ell}$,
\item $\V\left[Y_{\ell}\right] \leq  c_2 2^{-\beta \ell} N_{\ell}^{-1/\delta} $ \;and
\item $C_{\ell} \leq  c_3 2^{\gamma \ell}$.
\end{enumerate}

Then, there exists a positive constant $c_4$ such that for any $\epsilon < \exp(-1)$ there exists an  $L$ and  a sequence $\{N_{\ell}\}_{\ell=0}^L$ for which the multilevel estimator, $Q^{\mathrm{MLQMC}}_L$ has an $\mathrm{MSE}\leq\epsilon^2$, and
\begin{equation}
\scalebox{1}{$
\mathrm{cost}(Q^{\mathrm{MLQMC}}) \leq    \left\{
  \begin{aligned}
& c_4 \epsilon^{-2 \delta} && \mathrm{if} \quad  \delta\beta > \gamma, \\
& c_4 \epsilon^{-2 \delta}\left(\log \;\epsilon \right)^{1+\delta} && \mathrm{if} \quad  \delta\beta = \gamma, \\
& c_4 \epsilon^{-2\delta-\left(\gamma-\delta\beta\right)/\alpha} && \mathrm{if} \quad  \delta\beta < \gamma. \\
  \end{aligned}
  \right.$}
  \label{eq:Algo_regime}
\end{equation}
\label{Theorem_1}
\end{theorem} 

The factor $\alpha$, in assumption 1,  is the rate at which the expected value of the differences decreases with increasing  level. $\beta$, in assumption 2,  stands for the decay rate of the variance of the differences. The factor $\gamma$, in assumption 3, is determined by the efficiency of the solver. This factor will be different for the h-refinement scheme and the p-refinement scheme. All three factors will be estimated on the fly in our numerical experiments.

Following this theorem, the optimal cost of the MLMC estimator, is proportional to $\epsilon^{-2}$ when the variance over the levels decreases faster than the cost per level increases, i.e., $\beta > \gamma$, and $\delta = 1$. Similarly, for the MLQMC estimator, the optimal cost is proportional to $\epsilon^{-1}$. Note that this is only true in the limit, i.e., $\delta \to 1/2$. We will show in our numerical experiments that the theoretically derived asymptotic cost complexity is close to what we observe.

\subsection{Implementation details}\label{Implementation_aspects}

The MLMC and MLQMC methods are non-intrusive, requiring only an interface between the Finite Element solver routine and the multilevel routine. The Finite Element solver routines are written in \textsc{Matlab}, while the multilevel routine is written in  \textsc{Julia} \cite{PJ}.

All the computations are run in parallel, the computation of  the individual samples is parallelized. This is possible because of the \emph{embarrassingly parallel} nature of all the  Monte Carlo methods. In the aforementioned configuration, a number of $28$ samples can be computed concurrently. For more details on load balancing of MLMC/MLQMC samplers, we refer to \cite{Scheichl}.

When dealing with multiple quantities of interest (Qoi), the optimal amount of samples in Eq.\,\eqref{eq:nopt} is evaluated with the variance $V_\ell$ corresponding to the variance of the Qoi with the largest variance. By doing so, the variance constraint is guaranteed to be satisfied for all other Qoi's. Here, we consider only one Qoi. For the static cases, the Qoi  is the largest transversal deflection, which in this case is characterized by the largest variance. For the dynamic case, the Qoi is the response at the node which has the largest frequency response, which also is the node with largest variance.

For MLMC, the computation of the optimal number of samples per level according to Eq.\,\eqref{eq:nopt} is based on the variances of one degree of freedom (dof) on these levels. The selected dof is the one with the largest response variance. This ensures that the variance constraint is satisfied for all other dof's. For MLQMC, an adaptive algorithm is used to determine the number of samples needed, see \S\ref{Sec:MLQMC}. A first estimation of these variances is done by computing a trial sample set on levels 0, 1 and 2.
For MLMC, the size of this sample set is $40$ for both the elastic and the elastoplastic cases. Variances on additional levels are estimated according to the second condition from Theorem \ref{Theorem_1}, following \cite{Giles, Giles2}. For MLQMC, 2 samples with 10 shifts are taken not only on the initial levels (0, 1, 2) but also on all additional levels. 

\begin{figure}[H]
\begin{subfigure}[b]{0.32\textwidth}
\scalebox{0.32}{
	\makebox[\textwidth]{
	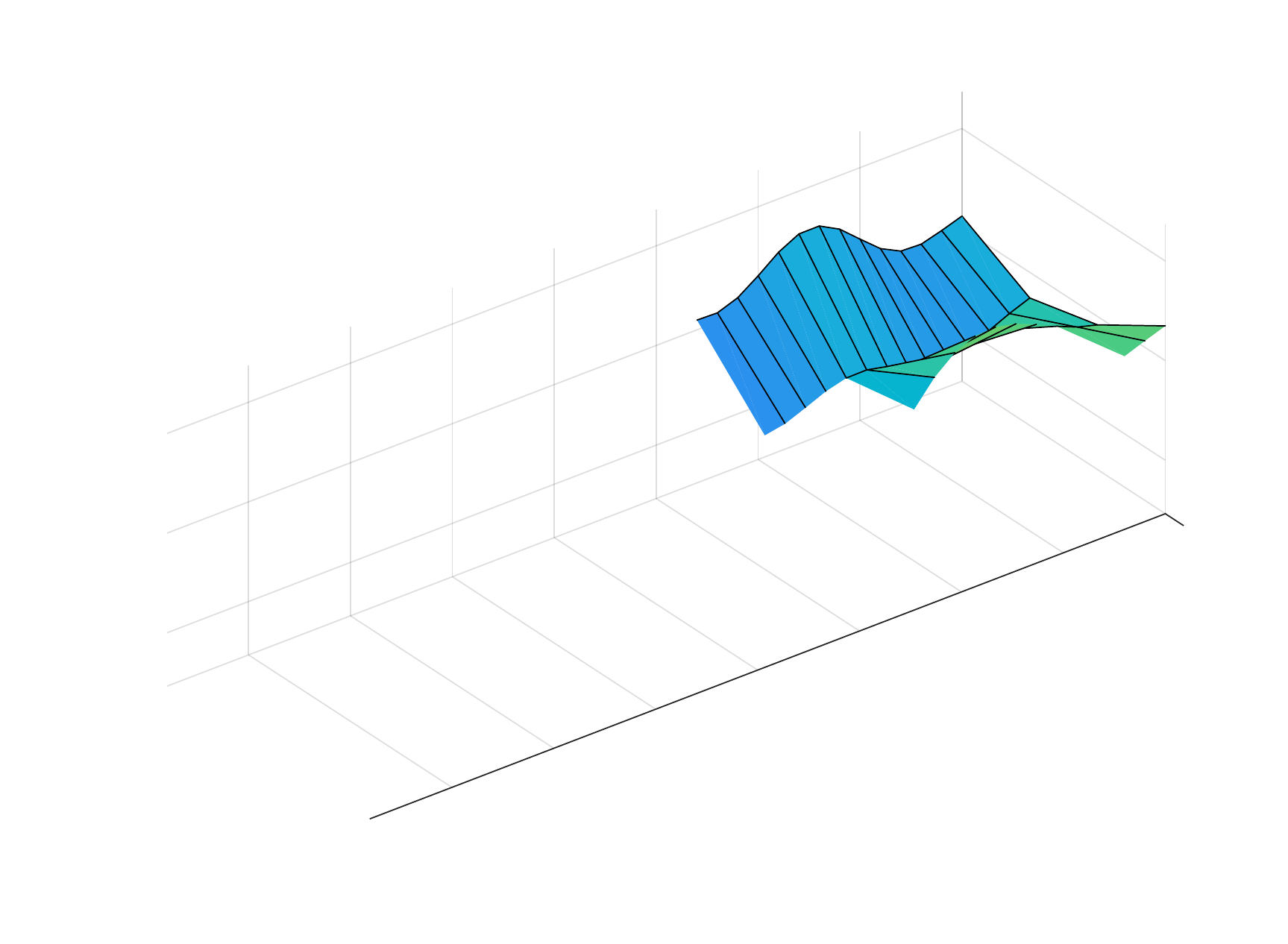}}
\end{subfigure}
\begin{subfigure}[b]{0.32\linewidth}
\centering
	\scalebox{0.32}{
		\makebox[\textwidth]{
	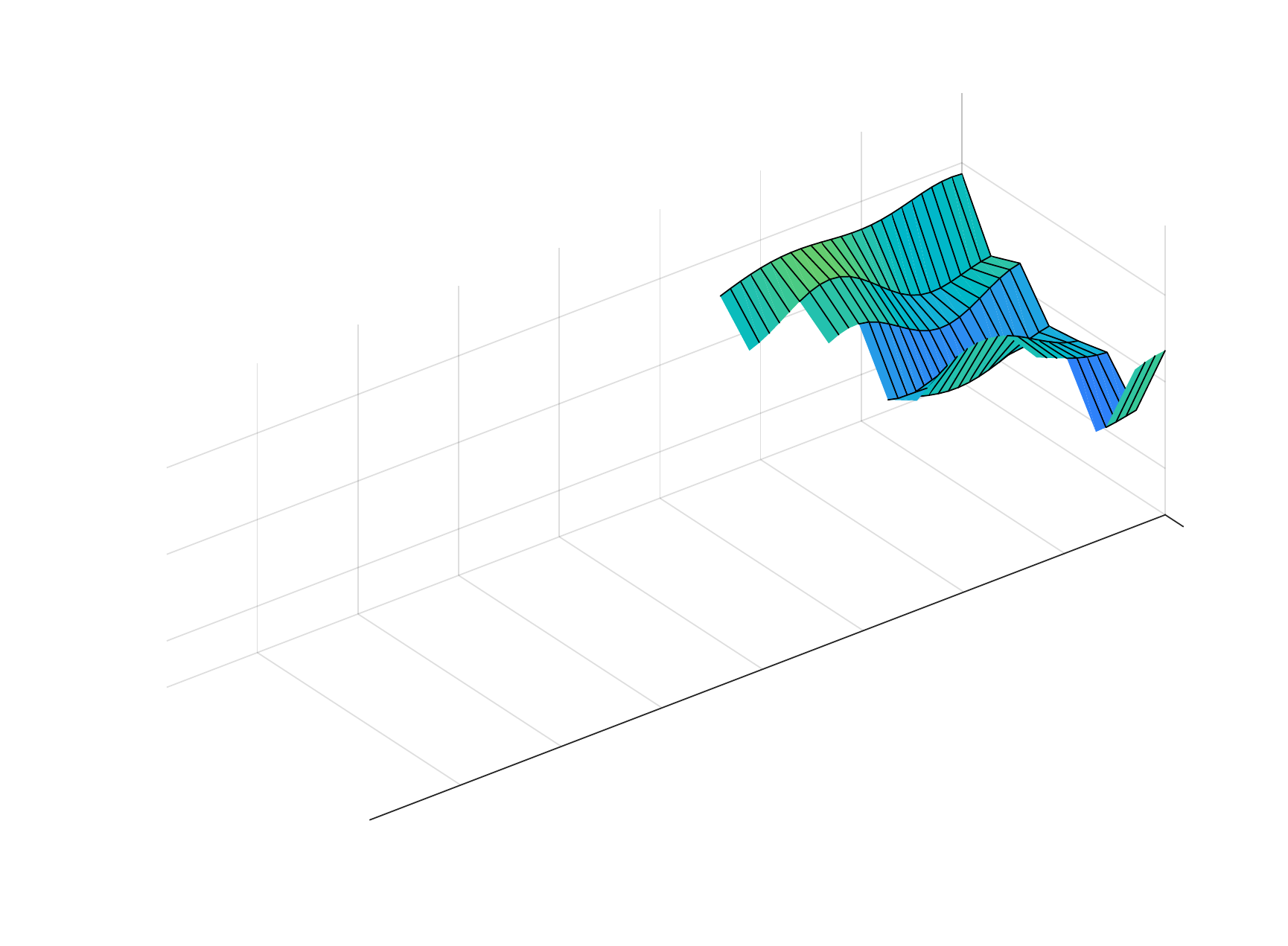}}
\end{subfigure}
\begin{subfigure}[b]{0.32\linewidth}
\hspace{2.5cm}
	\scalebox{0.32}{
			\makebox[\textwidth]{
	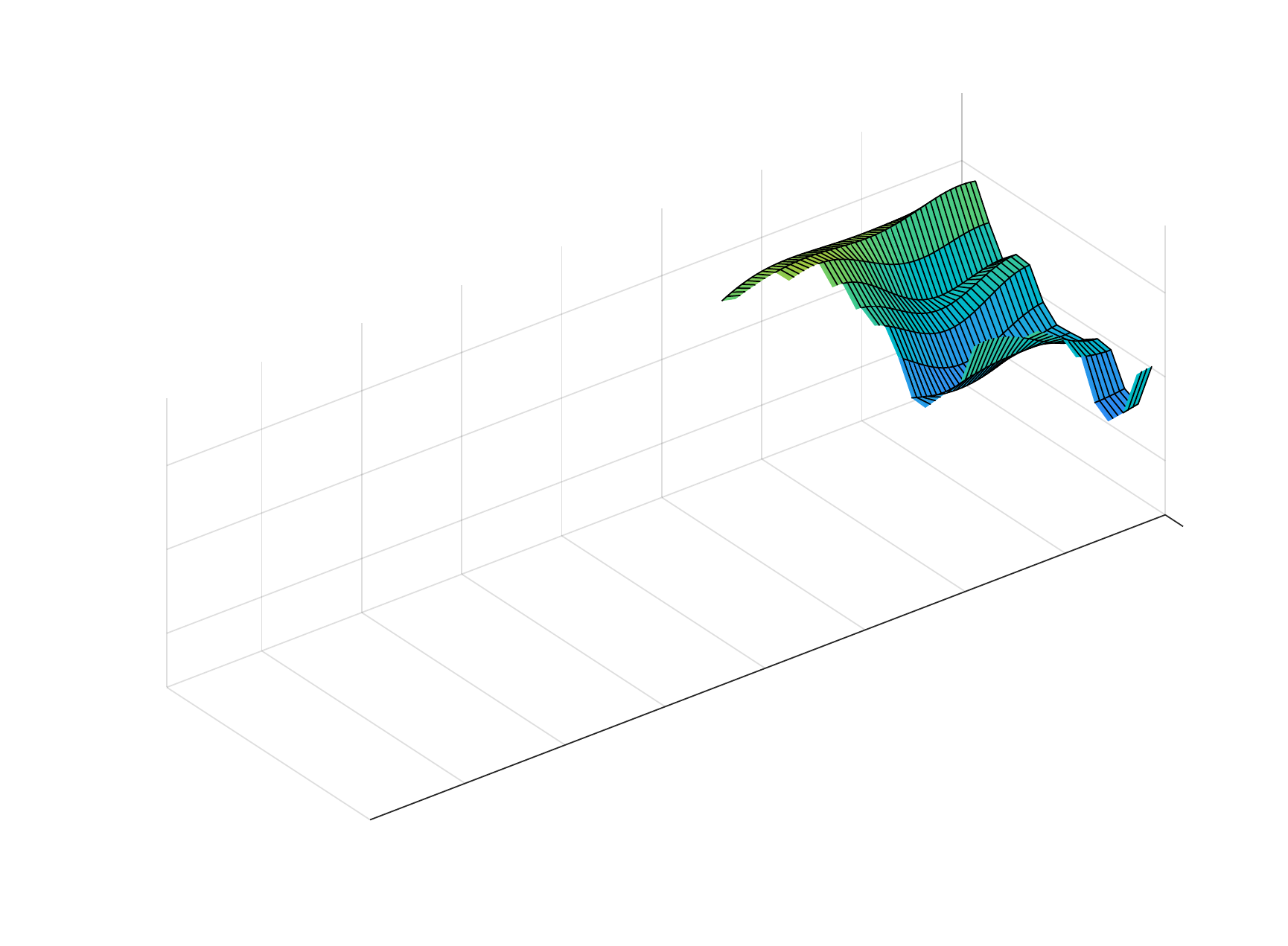}}
\end{subfigure}
	\caption{Realizations of a Gaussian random field on level 0 (left), level 1 (middle), level 2 (right). }
	\label{fig:level_overview}
\end{figure}

Each finite element is assigned a value of the Young's modulus. For the h-refinement mesh hierarchy, this is accomplished by means of the midpoint approach, i.e., the value is taken constant within each individual element and equal to the value of the realization of the random field at the center point of the element \cite{jie_lie}. For the p-refinement mesh hierarchy, this is accomplished by means of the integration point method, i.e., the Young's modulus is computed at the Gauss integration points when numerically computing the element stiffness matrix \cite{Brenner}. An illustrative example of a Gaussian random field for three successive h-levels can be seen in Fig.\,\ref{fig:level_overview}. 

\section{Numerical Results}

In this section, we discuss our numerical experiments with the MLMC and MLQMC method. We consider the static and dynamic cases, using both a homogeneous and a heterogeneous uncertain Young's modulus. First, we introduce the different simulated cases. For the static cases, the solution consists of the displacement of the beam in the spatial domain. For the dynamic case, the solution is a frequency response.  Secondly, we illustrate in which way the uncertainty model affects the uncertainty bounds of the solution for the static and dynamic cases. It should be noted that the static cases  require only one MLMC/MLQMC simulation and one MC simulation for comparison. In contrast, solutions for the dynamic case require multiple individual MLMC/MLQMC simulations and multiple MC simulations: one for each individual frequency of the frequency response function.  Thirdly, we compare MLMC, MLQMC and standard MC combined with h- and p-refinement for the static cases. We estimate the rates ($\alpha$, $\beta$, $\gamma$) from Theorem 1, present the number of samples for the finest considered tolerance and compare the different methods in terms of total simulation time.  

All the results have been computed  on a  workstation equipped with 14 physical (28 logical) cores, Intel Xeon E5645 CPU's, clocked at 2.40 GHz,  and a total of 128 GB RAM.  

\subsection{Presentation of the simulated cases}
\label{par:comparison}
The simulated cases have been presented in Tab.\,\ref{Tab:cases}. In order to model the uncertainty which is present in the Young's modulus of the beam, we opt for two different ways: a heterogeneous and a homogeneous Young's modulus. The load is modeled as a distributed load acting on each of the vertical middle nodes of the beam. The sum of all these individual loads is independent of the refinement of the mesh.  For the elastic case, the total load equals $10000\,\mathrm{kN}$. The load for the elastoplastic case is an incremental load starting at $0\,\mathrm{N}$ until $13.5\,\mathrm{kN}$ in steps of $135\,\mathrm{N}$. 

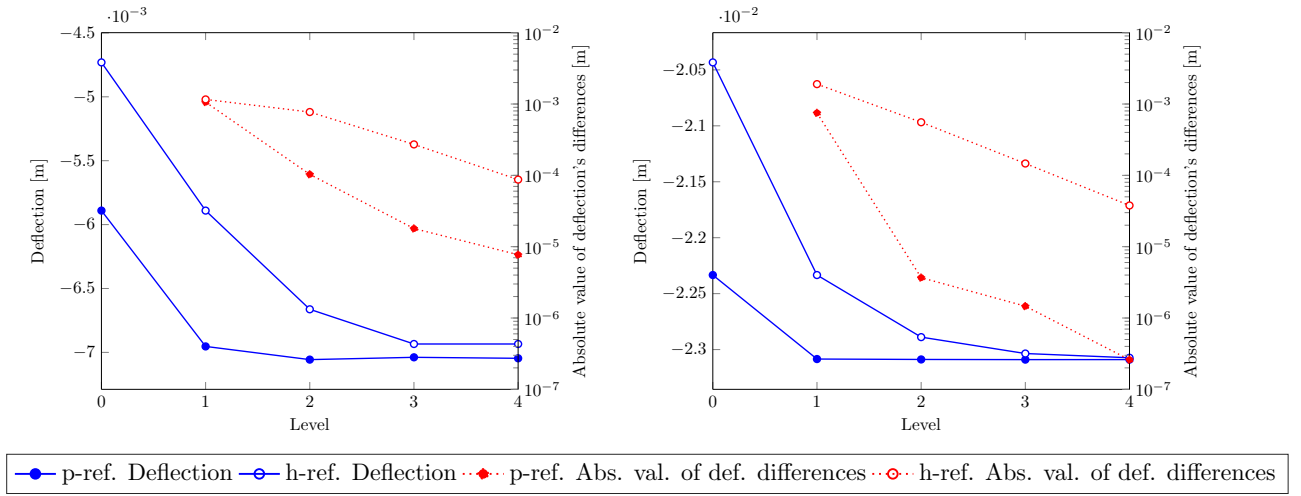
\begin{figure}[H]
\centering
        \begin{subfigure}[b]{0.475\textwidth}
        \scalebox{0.65}{
\begin{tikzpicture}
\pgfplotsset{
    scale only axis,
    xmin=0, xmax=4,
     xtick distance=1,
}

\begin{axis}[
  axis y line*=left,
  xlabel=Level,
  ylabel={Deflection [m]},
  xmin=0,
  xmax=4,
 xtick distance=1,every axis plot/.append style={thick}
]
\addplot[mark=*,blue]
  coordinates{
    (0,-0.00589041465874002)
    (1,-0.00695336058172622)
    (2,-0.00705660658168465)
    (3,-0.00703863951306348)
    (4,-0.00704635014844012)
}; \label{def1}

\addplot[mark=*,blue,mark options={solid,fill =white}]
  coordinates{
    (0,-0.00473074046750018)
    (1,-0.00589041465874002)
    (2,-0.00666283713193700)
    (3,-0.00693432675026149)
    (4,-0.00693432675026149)
}; \label{def2}

\end{axis}

\begin{axis}[
  axis y line*=right,
  axis x line=none,
  ymin=10^-7, ymax=10^-2,
  xmin=0,
  xmax=4,
ymode={log},
  ylabel={Absolute value of deflection's differences [m]},
  xtick distance=1,every axis plot/.append style={thick}
]

\addplot[mark=*,red,style={dotted}]
  coordinates{
    (1,0.00106294592298620)
    (2,0.000103245999958426)
    (3,1.79670686211653e-05)
    (4,7.71063537664037e-06)
};

\addplot[mark=*,red,style={dotted},mark options={solid,fill =white}]
  coordinates{
    (1,0.00115967419123984)
    (2,0.000772422473196979)
    (3,0.000271489618324491)
    (4,8.72170401129491e-05)
};

\end{axis}

\end{tikzpicture}}
 \end{subfigure}
         \begin{subfigure}[b]{0.475\textwidth}
        \scalebox{0.65}{
\begin{tikzpicture}
\pgfplotsset{
    scale only axis,
    xmin=0, xmax=4,
         xtick distance=1,
}

\begin{axis}[
  axis y line*=left,
  xlabel=Level,
  ylabel={Deflection [m$]$},every axis plot/.append style={thick}
]
\addplot[mark=*,blue]
  coordinates{
    (0,-0.0223337341764705)
    (1,-0.0230846367423070)
    (2,-0.0230883163806814)
    (3,-0.0230897813417939)
    (4,-0.0230895211165382)
}; \label{def1}
\addplot[mark=*,blue,mark options={solid,fill =white}]
  coordinates{
    (0,-0.0204330927443515)
    (1,-0.0223337341764705)
    (2,-0.0228887254988895)
    (3,-0.0230353318915608)
    (4,-0.0230731684871195)
}; \label{def2}
\end{axis}

\begin{axis}[
  axis y line*=right,
  axis x line=none,
  ymin=10^-7, ymax=10^-2,
ymode={log},
  ylabel={Absolute value of deflection's differences [m]},every axis plot/.append style={thick}
]
\addplot[mark=*,red,style={dotted}]
  coordinates{
    (1,0.000750902565836521)
    (2,3.67963837432164e-06)
    (3,1.46496111253988e-06)
    (4,2.60225255743068e-07)
};
\addplot[mark=*,red,mark options={solid,fill =white},style={dotted}]
  coordinates{
    (1,0.00190064143211900)
    (2,0.000554991322418975)
    (3,0.000146606392671257)
    (4,3.78365955587788e-05)
};

\end{axis}

\end{tikzpicture}}
 \end{subfigure}
 \begin{subfigure}[b]{0.99\textwidth}
 \vspace{2mm}
\centering
        \scalebox{0.9}{
\begin{tikzpicture}
\begin{axis}[%
    hide axis,
    xmin=10,
    xmax=50,
    ymin=0,
    ymax=0.4,
    legend style={draw=white!15!black,legend cell align=left,legend columns=-1},
    every axis plot/.append style={thick}
    ]
    \addlegendimage{mark=*,blue}
    \addlegendentry{p-ref. Deflection};
    
     \addlegendimage{mark=*,blue,mark options={solid,fill opacity =0}}
    \addlegendentry{h-ref. Deflection};
    
    \addlegendimage{mark=*,red,style={dotted}}
    \addlegendentry{p-ref. Abs. val. of def. differences};
    
    \addlegendimage{mark=*,red,mark options={solid,fill opacity =0},style={dotted}}
    \addlegendentry{h-ref. Abs. val. of def. differences};
    \end{axis}
\end{tikzpicture}}
 \end{subfigure}
 \caption{Deflection and difference of the deflection of the middle node of the beam's top layer of nodes for the elastoplastic case (left) and elastic case (right).}
 \label{fig:mesh_convergence}
\end{figure}

The coarsest mesh is chosen so as to discretize the beam by at least four elements over the height  and  forty elements over the length. For cases where h-refinement is applied, the coarsest finite element mesh (level 0) consists of 410 degrees of freedom and a square element size of $0.0625\,\mathrm{m}$, while the finest finite element mesh considered (level 3) consists of 21186 degrees of freedom and a square element size of $0.0078\,\mathrm{m}$. The cases where p-refinement is applied have 410 degrees of freedom on their coarsest finite element mesh and 5474 on their finest mesh (level 3, quartic elements). The amount of elements stays the same for p-refinement.  In order to ensure the correct representation of the dynamic response, and to determine the minimum number of finite elements required, the bending wavelength, $\lambda_{\mathrm{min}}$,  can be evaluated from the Euler--Bernoulli beam theory. From this we obtain
\begin{linenomath*}
\begin{equation}
\lambda_{\mathrm{min}} = \sqrt{\dfrac{2 \pi}{f_{\mathrm{max}}}} \sqrt[4]{\dfrac{E I}{\rho A}},
\label{eq:MinElements}
\end{equation}
\end{linenomath*}
with $E$ the mean Young's modulus, $I$ the moment of inertia, $A$ the area, $\rho$ the density, $f_{\mathrm{max}}$ the highest simulated frequency, and $\lambda_{\mathrm{min}}$ the smallest obtained wavelength for the highest input frequency. For the considered beam configuration, it has been checked that at least six elements are used to represent the wavelength on the coarsest grid for the highest simulated frequency, which in this case is 400 Hz.
\\
For all elastic calculations, the MLMC/MLQMC simulations are level adaptive. For the elastoplastic cases we chose to manually set the maximum level because of the considerable time cost it would require to compute a solution on these higher levels. This level is chosen based on a mesh convergence analysis. The results of this mesh convergence study are shown in Fig.\,\ref{fig:mesh_convergence}, for the elastoplastic case (left) and the elastic case (right) respectively. The figures show the transverse deflection of the middle node located on the beam's top layer of nodes (middle top side node) per level, represented as a full line, and the absolute value of its difference over the levels, represented as a  dashed  line. For the elastoplastic case, the deflection starts stagnating at around level 3. Following these results, we thus state that the bias condition for the elastoplastic case is fulfilled at level 3; no more than 4 MLMC/MLQMC levels are used.

\subsection{Uncertainty bounds on the solution}
\label{UncertaintyInSol}
In this part, we illustrate the effect of the uncertainty model, i.e., a homogeneous or a heterogeneous Young's modulus,  on the uncertainty bounds of the solution. We do this for the static cases, i.e., the displacement of the beam in the spatial domain and for the dynamic case, i.e., the frequency response functions (FRF). First, the static elastic case is presented. Here, we chose to show a visualization of the transverse displacement of the nodes along the top side of the beam.  Second, results for the static elastoplastic case will be shown. The results are visualized by a force deflection curve of the middle top side node of the beam. Third,  solutions for the dynamic case are  presented. These consist of frequency response functions for a single node of the finite element mesh. This node is chosen as the one that has the largest response of all the nodes that make up the mesh. All results are presented with their uncertainty bounds.

\subsubsection{Static Elastic case}\label{Sec:Viz_Sta_El}

\begin{figure}[H]
\begin{subfigure}[b]{0.54\textwidth}
\centering
\scalebox{0.45}{
	\makebox[\textwidth]{
	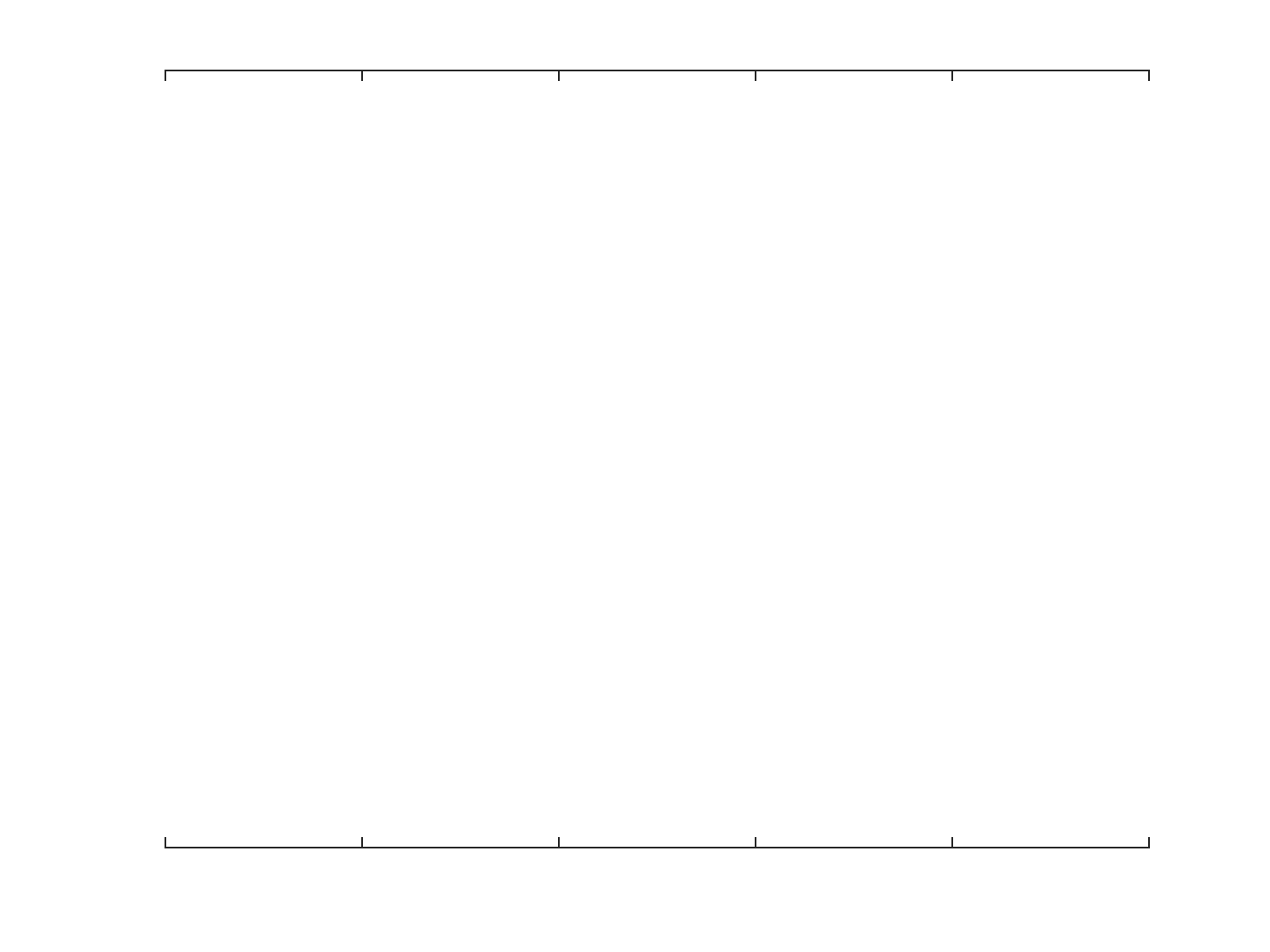}}
\end{subfigure}
\begin{subfigure}[b]{0.45\linewidth}
\centering
	\scalebox{0.45}{
		\makebox[\textwidth]{
	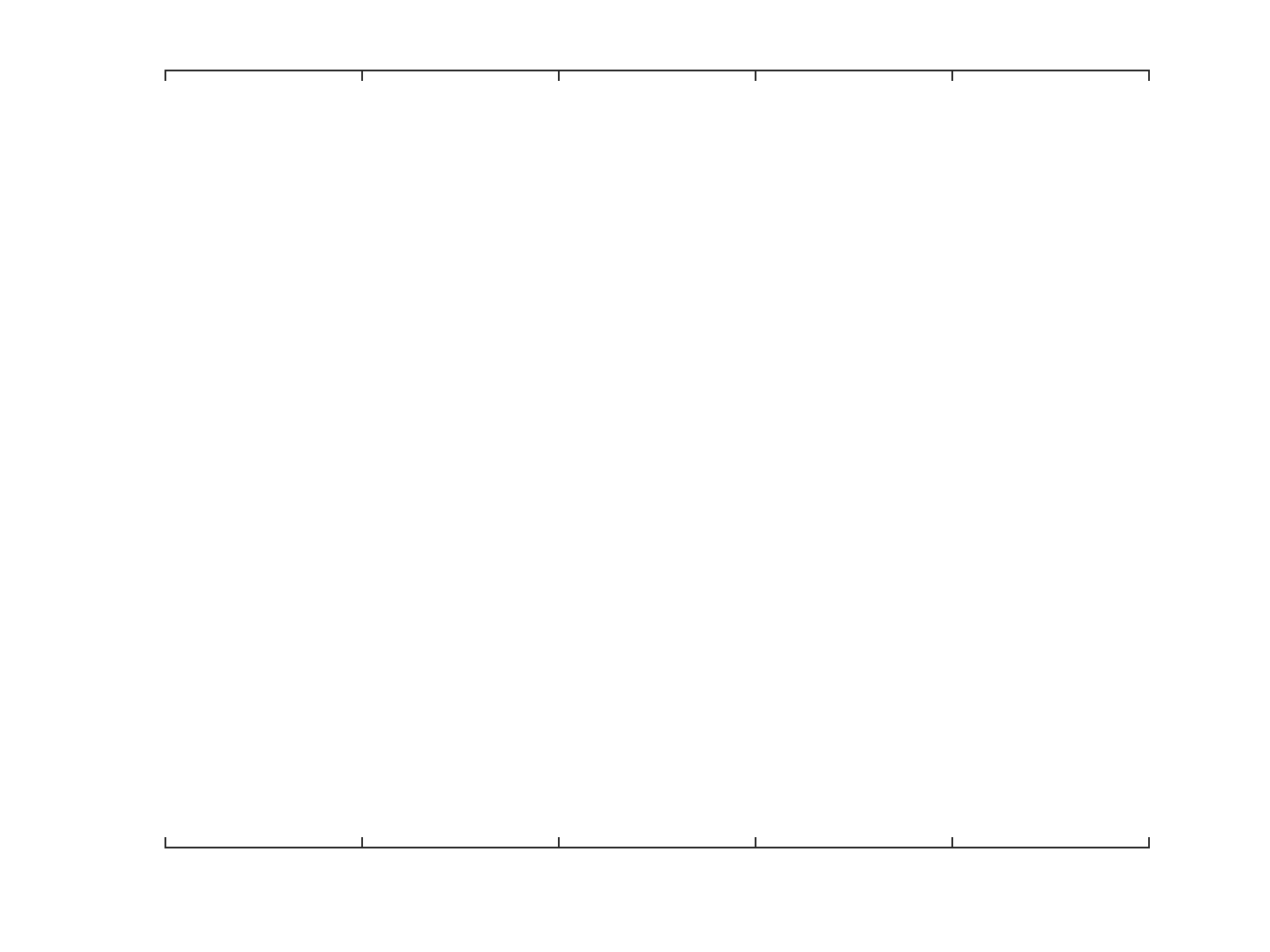}}
\end{subfigure}
\caption{Deflection of the beam for when the Young's modulus is homogeneous (left) and  heterogeneous (right).}
\label{fig:deflection_total}
\end{figure}

\begin{figure}[H]
\begin{subfigure}[b]{0.54\textwidth}
\centering
\scalebox{0.44}{
	\makebox[\textwidth]{
	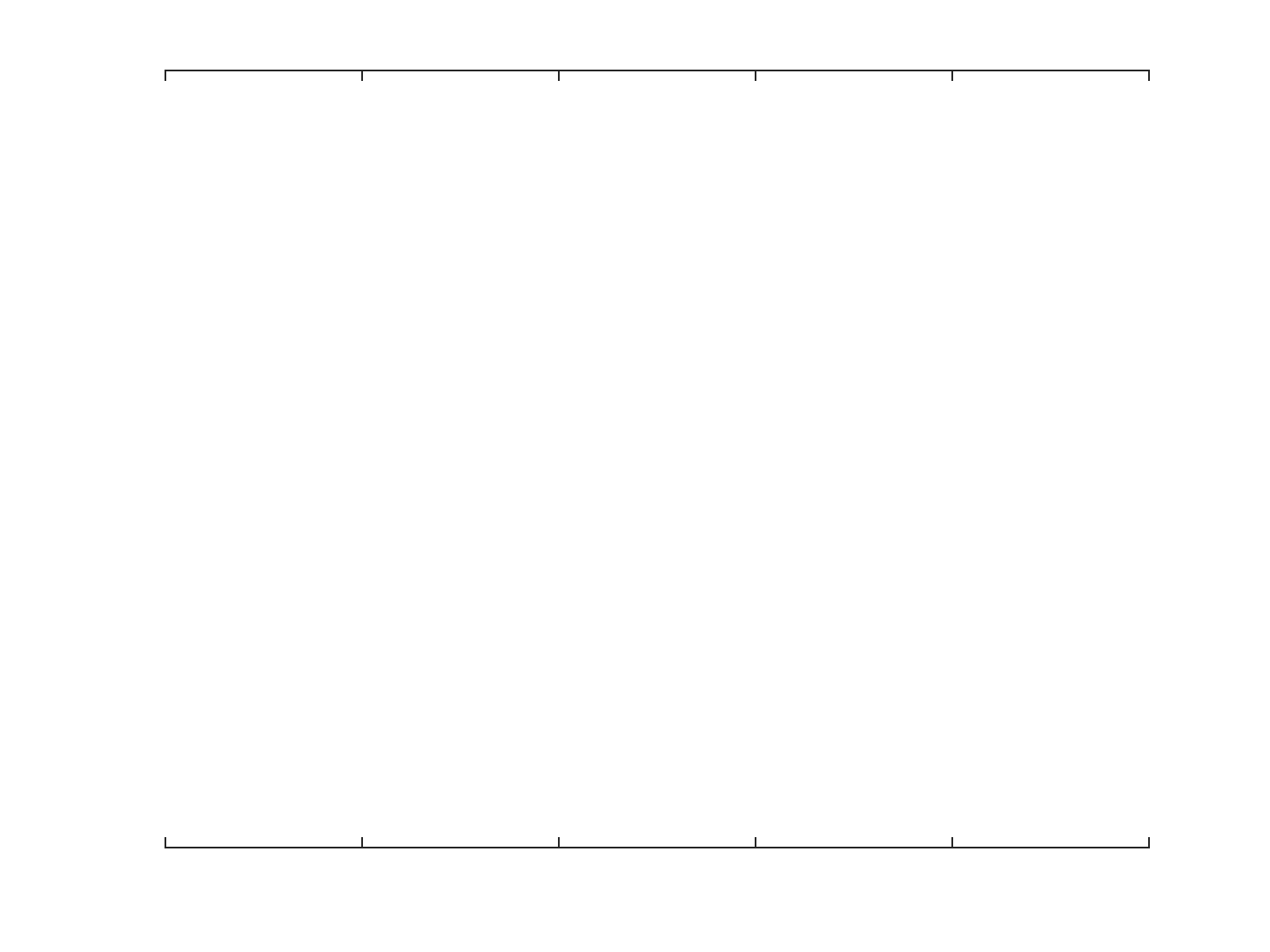}}
\end{subfigure}
\begin{subfigure}[b]{0.45\linewidth}
\centering
	\scalebox{0.44}{
		\makebox[\textwidth]{
	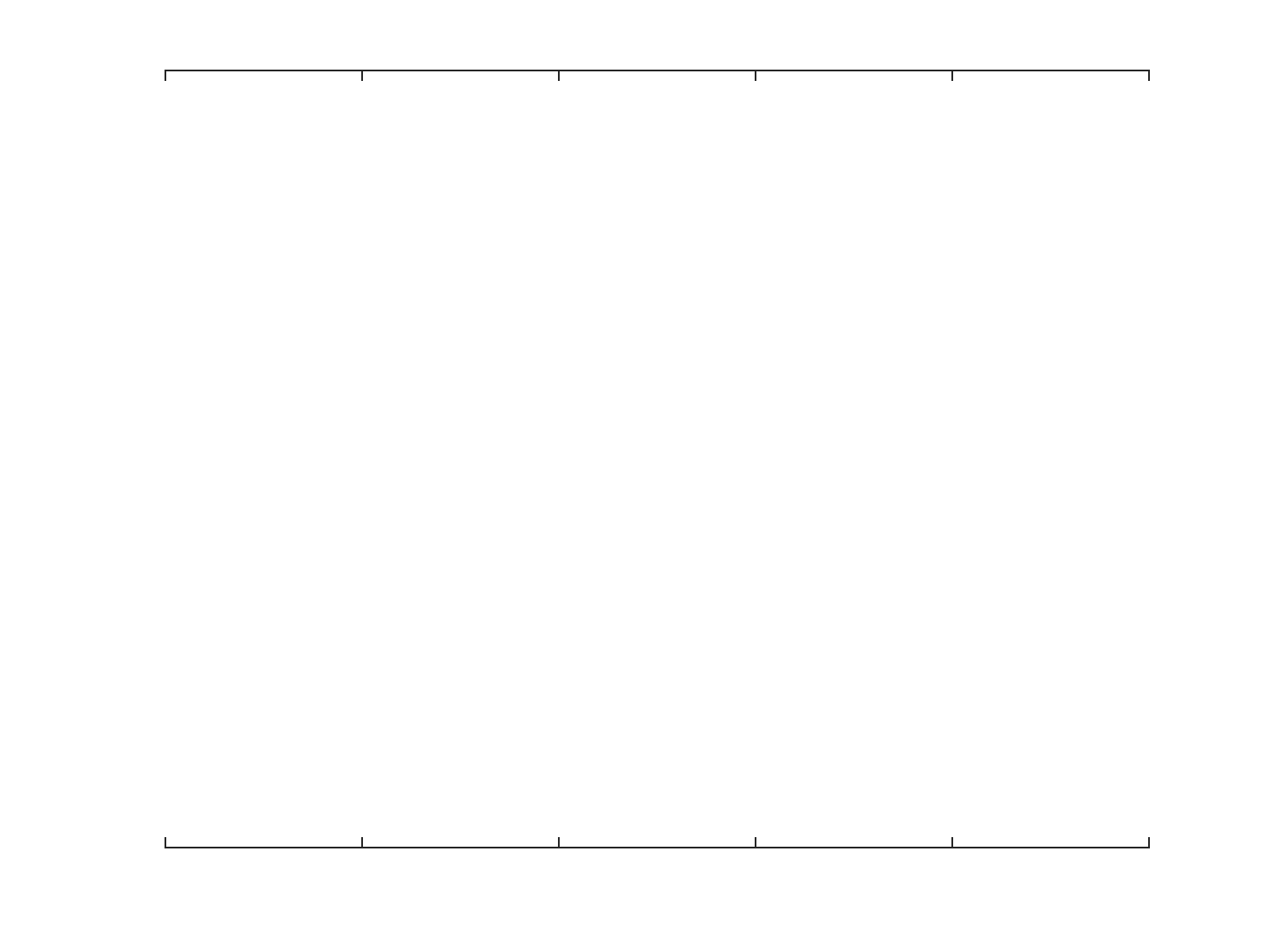}}
\end{subfigure}
\caption{Ten different deflection samples for when the Young's modulus is homogeneous (left) and heterogeneous (right).}
\label{fig:indiv_samples_elast}
\end{figure}

Fig.\,\ref{fig:deflection_total} shows the deflection of the beam with a homogeneous Young's modulus (left) and  with a heterogeneous Young's modulus (right). The orange solid line represents the average of the displacement, the  orange dashed lines are the 1$\sigma$ bounds equidistant around the average, which are primarily relevant in case of a normal distribution. The shades of blue represent the PDF, with the dark blue line corresponding to the most probable value.
Fig.\,\ref{fig:deflection_total} (right) shows that the average value and the most probable value of the  PDF tend to coincide for the heterogeneous case. Here, the PDF of the displacement closely resembles that of a normal distribution. This is however not the case when the Young's modulus is homogeneous. Then, the distribution of the solution has a non-negligible skewness. 
The homogeneous Young's modulus case exhibits a larger uncertainty on its displacement due to the fact that the Young's modulus is uncertain for each individual computed sample but uniform in each point for that sample. Averaging all these individual samples gives rise to wider uncertainty bounds. While for the heterogeneous Young's modulus case, each individual computed sample is also uncertain but non-uniform, in each individual point the value for the Young's modulus is different. This means that in different locations of the beam the Young's modulus will be different. These locations will tend to compensate each other so that the overall response of the beam does not become overly stiff or weak. Fig.\,\ref{fig:indiv_samples_elast} shows ten samples which illustrate this effect for a homogeneous Young's modulus (left) and a heterogeneous Young's modulus (right).

\begin{figure}[H]
\begin{subfigure}[b]{0.54\linewidth}
\centering
	\scalebox{0.44}{
		\makebox[\textwidth]{
	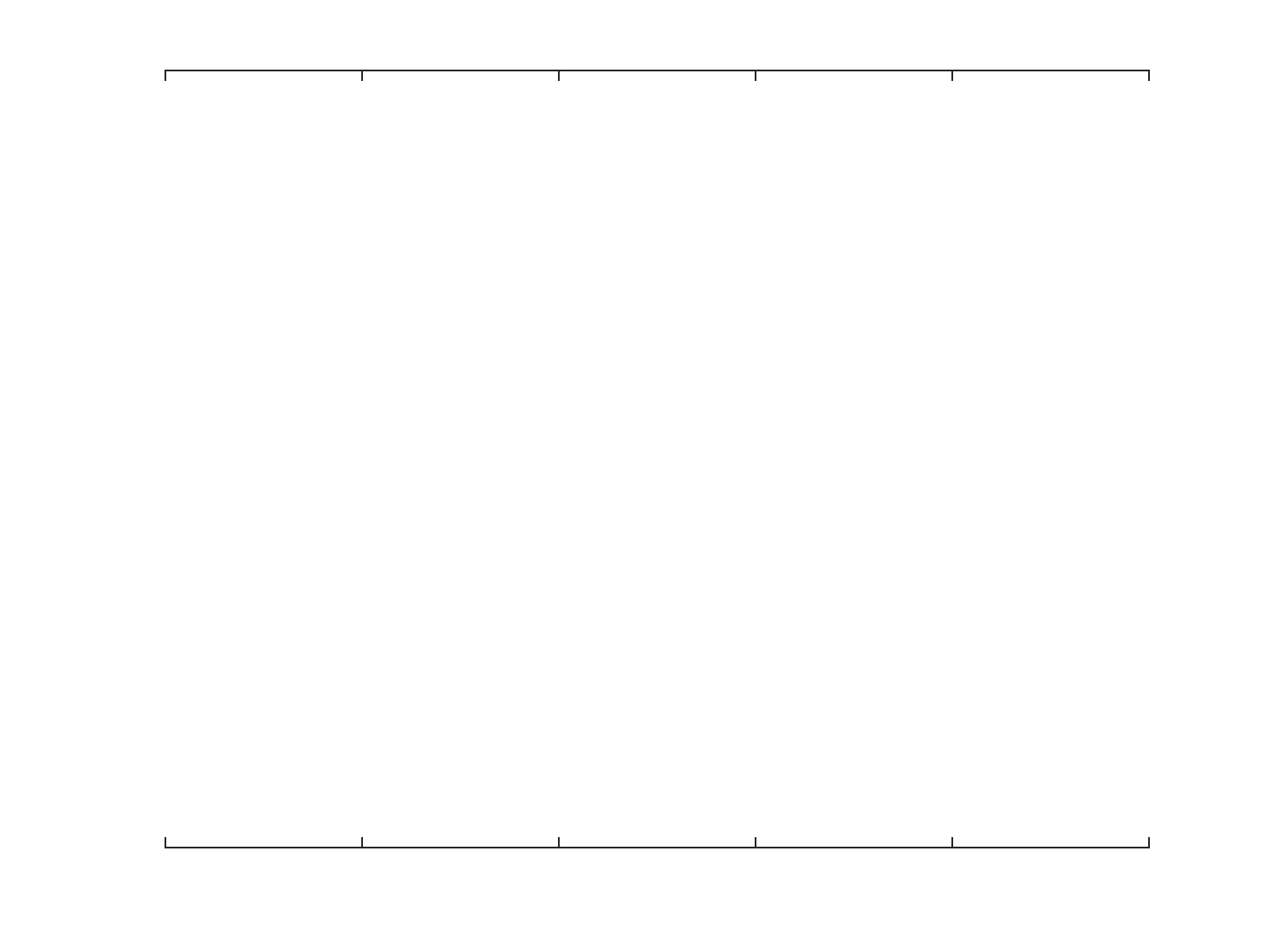}}
\end{subfigure}
\begin{subfigure}[b]{0.45\textwidth}
\centering
\scalebox{0.44}{
	\makebox[\textwidth]{
	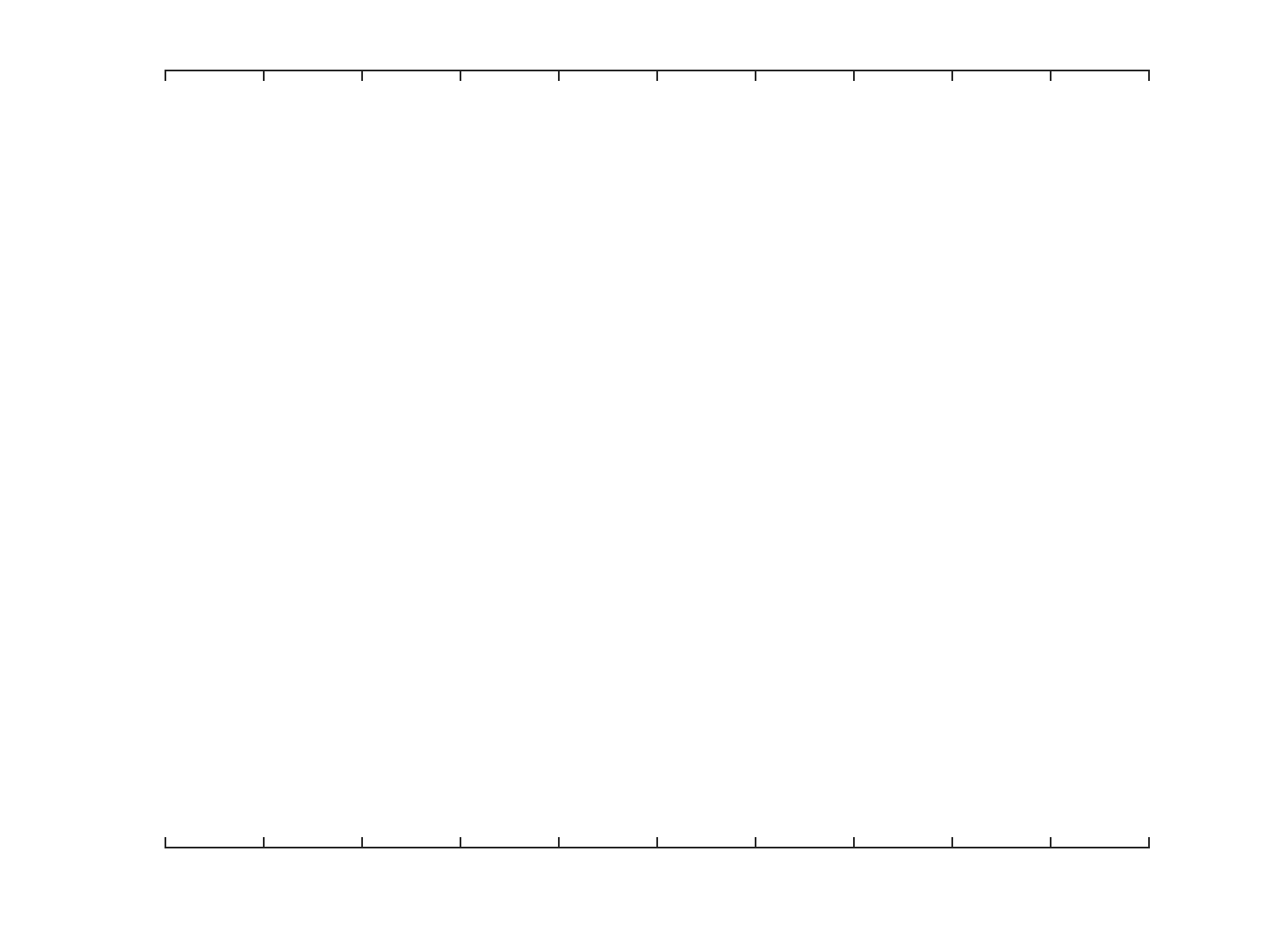}}
\end{subfigure}
\caption{Visualization of the PDF  when the Young's modulus is homogeneous, beam displacement (left), AB cut-through (right).}
\label{fig:Construction of the PDF Gam Dist}
\end{figure}

\begin{figure}[H]
\begin{subfigure}[b]{0.54\linewidth}
\centering
	\scalebox{0.44}{
		\makebox[\textwidth]{
	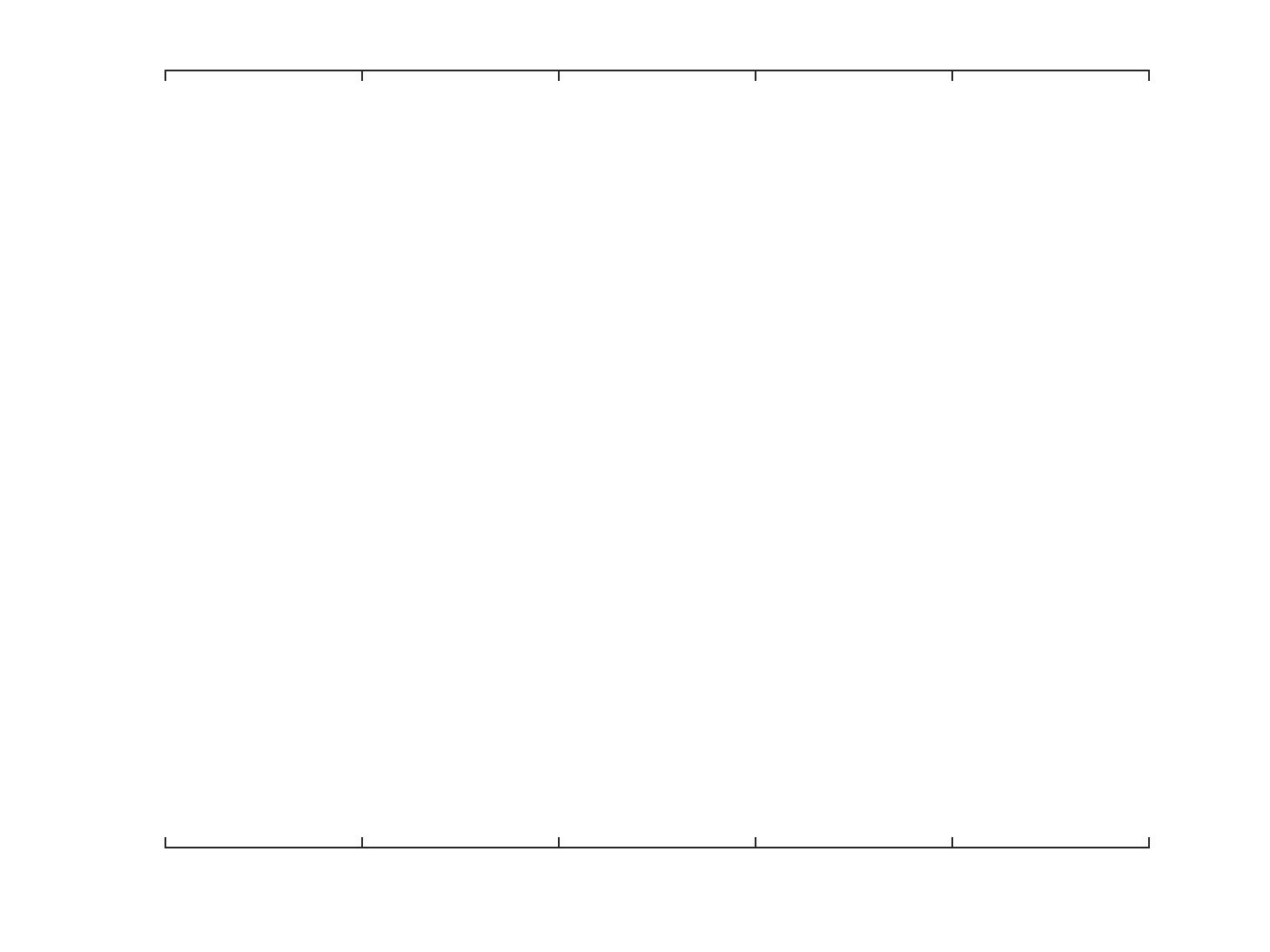}}
\end{subfigure}
\begin{subfigure}[b]{0.45\textwidth}
\centering
\scalebox{0.44}{
	\makebox[\textwidth]{
	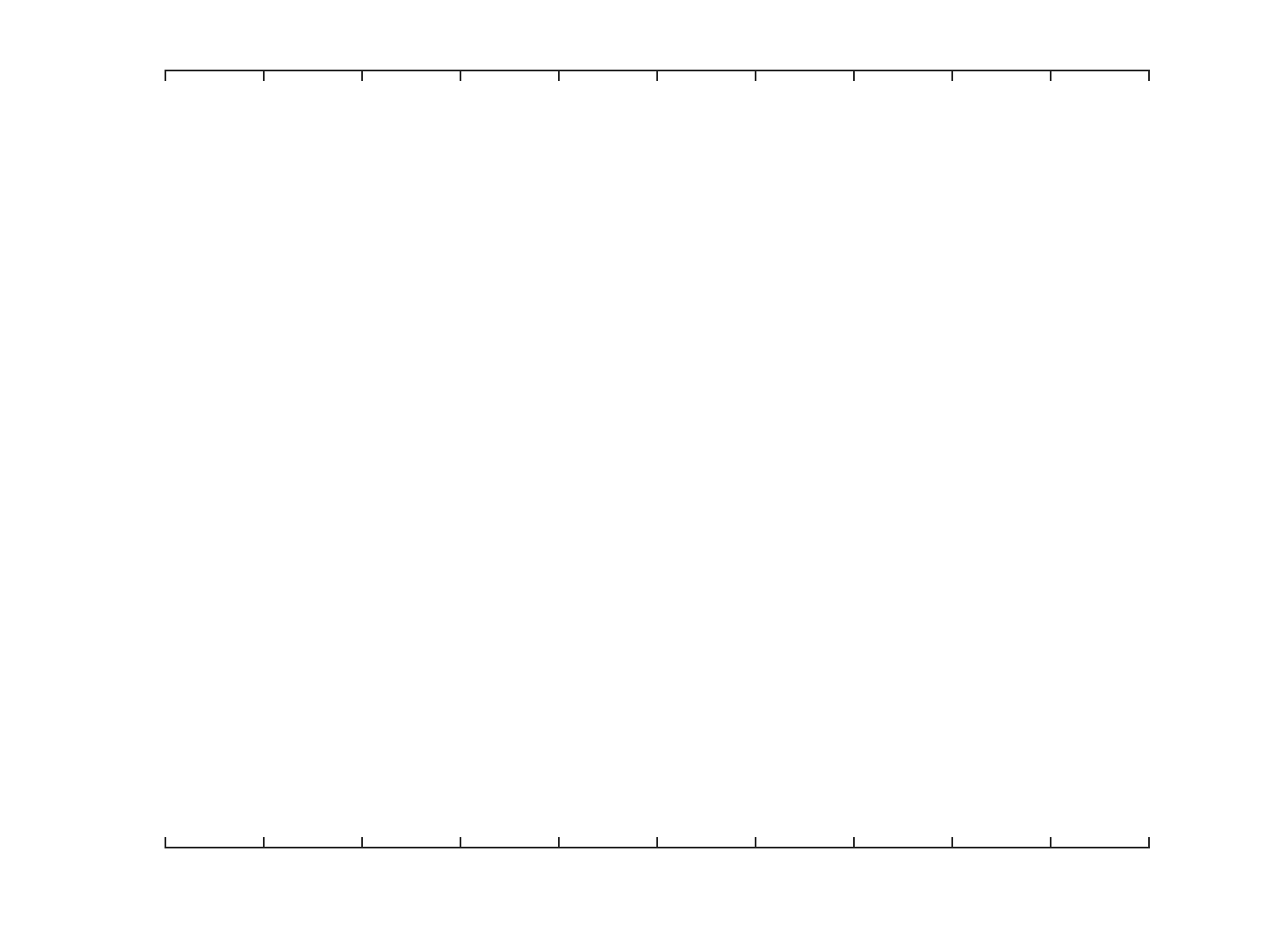}}
\end{subfigure}
\caption{Visualization of the PDF  when the Young's modulus is heterogeneous, beam displacement (left), AB cut-through (right).}
\label{fig:Construction of the PDF Gam Field}
\end{figure}

Fig.\,\ref{fig:Construction of the PDF Gam Dist} and Fig.\,\ref{fig:Construction of the PDF Gam Field} show a cut-through in order to better illustrate how the shades of blue represent the PDF.

\subsubsection{Static Elastoplastic case}

\begin{figure}[H]
\begin{subfigure}[b]{0.54\textwidth}
\centering
\scalebox{0.45}{
	\makebox[\textwidth]{
	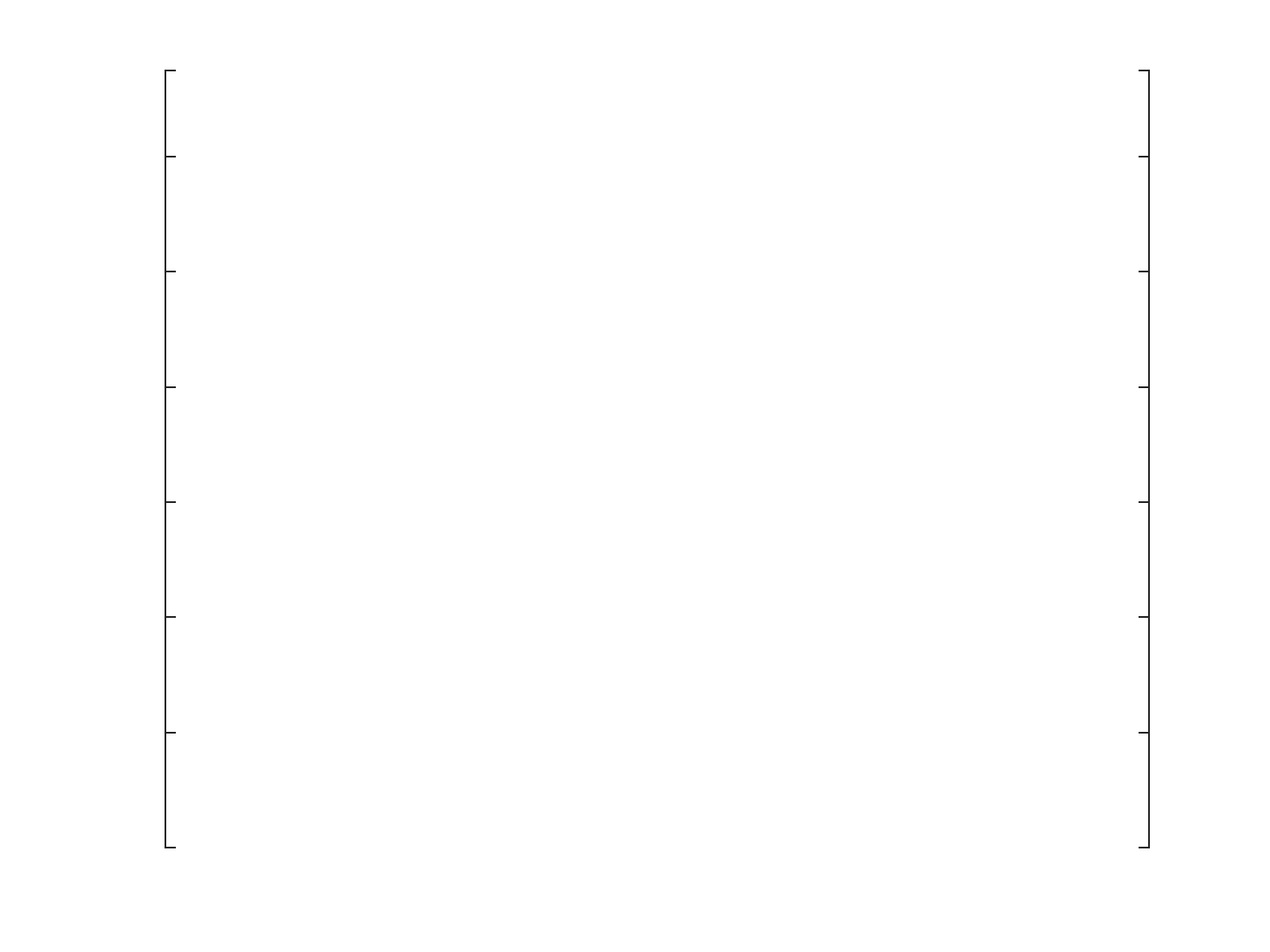}}
\end{subfigure}
\begin{subfigure}[b]{0.45\linewidth}
\centering
	\scalebox{0.45}{
		\makebox[\textwidth]{
	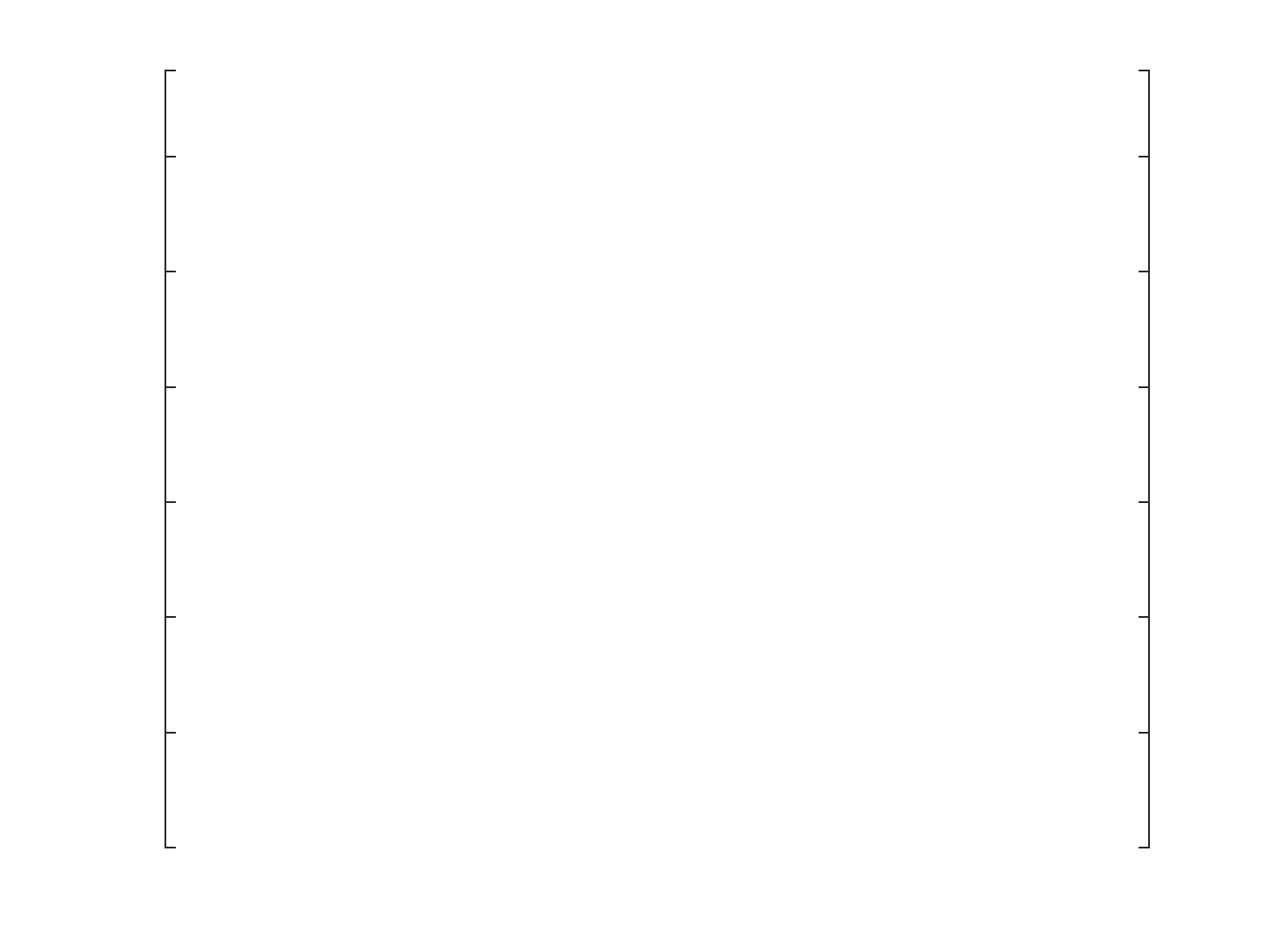}}
\end{subfigure}
\caption{Force deflection curve of the steel beam for when the Young's modulus is homogeneous (left) and  heterogeneous (right).}
\label{fig:deflection_total_plast}
\end{figure}
\begin{figure}[H]
\begin{subfigure}[b]{0.54\textwidth}
\centering
\scalebox{0.44}{
	\makebox[\textwidth]{
	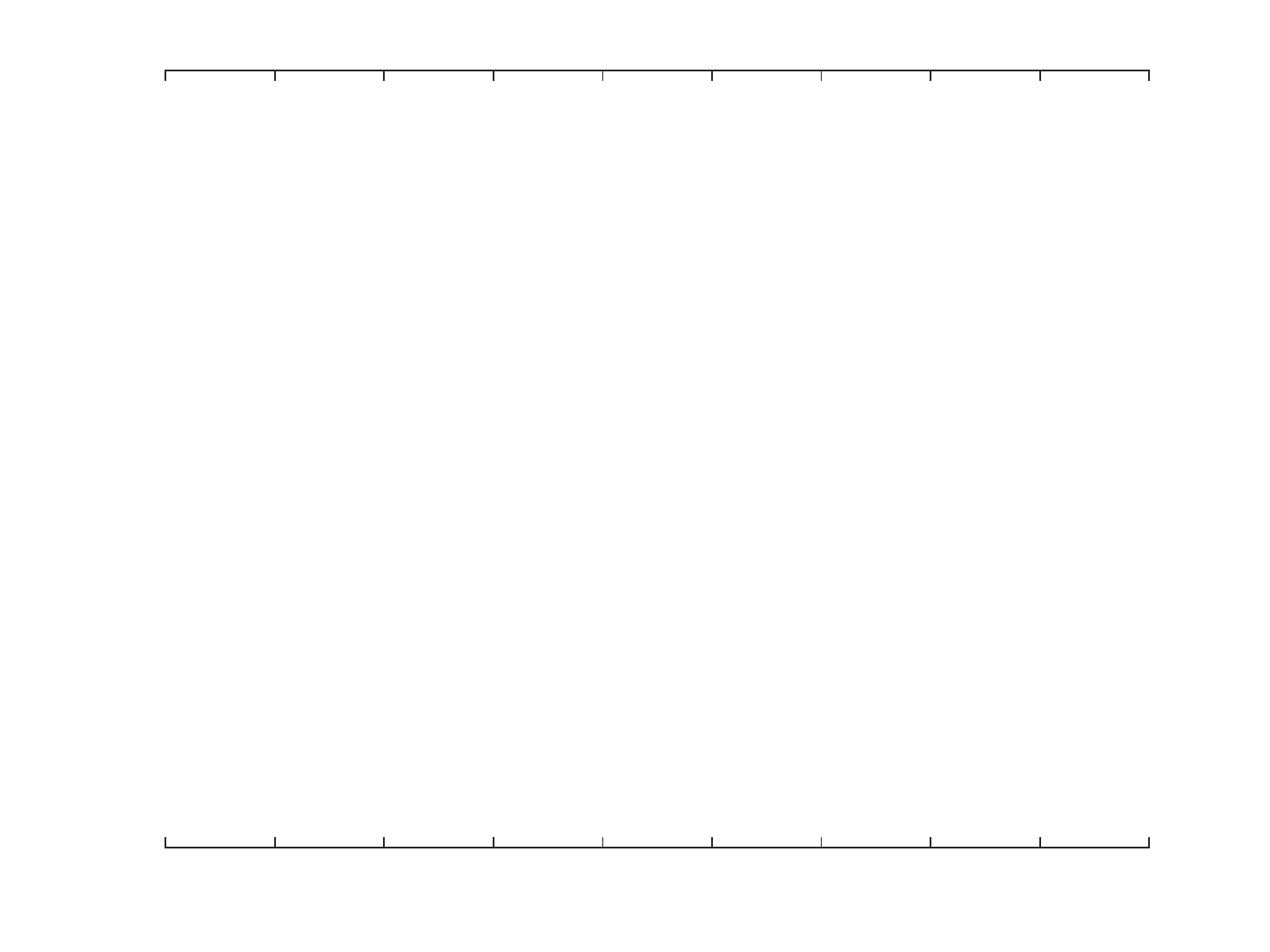}}
\end{subfigure}
\begin{subfigure}[b]{0.45\linewidth}
\centering
	\scalebox{0.44}{
		\makebox[\textwidth]{
	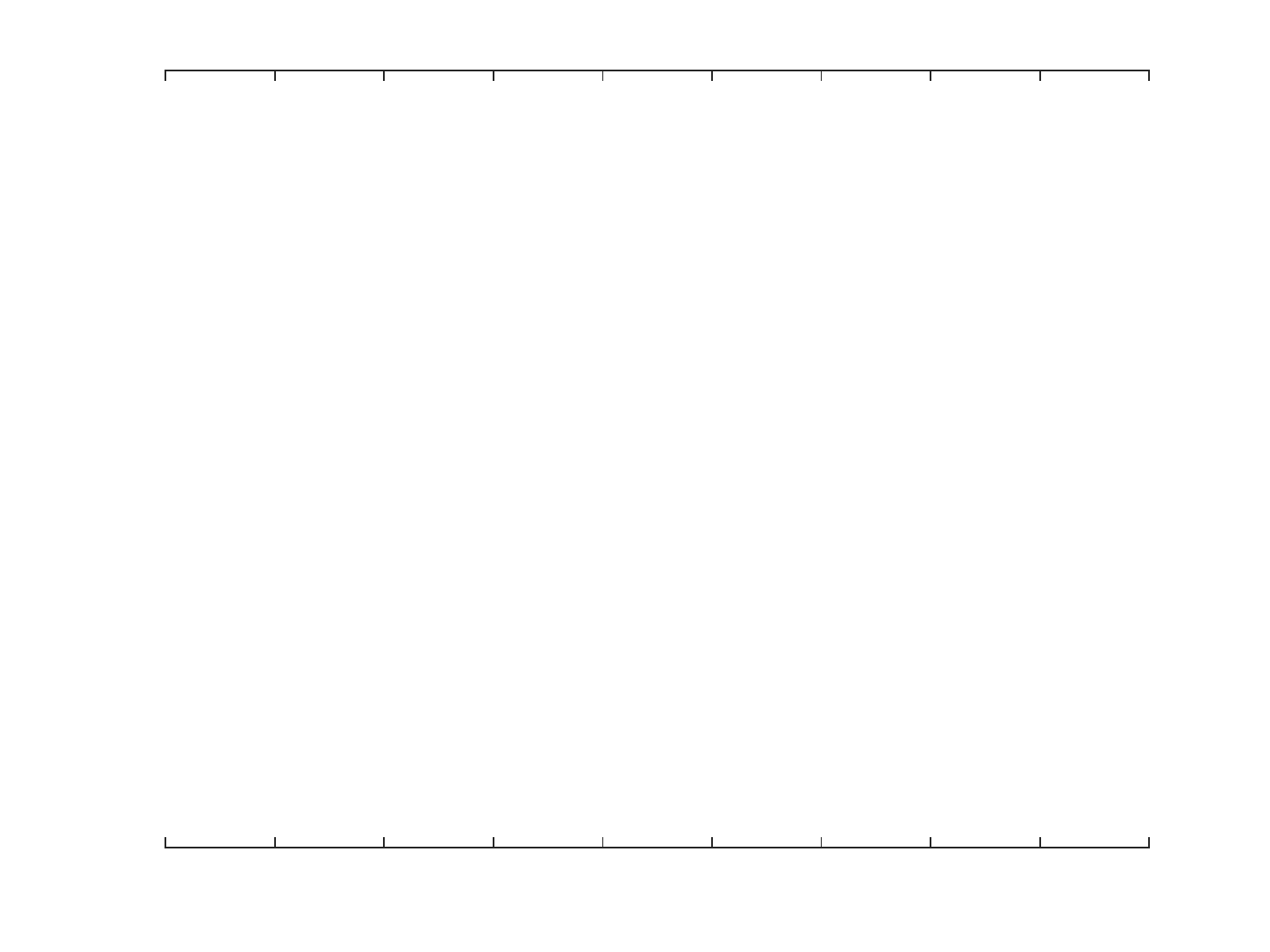}}
\end{subfigure}
\caption{Ten different force deflection curve samples for when the Young's modulus is homogeneous (left) and heterogeneous (right).}
\label{fig:indiv_samples_plast}
\end{figure}

Fig.\,\ref{fig:deflection_total_plast} shows the force deflection curve of the middle top side in case of a homogeneous Young's modulus (left) and a heterogeneous Young's modulus (right).
The line style and color convention is the same as for the static elastic case. As can be observed, the uncertainty bounds in case of a homogeneous Young's modulus, Fig.\,\ref{fig:deflection_total_plast} (left) are wider and more spread out than in case of a heterogeneous modulus (right). This behavior corroborates the one from the  static elastic case. In Fig.\,\ref{fig:indiv_samples_plast}, ten individual samples are shown for a homogeneous Young's modulus (left) and a heterogeneous Young's modulus (right).

\subsubsection{Dynamic Elastic case} 

\begin{figure}[H]
\begin{subfigure}[b]{0.54\textwidth}
\centering
\scalebox{0.46}{
	\makebox[\textwidth]{
	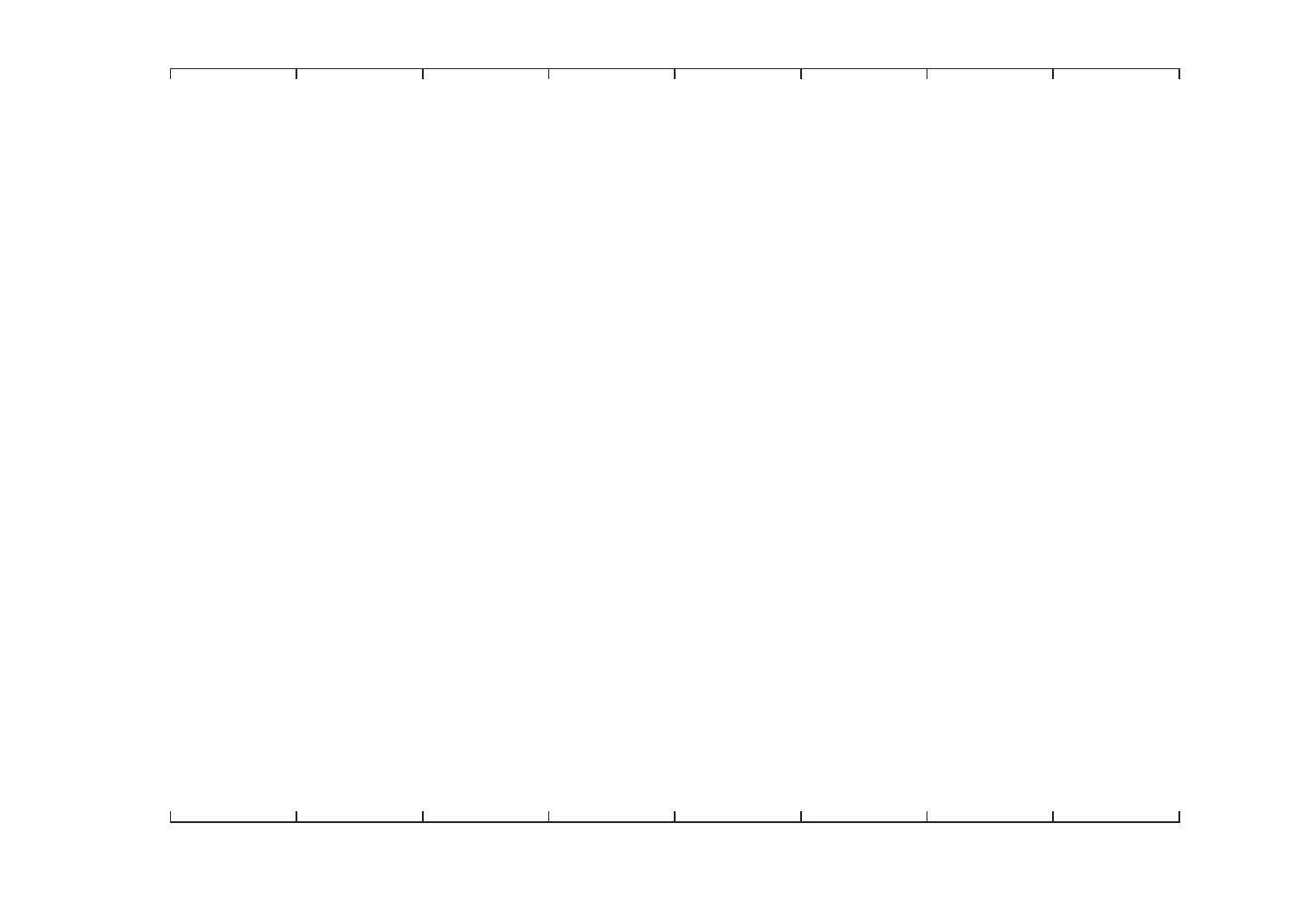}}
\end{subfigure}
\begin{subfigure}[b]{0.45\linewidth}
\centering
	\scalebox{0.43}{
		\makebox[\textwidth]{
	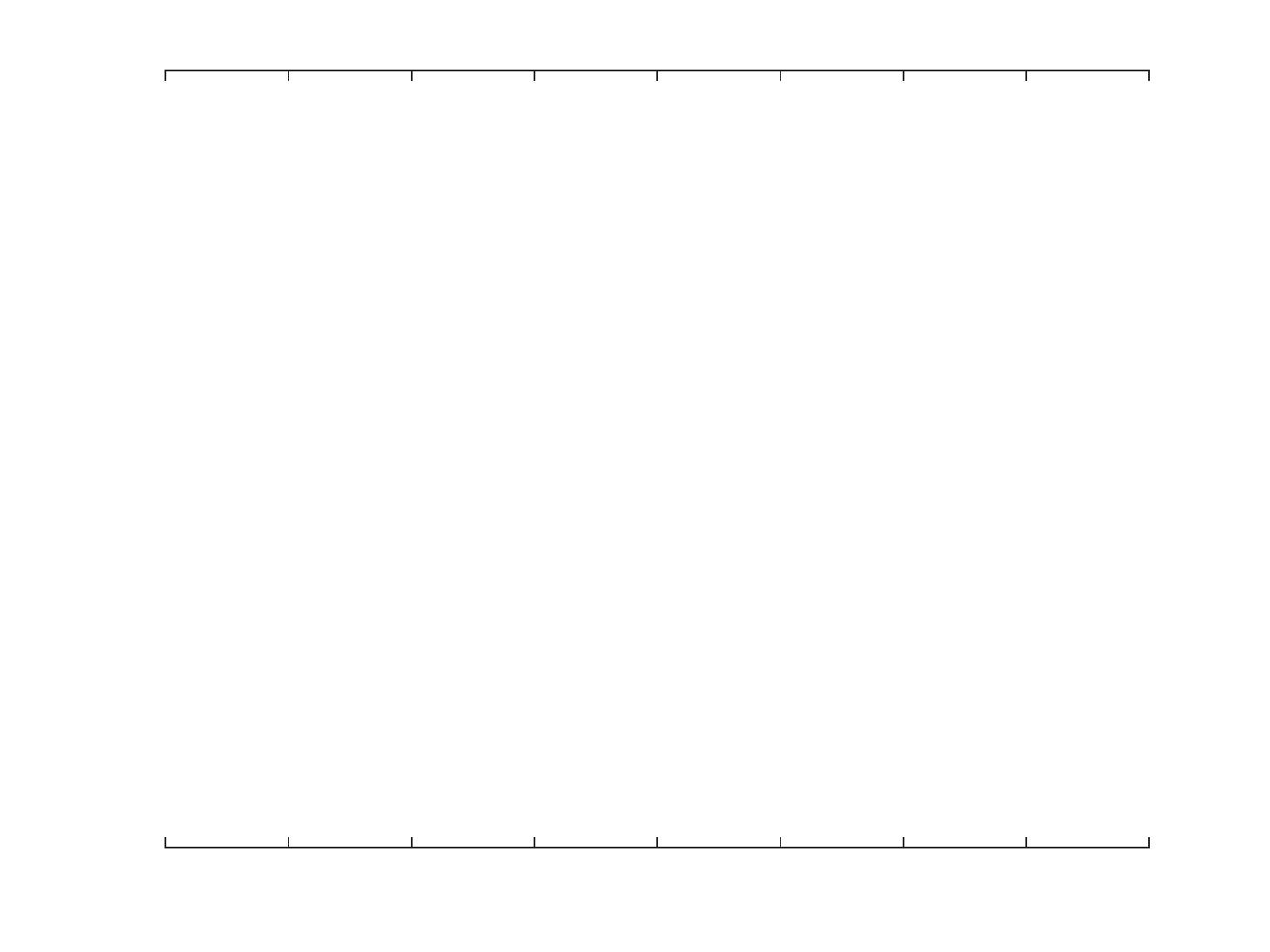}}
\end{subfigure}
\caption{Dynamic Responses of the beam for when the Young's modulus is  homogeneous (left) and  heterogeneous (right).}
\label{fig:dynamic_total}
\end{figure}
\begin{figure}[H]
\begin{subfigure}[b]{0.54\textwidth}
\centering
\scalebox{0.44}{
	\makebox[\textwidth]{
	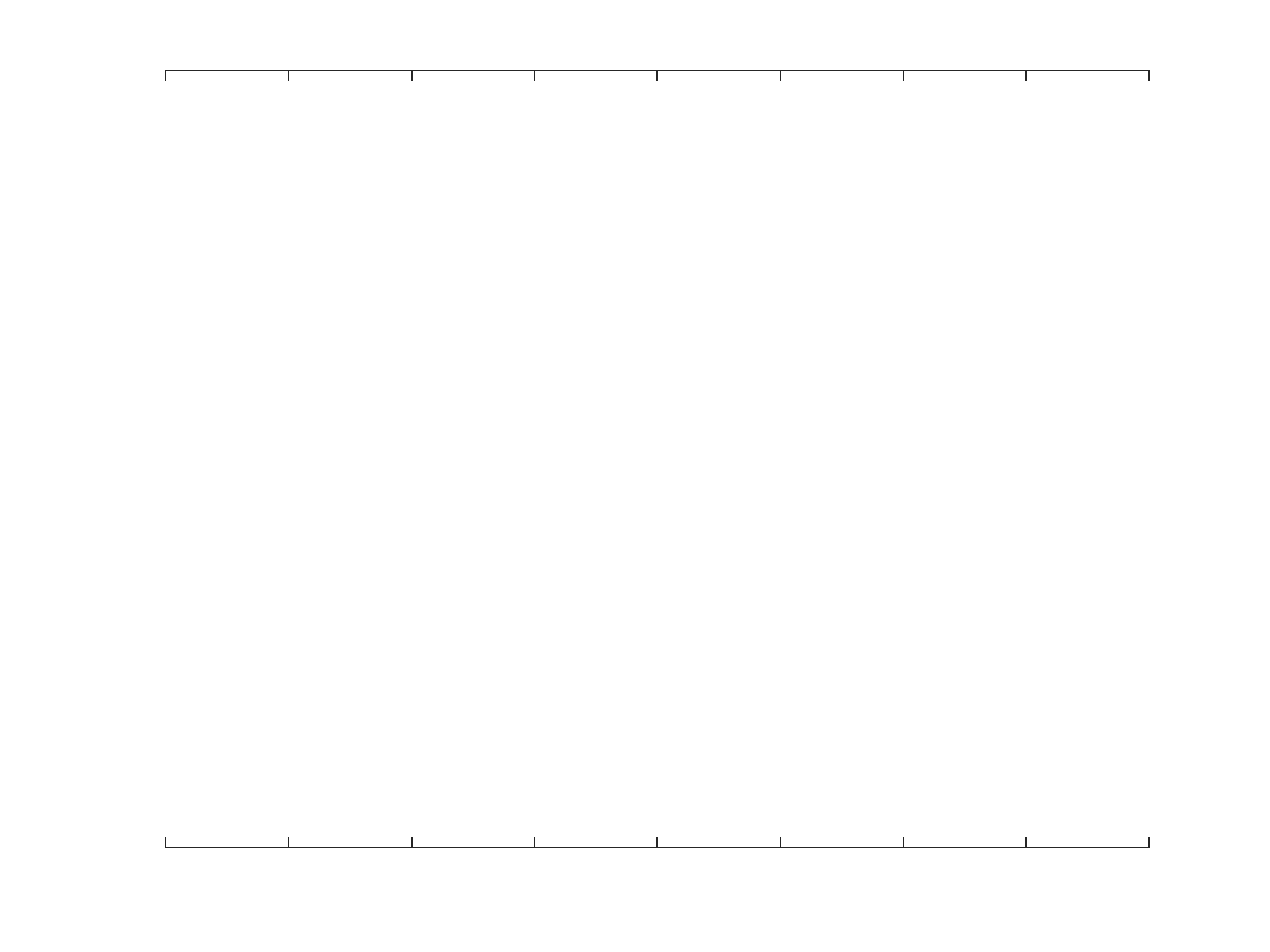}}
\end{subfigure}
\begin{subfigure}[b]{0.45\linewidth}
\centering
	\scalebox{0.44}{
		\makebox[\textwidth]{
	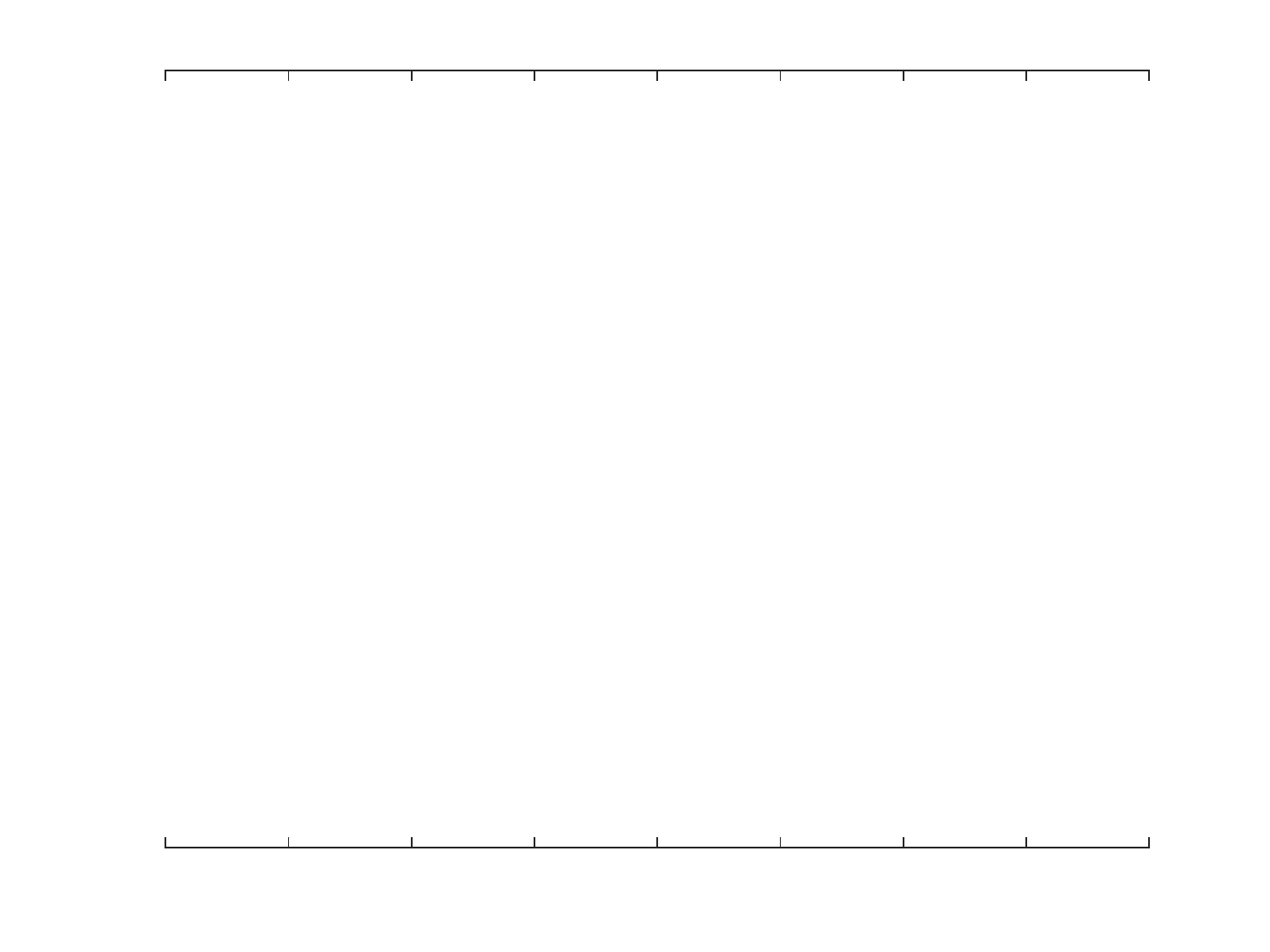}}
\end{subfigure}
\caption{Ten different FRF samples for when the Young's modulus is homogeneous (left) and heterogeneous (right).}
\label{fig:dynamic_samples}
\end{figure}

The FRF results are presented in Fig.\,\ref{fig:dynamic_total}. As was the case for the static elastic and elastoplastic case, the shades of blue represent the PDF, with the blue line being the most probable value, and the orange line the average value. As can be observed, the uncertainty bounds for the FRF  are wider and more spread out when the Young's modulus is homogeneous, Fig.\,\ref{fig:dynamic_total} (left), as opposed to a heterogeneous Young's modulus,  Fig.\,\ref{fig:dynamic_total} (right). This discrepancy is due to the fact that in case of a homogeneous Young's modulus, the resonance frequency will be shifted for each different sample. Averaging all these samples gives rise to a broad and wide uncertainty bound.  In case of a heterogeneous Young's modulus, the different samples compensate each other, in analogy with the explanation given in $\S$\ref{Sec:Viz_Sta_El}. This gives rise to much smaller uncertainty bounds. Fig.\,\ref{fig:dynamic_samples} shows the resulting FRF for ten  realizations.

An important inequality that must hold for the multilevel methods to work well is
\begin{linenomath*}
\begin{equation}
\V[P_1-P_0] \ll \V[P_1].
\label{eq:Variance_ineq}
\end{equation}
\end{linenomath*}
It has been observed empirically that Eq.\,\eqref{eq:Variance_ineq} is not necessarily fulfilled near resonance frequencies for the dynamic elastic case in case of a heterogeneous Young's modulus. In order to remedy to this, we use the following strategy. Using a small number of samples, the magnitude of $\V[P_1] - \V[P_1-P_0]$ is estimated. When the estimation is below a certain threshold, the coarsest level is discarded  and the algorithm is restarted on a finer mesh. We write this condition as
\begin{linenomath*}
\begin{equation}
\log_2\left(\dfrac{\V[P_1]}{\V[P_1-P_0]}\right) > {T}.
\label{eq:tresh}
\end{equation}  
\end{linenomath*} 
For the experiments reported here, we selected ${T}$ to be equal to $2.3$.
\subsection{Benchmark analysis}
Having illustrated the uncertainty propagation towards the solution, we now present a benchmark analysis where we compare the different Monte Carlo methods combined with both  refinement schemes in terms of computational cost.

\subsubsection{Rates}
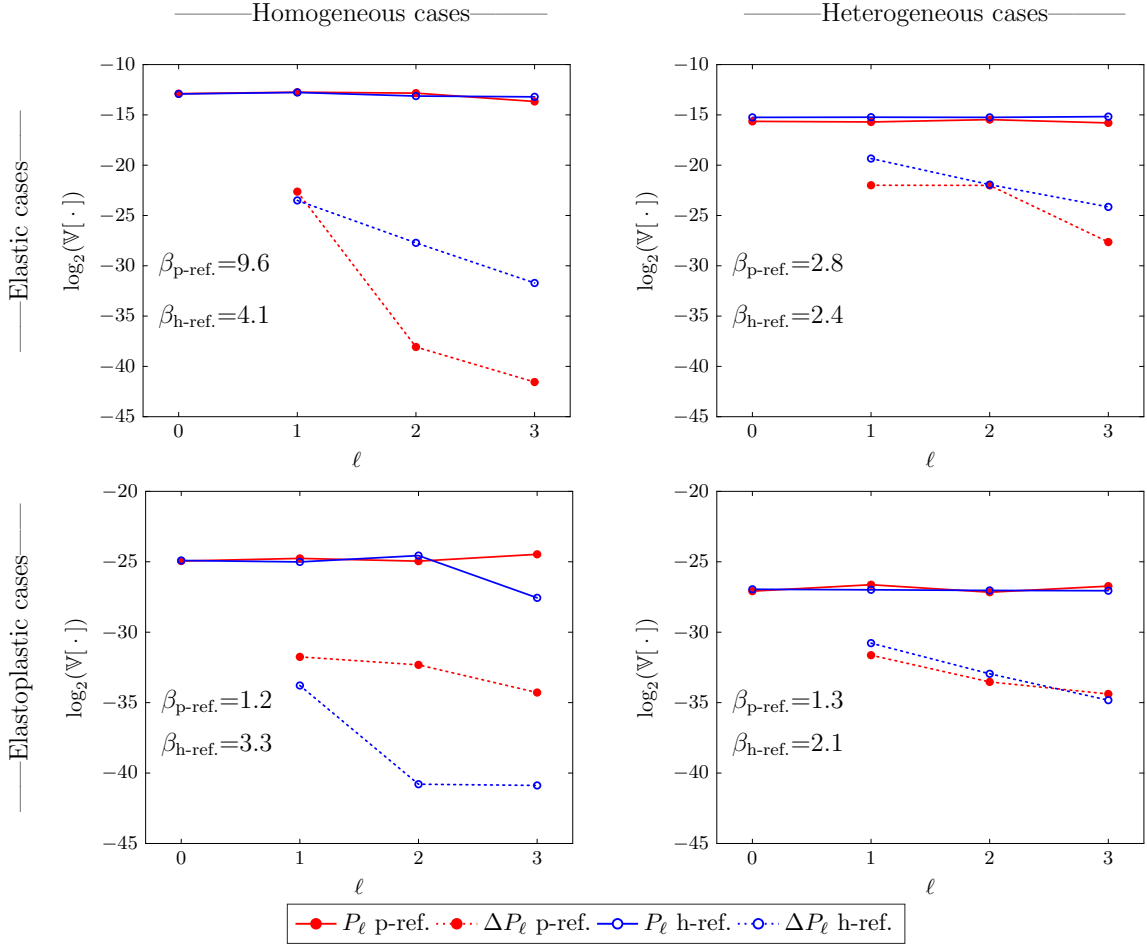
\begin{figure}[H]
\centering
        \begin{subfigure}[b]{0.475\textwidth}
        \scalebox{0.82}{
%
%
%
%
%
%
%
\begin{tikzpicture}
\begin{axis}[ticklabel style={{font=\small}}, major tick length={2pt}, every tick/.style={{black, line cap=round}}, axis on top, legend style={{draw=none, font=\small, at={(0.03,0.03)}, anchor=south west, fill=none, legend cell align=left}}, xlabel={$\ell$}, xtick distance={1},ytick distance= 5, ylabel={$\log_2(\mathbb{V}[\;\cdot\;])$},ymin=-45,ymax=-10,every axis plot/.append style={thick}]
    \addplot[mark={*}, mark size={1.5pt}, line cap={round}, mark options={solid}, color={red}]
        table[row sep={\\}]
        {
            \\
            0.0  -12.890092613421885  \\
            1.0  -12.745952188608218  \\
            2.0  -12.838092515773704  \\
            3.0  -13.667077690419077  \\
        }
        ;
    \addplot[mark={*}, mark size={1.5pt}, line cap={round}, mark options={solid}, color={red}, style={dotted}]
        table[row sep={\\}]
        {
            \\
            1.0  -22.630279186083595  \\
            2.0  -38.06871520608006  \\
            3.0  -41.555286375265965  \\
        }
        ;
    \addplot[mark={*}, mark size={1.5pt}, line cap={round}, mark options={solid, fill opacity=0}, color={blue}]
        table[row sep={\\}]
        {
            \\
            0.0  -12.919472599057826  \\
            1.0  -12.770465394114007  \\
            2.0  -13.127822512054623  \\
            3.0  -13.206448642743442  \\
        }
        ;
    \addplot[mark={*}, mark size={1.5pt}, line cap={round}, mark options={solid,fill opacity=0}, color={blue}, style={dotted}]
        table[row sep={\\}]
        {
            \\
            1.0  -23.502521609507205  \\
            2.0  -27.719335637642615  \\
            3.0  -31.710884945221107  \\
        }
        ;

\node[black] at (axis cs:0.3,-30){\large{$\beta_{\text{p-ref.}}$=9.6}};
\node[black] at (axis cs:0.3,-35){\large{$\beta_{\text{h-ref.}}$=4.1}};

\end{axis}
\node at (3.5,6.5) {\large---------Homogeneous cases---------};
\node at (-2,3) {\rotatebox{90}{\large------Elastic cases------}};

\end{tikzpicture}
}
 \end{subfigure}
 \hspace{2mm}
                \begin{subfigure}[b]{0.475\textwidth}
        \scalebox{0.82}{
%
%
%
%
%
%
%
\begin{tikzpicture}
\begin{axis}[ticklabel style={{font=\small}}, major tick length={2pt}, every tick/.style={{black, line cap=round}}, axis on top, legend style={{draw=none, font=\small, at={(0.03,0.03)}, anchor=south west, fill=none, legend cell align=left}}, xlabel={$\ell$}, xtick distance={1},ytick distance= 5, ylabel={$\log_2(\mathbb{V}[\;\cdot\;])$},ymin=-45,ymax=-10,every axis plot/.append style={thick}]
    \addplot[mark={*}, mark size={1.5pt}, line cap={round}, mark options={solid}, color={red}]
        table[row sep={\\}]
        {
            \\
            0.0  -15.643461946749243  \\
            1.0  -15.703801220530503  \\
            2.0  -15.470294314679464  \\
            3.0  -15.798780185517519  \\
        }
        ;
    \addplot[mark={*}, mark size={1.5pt}, line cap={round}, mark options={solid}, color={red}, style={dotted}]
        table[row sep={\\}]
        {
            \\
            1.0  -21.993959024821653  \\
            2.0  -22.007008922232018  \\
            3.0  -27.638947707083936  \\
        }
        ;
    \addplot[mark={*}, mark size={1.5pt}, line cap={round}, mark options={solid,fill opacity=0}, color={blue}]
        table[row sep={\\}]
        {
            \\
            0.0  -15.248382774603352  \\
            1.0  -15.230839878326721  \\
            2.0  -15.247680658948546  \\
            3.0  -15.170668149672565  \\
        }
        ;
    \addplot[mark={*}, mark size={1.5pt}, line cap={round}, mark options={solid,fill opacity=0}, color={blue}, style={dotted}]
        table[row sep={\\}]
        {
            \\
            1.0  -19.33901018094879  \\
            2.0  -21.932817654600285  \\
            3.0  -24.145949520236055  \\
        }
        ;
\node[black] at (axis cs:0.3,-30){\large{$\beta_{\text{p-ref.}}$=2.8}};
\node[black] at (axis cs:0.3,-35){\large{$\beta_{\text{h-ref.}}$=2.4}};
\end{axis}
\node at (3.5,6.5) {\large---------Heterogeneous cases---------};

\end{tikzpicture}
}
 \end{subfigure}
         \begin{subfigure}[b]{0.475\textwidth}
        \scalebox{0.82}{
%
%
%
%
%
%
%
\begin{tikzpicture}
\begin{axis}[ticklabel style={{font=\small}}, major tick length={2pt}, every tick/.style={{black, line cap=round}}, axis on top, legend style={{draw=none, font=\small, at={(0.03,0.03)}, anchor=south west, fill=none, legend cell align=left}}, xlabel={$\ell$}, xtick distance={1}, ylabel={$\log_2(\mathbb{V}[\;\cdot\;])$},ymax=-20,ymin=-45,every axis plot/.append style={thick}]
    \addplot[mark={*}, mark size={1.5pt}, line cap={round}, mark options={solid}, color={red}]
        table[row sep={\\}]
        {
            \\
            0.0  -24.947189560222125  \\
            1.0  -24.760630842404446  \\
            2.0  -24.95597157400712  \\
            3.0  -24.470471677028456  \\
        }
        ;
    \addplot[mark={*}, mark size={1.5pt}, line cap={round}, mark options={solid}, color={red}, style={dotted}]
        table[row sep={\\}]
        {
            \\
            1.0  -31.751730612769403  \\
            2.0  -32.31877675853436  \\
            3.0  -34.28312879317675  \\
        }
        ;
        \addplot[mark={*}, mark size={1.5pt}, line cap={round}, mark options={solid,fill opacity=0}, color={blue}]
        table[row sep={\\}]
        {
            \\
            0.0  -24.908166956950865  \\
            1.0  -25.00964143391462  \\
            2.0  -24.56108579300735  \\
            3.0  -27.559492960056637  \\
        }
        ;
    \addplot[mark={*}, mark size={1.5pt}, line cap={round}, mark options={solid,fill opacity=0}, color={blue}, style={dotted}]
        table[row sep={\\}]
        {
            \\
            1.0  -33.783614061344124  \\
            2.0  -40.792383795488654  \\
             3.0  -40.8785579956400  \\
        }
        ;
  \node[black] at (axis cs:0.3,-35){\large{$\beta_{\text{p-ref.}}$=1.2}};
\node[black] at (axis cs:0.3,-38){\large{$\beta_{\text{h-ref.}}$=3.3}};
\end{axis}
\node at (-2,3) {\rotatebox{90}{\large------Elastoplastic cases------}};

\end{tikzpicture}
}
 \end{subfigure}
  \hspace{2mm}
      \begin{subfigure}[b]{0.475\textwidth}
        \scalebox{0.82}{
%
%
%
%
%
%
%
\begin{tikzpicture}
\begin{axis}[ticklabel style={{font=\small}}, major tick length={2pt}, every tick/.style={{black, line cap=round}}, axis on top, legend style={{draw=none, font=\small, at={(0.03,0.03)}, anchor=south west, fill=none, legend cell align=left}}, xlabel={$\ell$}, xtick distance={1}, ylabel={$\log_2(\mathbb{V}[\;\cdot\;])$},ymax=-20,ymin=-45,every axis plot/.append style={thick}]
    \addplot[mark={*}, mark size={1.5pt}, line cap={round}, mark options={solid}, color={red}]
        table[row sep={\\}]
        {
            \\
            0.0  -27.085154474159093  \\
            1.0  -26.629636392570294  \\
            2.0  -27.160375381268313  \\
            3.0  -26.73089949020011  \\
        }
        ;
    \addplot[mark={*}, mark size={1.5pt}, line cap={round}, mark options={solid}, color={red}, style={dotted}]
        table[row sep={\\}]
        {
            \\
            1.0  -31.63187872257093  \\
            2.0  -33.52941273871351  \\
            3.0  -34.385908673615155  \\
        }
        ;
    
        \addplot[mark={*}, mark size={1.5pt}, line cap={round}, mark options={solid,fill opacity=0}, color={blue}]
        table[row sep={\\}]
        {
            \\
            0.0  -26.95274547222109  \\
            1.0  -26.990303587787206  \\
            2.0  -27.03401112791239  \\
            3.0  -27.053705713428492  \\
        }
        ;
    \addplot[mark={*}, mark size={1.5pt}, line cap={round}, mark options={solid,fill opacity=0}, color={blue}, style={dotted}]
        table[row sep={\\}]
        {
            \\
            1.0  -30.77681694755675  \\
            2.0  -32.9526543208583  \\
            3.0  -34.82571604921131  \\
        }
        ;
 \node[black] at (axis cs:0.3,-35){\large{$\beta_{\text{p-ref.}}$=1.3}};
\node[black] at (axis cs:0.3,-38){\large{$\beta_{\text{h-ref.}}$=2.1}};
\end{axis}
\end{tikzpicture}
}
 \end{subfigure}
 \begin{subfigure}[b]{0.99\textwidth}
\centering
        \scalebox{0.9}{
\begin{tikzpicture}
\begin{axis}[%
    hide axis,
    xmin=10,
    xmax=50,
    ymin=0,
    ymax=0.4,
    legend style={draw=white!15!black,legend cell align=left,legend columns=-1},
    every axis plot/.append style={thick}
    ]
    \addlegendimage{red,mark=*}
    \addlegendentry{$P_{\ell}$ p-ref.};
    
     \addlegendimage{red,mark=*,dotted,mark options = {solid}}
    \addlegendentry{$\Delta P_{\ell}$ p-ref.};
    
    \addlegendimage{blue,mark options = {solid,fill=white}, mark=*}
    \addlegendentry{$P_{\ell}$ h-ref.};
    
    \addlegendimage{blue,dotted,mark options = {solid,fill=white}, mark=*}
    \addlegendentry{$\Delta P_{\ell}$ h-ref.};
    \end{axis}
\end{tikzpicture}}
\vspace{5mm}
 \end{subfigure}
 \vspace{-3mm}
 \caption{Variance of the quantity of interest and variance of its differences for the homogeneous elastic and elastoplastic cases (top and bottom left) and the heterogeneous elastic and elastoplastic cases (top and bottom right).}
 \label{fig:variances}
\end{figure}

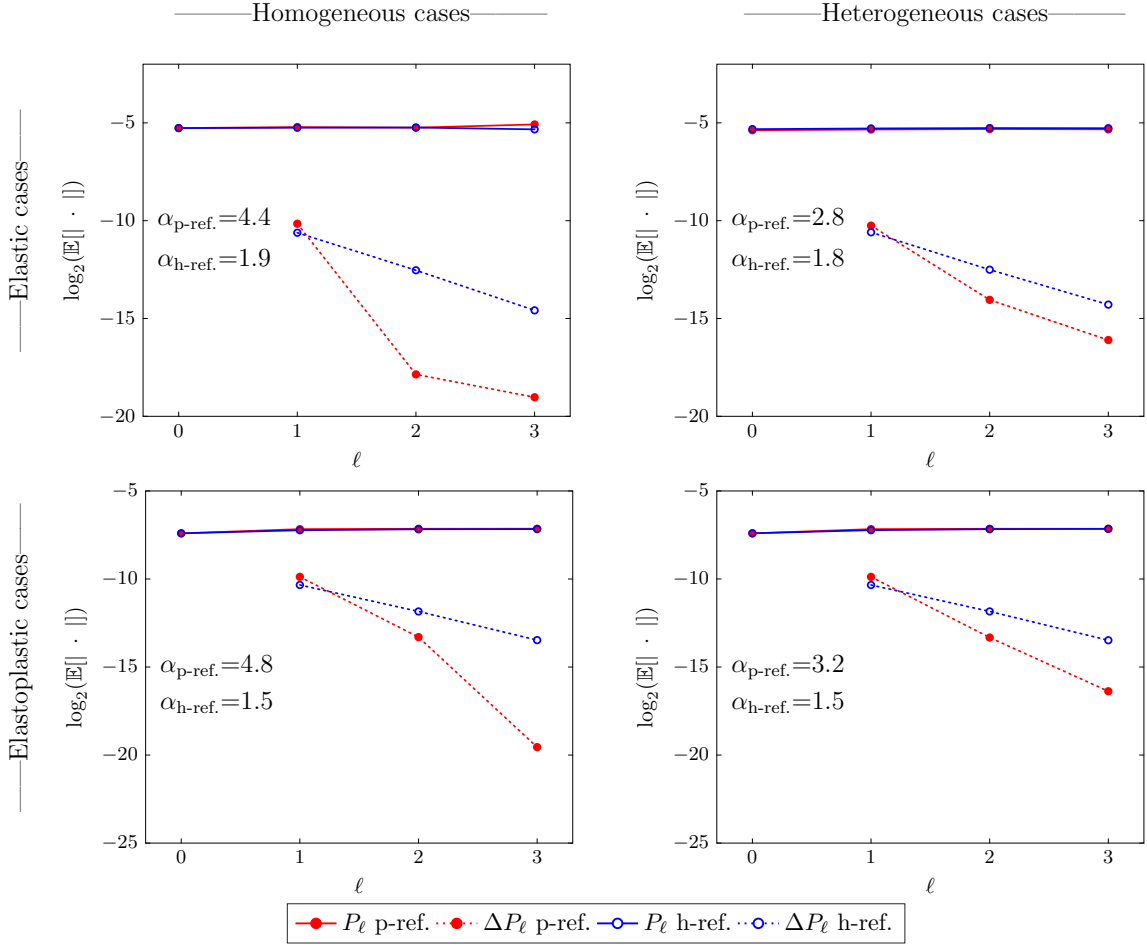
\begin{figure}[H]
\centering
         \begin{subfigure}[b]{0.475\textwidth}
        \scalebox{0.82}{
%
%
%
%
%
%
%
\begin{tikzpicture}
\begin{axis}[ticklabel style={{font=\small}}, major tick length={2pt}, every tick/.style={{black, line cap=round}}, axis on top, legend style={{draw=none, font=\small, at={(0.03,0.03)}, anchor=south west, fill=none, legend cell align=left}}, xlabel={$\ell$}, xtick distance={1}, ylabel={$\log_2(\mathbb{E}[|\;\cdot\;|])$},ymax=-2,ymin=-20,every axis plot/.append style={thick}]
    \addplot[mark={*}, mark size={1.5pt}, line cap={round}, mark options={solid}, color={red}]
        table[row sep={\\}]
        {
            \\
            0.0  -5.2664317321355165  \\
            1.0  -5.213370338495034  \\
            2.0  -5.242531886486957  \\
            3.0  -5.080649222007  \\
        }
        ;
    \addplot[mark={*}, mark size={1.5pt}, line cap={round}, mark options={solid}, color={red}, style={dotted}]
        table[row sep={\\}]
        {
            \\
            1.0  -10.155533837232326  \\
            2.0  -17.857843242279362  \\
            3.0  -19.02475386443296  \\
        }
        ;
        \addplot[mark={*}, mark size={1.5pt}, line cap={round}, mark options={solid,fill opacity=0}, color={blue}]
        table[row sep={\\}]
        {
            \\
            0.0  -5.267697978637517  \\
            1.0  -5.249160560259795  \\
            2.0  -5.238792942191463  \\
            3.0  -5.331623440993495  \\
        }
        ;
    \addplot[mark={*}, mark size={1.5pt}, line cap={round}, mark options={solid,fill opacity=0}, color={blue}, style={dotted}]
        table[row sep={\\}]
        {
            \\
            1.0  -10.615188667948205  \\
            2.0  -12.534549504888147  \\
            3.0  -14.583841589116316  \\
        }
        ;
   \node[black] at (axis cs:0.3,-10){\large{$\alpha_{\text{p-ref.}}$=4.4}};
\node[black] at (axis cs:0.3,-12){\large{$\alpha_{\text{h-ref.}}$=1.9}};
\end{axis}
\node at (3.5,6.5) {\large---------Homogeneous cases---------};
\node at (-2,3) {\rotatebox{90}{\large------Elastic cases------}};
\end{tikzpicture}
}
 \end{subfigure}
  \hspace{2mm}
         \begin{subfigure}[b]{0.475\textwidth}
        \scalebox{0.82}{
%
%
%
%
%
%
%
\begin{tikzpicture}
\begin{axis}[ticklabel style={{font=\small}}, major tick length={2pt}, every tick/.style={{black, line cap=round}}, axis on top, legend style={{draw=none, font=\small, at={(0.03,0.03)}, anchor=south west, fill=none, legend cell align=left}}, xlabel={$\ell$}, xtick distance={1}, ylabel={$\log_2(\mathbb{E}[|\;\cdot\;|])$},ymax=-2,ymin=-20,every axis plot/.append style={thick}]
    \addplot[mark={*}, mark size={1.5pt}, line cap={round}, mark options={solid}, color={red}]
        table[row sep={\\}]
        {
            \\
            0.0  -5.376794305696258  \\
            1.0  -5.340330095662012  \\
            2.0  -5.315574036164285  \\
            3.0  -5.329236410638712  \\
        }
        ;
    \addplot[mark={*}, mark size={1.5pt}, line cap={round}, mark options={solid}, color={red}, style={dotted}]
        table[row sep={\\}]
        {
            \\
            1.0  -10.250717606997387  \\
            2.0  -14.04756741929724  \\
            3.0  -16.10109140968389  \\
        }
        ;
        \addplot[mark={*}, mark size={1.5pt}, line cap={round}, mark options={solid,fill opacity=0}, color={blue}]
        table[row sep={\\}]
        {
            \\
            0.0  -5.319850192198755  \\
            1.0  -5.284904452600748  \\
            2.0  -5.265775574097016  \\
            3.0  -5.27291122304407  \\
        }
        ;
    \addplot[mark={*}, mark size={1.5pt}, line cap={round}, mark options={solid,fill=white}, color={blue}, style={dotted}]
        table[row sep={\\}]
        {
            \\
            1.0  -10.594102345496124  \\
            2.0  -12.505514791808867  \\
            3.0  -14.29053762324863  \\
        }
        ;
  \node[black] at (axis cs:0.3,-10){\large{$\alpha_{\text{p-ref.}}$=2.8}};
\node[black] at (axis cs:0.3,-12){\large{$\alpha_{\text{h-ref.}}$=1.8}};
\end{axis}
\node at (3.5,6.5) {\large---------Heterogeneous cases---------};

\end{tikzpicture}
}
 \end{subfigure}
         \begin{subfigure}[b]{0.475\textwidth}
        \scalebox{0.82}{
%
%
%
%
%
%
%
\begin{tikzpicture}
\begin{axis}[ticklabel style={{font=\small}}, major tick length={2pt}, every tick/.style={{black, line cap=round}}, axis on top, legend style={{draw=none, font=\small, at={(0.03,0.03)}, anchor=south west, fill=none, legend cell align=left}}, xlabel={$\ell$}, xtick distance={1}, ylabel={$\log_2(\mathbb{E}[|\;\cdot\;|])$},ymin=-25,ymax=-5,every axis plot/.append style={thick}]
    \addplot[mark={*}, mark size={1.5pt}, line cap={round}, mark options={solid}, color={red}]
        table[row sep={\\}]
        {
            \\
            0.0  -7.406098701423242  \\
            1.0  -7.166924521911086  \\
            2.0  -7.14692066030635  \\
            3.0  -7.143401344457465  \\
        }
        ;
    \addplot[mark={*}, mark size={1.5pt}, line cap={round}, mark options={solid}, color={red}, style={dotted}]
        table[row sep={\\}]
        {
            \\
            1.0  -9.87719626614519  \\
            2.0  -13.300236741972325  \\
            3.0  -19.54939043058653  \\
        }
        ;
    
        \addplot[mark={*}, mark size={1.5pt}, line cap={round}, mark options={solid,fill opacity=0}, color={blue}]
        table[row sep={\\}]
        {
            \\
            0.0  -7.406359470831839  \\
            1.0  -7.232077548274022  \\
            2.0  -7.175914403543344  \\
            3.0  -7.166334483915969  \\
        }
        ;
    \addplot[mark={*}, mark size={1.5pt}, line cap={round}, mark options={solid,fill opacity =0}, color={blue}, style={dotted}]
        table[row sep={\\}]
        {
            \\
            1.0  -10.339854607845439  \\
            2.0  -11.843541177611852  \\
            3.0  -13.471735041070573  \\
        }
        ;
   \node[black] at (axis cs:0.3,-15){\large{$\alpha_{\text{p-ref.}}$=4.8}};
\node[black] at (axis cs:0.3,-17){\large{$\alpha_{\text{h-ref.}}$=1.5}};
\end{axis}
\node at (-2,3) {\rotatebox{90}{\large------Elastoplastic cases------}};

\end{tikzpicture}
}
 \end{subfigure}
  \hspace{2mm}
         \begin{subfigure}[b]{0.475\textwidth}
        \scalebox{0.82}{
%
%
%
%
%
%
%
\begin{tikzpicture}
\begin{axis}[ticklabel style={{font=\small}}, major tick length={2pt}, every tick/.style={{black, line cap=round}}, axis on top, legend style={{draw=none, font=\small, at={(0.03,0.03)}, anchor=south west, fill=none, legend cell align=left}}, xlabel={$\ell$}, xtick distance={1}, ylabel={$\log_2(\mathbb{E}[|\;\cdot\;|])$},ymin=-25,ymax=-5,every axis plot/.append style={thick}]
    \addplot[mark={*}, mark size={1.5pt}, line cap={round}, mark options={solid}, color={red}]
        table[row sep={\\}]
        {
            \\
            0.0  -7.4061322890140575  \\
            1.0  -7.1676924609317725  \\
            2.0  -7.147544843702568  \\
            3.0  -7.1482956066546866  \\
        }
        ;
    \addplot[mark={*}, mark size={1.5pt}, line cap={round}, mark options={solid}, color={red}, style={dotted}]
        table[row sep={\\}]
        {
            \\
            1.0  -9.879409480364188  \\
            2.0  -13.325249435027848  \\
            3.0  -16.375054985882123  \\
        }
        ;
    
        \addplot[mark={*}, mark size={1.5pt}, line cap={round}, mark options={solid,fill opacity =0}, color={blue}]
        table[row sep={\\}]
        {
            \\
            0.0  -7.405969373112486  \\
            1.0  -7.228402295606055  \\
            2.0  -7.1709294622713875  \\
            3.0  -7.150905849804181  \\
        }
        ;
    \addplot[mark={*}, mark size={1.5pt}, line cap={round}, mark options={solid,fill opacity=0}, color={blue}, style={dotted}]
        table[row sep={\\}]
        {
            \\
            1.0  -10.340569232182093  \\
            2.0  -11.842593403868872  \\
            3.0  -13.479402186457829  \\
        }
        ;
  \node[black] at (axis cs:0.3,-15){\large{$\alpha_{\text{p-ref.}}$=3.2}};
\node[black] at (axis cs:0.3,-17){\large{$\alpha_{\text{h-ref.}}$=1.5}};
\end{axis}
\end{tikzpicture}
}
 \end{subfigure}
 \begin{subfigure}[b]{0.99\textwidth}
\centering
        \scalebox{0.9}{
\begin{tikzpicture}
\begin{axis}[%
    hide axis,
    xmin=10,
    xmax=50,
    ymin=0,
    ymax=0.4,
    legend style={draw=white!15!black,legend cell align=left,legend columns=-1},
    every axis plot/.append style={thick}
    ]
    \addlegendimage{red,mark=*}
    \addlegendentry{$P_{\ell}$ p-ref.};
    
     \addlegendimage{red,mark=*,dotted,mark options = {solid}}
    \addlegendentry{$\Delta P_{\ell}$ p-ref.};
    
    \addlegendimage{blue,mark options = {solid,fill=white}, mark=*}
    \addlegendentry{$P_{\ell}$ h-ref.};
    
    \addlegendimage{blue,dotted,mark options = {solid,fill=white}, mark=*}
    \addlegendentry{$\Delta P_{\ell}$ h-ref.};
    \end{axis}
\end{tikzpicture}}
\vspace{5mm}
 \end{subfigure}
 \vspace{-3mm}
 \caption{Expected value of the quantity of interest and expected value of its differences for the homogeneous elastic and elastoplastic cases (top and bottom left) and the heterogeneous elastic and elastoplastic cases (top and bottom right).}
 \label{fig:expval}
\end{figure}

We first give the parameter $\gamma$ for both refinement schemes. For p-refinement, this parameter has been measured to be equal to 1.5 while for h-refinement it equals 2.0. This means that the cost increase for one solve per increasing level is larger when using h-refinement than p-refinement. This is because less dof's are added per increasing level.
\\
Fig.\,\ref{fig:variances} shows the behavior of the variance of the quantity of interest $P_\ell$, and of the difference $P_\ell-P_{\ell-1}$, in case of a tolerance $\epsilon$ equal to 3.8E-5 for the elastic cases and 2.5E-6 for the elastoplastic cases. Note that the variance of $P_\ell$ over the different levels remains constant while the variance of the differences between two successive levels continuously decreases. The rates $\beta$ are included in the figures for h- and p-refinement. These rates represent the slopes of the differences, $\Delta P_\ell$. For all but two  cases we find that $\beta > \gamma$, and thus we expect the MLMC cost to be proportional to $\epsilon^{-2}$, see Theorem 1. Only the elastoplastic cases where p-refinement is used, Fig.\,\ref{fig:variances} (bottom left and right) we find that $\beta < \gamma$. Following Theorem 1 and the results from Fig.\,\ref{fig:expval}, we calculate the cost according to $\epsilon^{-2\delta-\left(\gamma-\delta\beta\right)/\alpha}$, with $\delta = 1$. We find that for both cases the cost is approximately proportional to $\epsilon^{-2}$. We will show this in $\S$\ref{Runtime}. The MLQMC cost cannot be easily predicted due to its dependence on the factor $\delta$. We will empirically show the cost proportionality in $\S$\ref{Runtime}.
\\
\\
Furthermore we observe that the value of the variance of the differences of the p-refinement cases is larger than those of the h-refinement cases (red dashed line is lower than blue dashed line) except for the homogeneous elastoplastic case, Fig\,\ref{fig:variances} (bottom left). There, the value of the variance of the differences  is much lower for p-refinement than for h-refinement. This will lead to a larger number of samples for p-refinement and could lead to a larger computational time with respect to h-refinement. This insight can be gained by investigating Eq.\,\eqref{eq:nopt}, which calculates the optimal amount of samples on a level $\ell$ given $V_\ell$ and $C_\ell$. A lower $V_\ell$ will result in a lower optimal number of samples. This  also follows from Fig.\,\ref{fig:Times}, where we present the simulation times needed to achieve a user defined tolerance on the RMSE.


In Fig.\,\ref{fig:expval} the expected values of the quantity of interest $P_\ell$, and of the difference $P_\ell-P_{\ell-1}$ are presented. The considered tolerances are the same as for Fig.\,\ref{fig:variances}. From these figures, it is clear that the expected values of the differences over the levels, $\Delta P_\ell$, decreases faster when using p-refinement. Even while using less degrees of freedom per level, in case of p-refinement, we still obtain a very good decrease of the expected values of the differences.

\subsubsection{Number of samples}
\begin{figure}[H]
\centering
        \begin{subfigure}[b]{0.475\textwidth}
        \scalebox{0.85}{
%
%
%
%
%

\begin{tikzpicture}
\begin{axis}[ticklabel style={{font=\small}}, major tick length={2pt}, every tick/.style={{black, line cap=round}}, axis on top, legend style={{draw=none, font=\small, at={(0.03,0.03)}, anchor=south west, fill=none, legend cell align=left}}, xlabel={$\ell$}, xtick distance={1}, ylabel={$N_\ell$}, ymode={log}, legend style={{draw=none, font=\small, at={(0.97,0.97)}, anchor=north east, fill=none, legend cell align=left}},ybar,ymin=1,ymax=10^6,enlarge x limits  = 0.001,xmin=-0.5,xmax=3.5,bar width=7pt]
    \addplot[ mark size={1pt}, line cap={round}, mark options={solid}, color={red!60},style=fill]
        table[row sep={\\}]
        {
            \\
            0  187480  \\
            1  3888  \\
            2  40  \\
            3  2  \\
        }
        ;
    \addplot[ mark size={1pt}, line cap={round}, mark options={solid}, color={red!60},pattern color = red!60, pattern={north west lines}]
        table[row sep={\\}]
        {
            \\
            0  188536  \\
            1  2406  \\
            2  287  \\
            3  39  \\
        }
        ;
    
        \addplot[mark size={1pt}, line cap={round}, mark options={solid}, color={blue!60!cyan},style=fill]
        table[row sep={\\}]
        {
            \\
            0  13290  \\
            1  80  \\
            2  20  \\
            3  20  \\
        }
        ;
    \addplot[ mark size={1pt}, line cap={round}, mark options={solid}, color={blue!60!cyan},pattern color = blue!60!cyan,pattern={north west lines}]
        table[row sep={\\}]
        {
            \\
           0  15950  \\
            1  60  \\
            2  20  \\
            3  20  \\
        }
        ;
\end{axis}
\node at (3.5,6.5) {\large---------Homogeneous cases---------};
\node at (-2,3) {\rotatebox{90}{\large------Elastic cases------}};
\node at (-1.5,3){\rotatebox{90}{\large$\varepsilon=$3.8e-05}};
\end{tikzpicture}
}
 \end{subfigure}
         \begin{subfigure}[b]{0.475\textwidth}
                 \scalebox{0.85}{
%
%
%
%
%

\begin{tikzpicture}
\begin{axis}[ticklabel style={{font=\small}}, major tick length={2pt}, every tick/.style={{black, line cap=round}}, axis on top, legend style={{draw=none, font=\small, at={(0.03,0.03)}, anchor=south west, fill=none, legend cell align=left}}, xlabel={$\ell$}, xtick distance={1}, ylabel={$N_\ell$}, ymode={log}, legend style={{draw=none, font=\small, at={(0.97,0.97)}, anchor=north east, fill=none, legend cell align=left}},ymin=1,ymax=10^6,ybar,enlarge x limits  = 0.001,xmin=-0.5,xmax=3.5,bar width=7pt]
    \addplot[ mark size={1pt}, line cap={round}, mark options={solid}, color={red!60},style=fill]
        table[row sep={\\}]
        {
            \\
            0  41011  \\
            1  2697  \\
            2  1602  \\
            3  122  \\
        }
        ;
    \addplot[ mark size={1pt}, line cap={round}, mark options={solid}, color={red!60},pattern color = red!60, pattern={north west lines}]
        table[row sep={\\}]
        {
            \\
            0  77330  \\
            1  9367  \\
            2  1917  \\
            3  435  \\
        }
        ;
    
        \addplot[mark size={1pt}, line cap={round}, mark options={solid}, color={blue!60!cyan},style=fill]
        table[row sep={\\}]
        {
            \\
            0  7680  \\
            1  580  \\
            2  180  \\
            3  20  \\
        }
        ;
    \addplot[ mark size={1pt}, line cap={round}, mark options={solid}, color={blue!60!cyan},pattern color = blue!60!cyan,pattern={north west lines}]
        table[row sep={\\}]
        {
            \\
            0  22970  \\
            1  3700  \\
            2  330  \\
            3  60  \\
        }
        ;
\end{axis}
\node at (3.5,6.5) {\large---------Heterogeneous cases---------};

\end{tikzpicture}
}
 \end{subfigure}
         \begin{subfigure}[b]{0.475\textwidth}
                 \scalebox{0.85}{
%
%
%
%
%

\begin{tikzpicture}
\begin{axis}[ticklabel style={{font=\small}}, major tick length={2pt}, every tick/.style={{black, line cap=round}}, axis on top, legend style={{draw=none, font=\small, at={(0.03,0.03)}, anchor=south west, fill=none, legend cell align=left}}, xlabel={$\ell$}, xtick distance={1}, ylabel={$N_\ell$}, ymode={log}, legend style={{draw=none, font=\small, at={(0.97,0.97)}, anchor=north east, fill=none, legend cell align=left}},ymin=1,ymax=10^5,ybar,ymin=1,ymax=10^5,enlarge x limits  = 0.001,xmin=-0.5,xmax=3.5,bar width=7pt]
    \addplot[ mark size={1pt}, line cap={round}, mark options={solid}, color={red!60},style=fill]
        table[row sep={\\}]
        {
            \\
            0  15457  \\
            1  856  \\
            2  425  \\
            3  127  \\
        }
        ;
    \addplot[ mark size={1pt}, line cap={round}, mark options={solid}, color={red!60},pattern color = red!60, pattern={north west lines}]
        table[row sep={\\}]
        {
            \\
            0  11592  \\
            1  272  \\
            2  40  \\
            3  6  \\
        }
        ;
    
        \addplot[mark size={1pt}, line cap={round}, mark options={solid}, color={blue!60!cyan},style=fill]
        table[row sep={\\}]
        {
            \\
            0  6400  \\
            1  220  \\
            2  50  \\
            3  40  \\
        }
        ;
    \addplot[ mark size={1pt}, line cap={round}, mark options={solid}, color={blue!60!cyan},pattern color = blue!60!cyan,pattern={north west lines}]
        table[row sep={\\}]
        {
            \\
            0  2130  \\
            1  30  \\
            2  20 \\
            3  20  \\
        }
        ;
\end{axis}
\node at (-2,3) {\rotatebox{90}{\large------Elastoplastic cases------}};
\node at (-1.5,3){\rotatebox{90}{\large$\varepsilon=$2.5e-06}};
\end{tikzpicture}
}
 \end{subfigure}
         \begin{subfigure}[b]{0.475\textwidth}
                 \scalebox{0.85}{
%
%
%
%
%

\begin{tikzpicture}
\begin{axis}[ticklabel style={{font=\small}}, major tick length={2pt}, every tick/.style={{black, line cap=round}}, axis on top, legend style={{draw=none, font=\small, at={(0.03,0.03)}, anchor=south west, fill=none, legend cell align=left}}, xlabel={$\ell$}, xtick distance={1}, ylabel={$N_\ell$}, ymode={log}, legend style={{draw=none, font=\small, at={(0.97,0.97)}, anchor=north east, fill=none, legend cell align=left}},ymin=1,ymax=10^5,ybar,ymin=1,ymax=10^5,enlarge x limits  = 0.001,xmin=-0.5,xmax=3.5,bar width=7pt]
    \addplot[ mark size={1pt}, line cap={round}, mark options={solid}, color={red!60},style=fill]
        table[row sep={\\}]
        {
            \\
           0  4613  \\
            1  558  \\
            2  175  \\
            3  81  \\
        }
        ;
    \addplot[ mark size={1pt}, line cap={round}, mark options={solid}, color={red!60},pattern color = red!60, pattern={north west lines}]
        table[row sep={\\}]
        {
            \\
            0  6504  \\
            1  859  \\
            2  202  \\
            3  62  \\
        }
        ;
    
        \addplot[mark size={1pt}, line cap={round}, mark options={solid}, color={blue!60!cyan},style=fill]
        table[row sep={\\}]
        {
            \\
            0  1010  \\
            1  80  \\
            2  30  \\
            3  20  \\
        }
        ;
    \addplot[ mark size={1pt}, line cap={round}, mark options={solid}, color={blue!60!cyan},pattern color = blue!60!cyan,pattern={north west lines}]
        table[row sep={\\}]
        {
            \\
            0  840  \\
            1  180  \\
            2  40  \\
            3  20  \\
        }
        ;
\end{axis}
\end{tikzpicture}
}
 \end{subfigure}
 \begin{subfigure}[b]{0.99\textwidth}
\centering
\pgfplotsset{compat=1.11,
    /pgfplots/ybar legend/.style={
    /pgfplots/legend image code/.code={%
       \draw[##1,/tikz/.cd,yshift=-0.25em]
        (0cm,0cm) rectangle (3pt,0.8em);},
   },
}
\begin{tikzpicture}
\begin{axis}[%
    hide axis,
    xmin=10,
    xmax=50,
    ymin=0,
    ymax=0.4,
    ybar,
    legend style={draw=white!15!black,legend cell align=left,legend columns=-1}
    ]
    \addlegendimage{mark options={solid}, color={red!60},style=fill}
    \addlegendentry{MLMC p-ref.};
    
     \addlegendimage{mark options={solid}, color={red!60},pattern color = red!60, pattern={north west lines}}
    \addlegendentry{MLMC h-ref.};
    
    \addlegendimage{mark options={solid}, color={blue!60!cyan},style=fill}
    \addlegendentry{MLQMC p-ref.};
    
    \addlegendimage{ mark options={solid}, color={blue!60!cyan},pattern color = blue!60!cyan,pattern={north west lines}}
    \addlegendentry{MLQMC h-ref.};
    
%
    
    \end{axis}
\end{tikzpicture}
\end{subfigure}
\caption{Number of samples of MLMC and MLQMC with h- and p-refinement for the homogeneous elastic and elastoplastic cases (top and bottom left) and the heterogeneous elastic and elastoplastic cases (top and bottom right) for the finest considered tolerance.}
\label{fig:sample_sizes}
\end{figure}
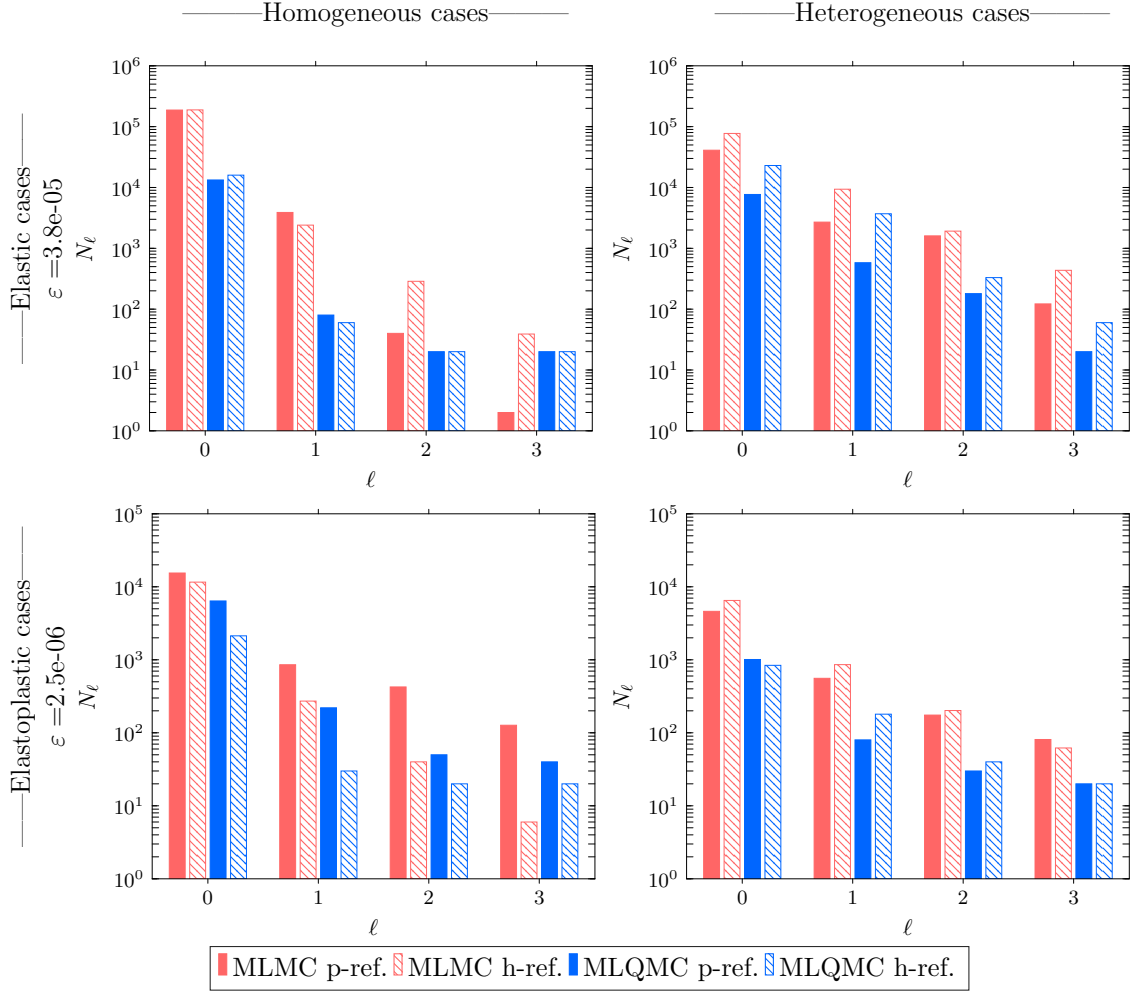

Fig.\,\ref{fig:sample_sizes} shows the total number of samples over the different levels for the elastic and elastoplastic case, for both MLMC and MLQMC. Note that the number of samples is decreasing as the level $\ell$ increases, as required. Numerical values for $N_\ell$ are repeated in Tab.\,\ref{tab:Samples_elast} and in Tab.\,\ref{tab:Samples_elastplast}. Observe that the number of samples on lower levels is higher than on the higher levels. Samples on lower levels are computationally less expensive. It is therefore advantageous to have a high number of samples on lower levels. Also, it can be seen that the sample sizes for MLQMC are lower than for MLMC. This is because the MLQMC sample points are chosen deterministically in an optimal way, see Fig.\,\ref{fig:points}. This will result in a lower computation time for MLQMC. Observe that in Tab.\,\ref{tab:Samples_elast} and in Tab.\,\ref{tab:Samples_elastplast}, all MC samples are taken on the highest level of MLMC/MLQMC. Furthermore we observe that for MC the number of samples when modeling the uncertainty as a homogeneous Young's modulus is much higher compared to a heterogeneous Young's modulus. This will lead to a considerably higher computing time.

\begin{table}[H]
\vspace{-3cm}
\scalebox{0.8}{
 \begin{tabular}{c c c c c c c c c c c c c c c}
 \toprule
\parbox[t]{10mm}{\multirow{43}{*}{\rotatebox[origin=c]{90}{\scalebox{1.4}{--------------------------------------------------------------------------------Elastic cases--------------------------------------------------------------------------------}}}} 
&\multirow{3}{*}{Level}& \multicolumn{6}{c}{Homogeneous Young's modulus} & \multicolumn{6}{c}{Heterogeneous Young's modulus}\\
\cmidrule(rl{4pt}){3-8} \cmidrule(rl{4pt}){9-14}
&&\multicolumn{3}{c}{p-ref.} & \multicolumn{3}{c}{h-ref.} & \multicolumn{3}{c}{p-ref.} & \multicolumn{3}{c}{h-ref.}\\
\cmidrule(rl{4pt}){3-5} \cmidrule(rl{4pt}){6-8} \cmidrule(rl{4pt}){9-11} \cmidrule(rl{4pt}){12-14}
&&\multicolumn{1}{c}{MLMC} & \multicolumn{1}{c}{MLQMC} & \multicolumn{1}{c}{MC} &\multicolumn{1}{c}{MLMC} & \multicolumn{1}{c}{MLQMC} & \multicolumn{1}{c}{MC}&\multicolumn{1}{c}{MLMC} & \multicolumn{1}{c}{MLQMC} & \multicolumn{1}{c}{MC}&\multicolumn{1}{c}{MLMC} & \multicolumn{1}{c}{MLQMC} & \multicolumn{1}{c}{MC}\\
 \cmidrule(){3-14}
& &\multicolumn{12}{c}{Tolerance on RMSE of 3.8E-5}\\
 \cmidrule(){3-14}
&0 & 187480 & 13290 & -    & 188536 & 15950 & -     & 41011  & 7680  & -    & 77330 & 22970 & -\\
&1 & 3888   & 80   & -     & 2406   & 60   & -      & 2697   & 580    & -    & 9367  & 3700 & -\\
&2 & 40    & 20   & -      & 287    & 20   & -      & 1602   & 180   & -    & 1917  & 330 & -\\
&3 & 2     & 20   & 188903  & 39     & 20   & /     & 122    & 20    & 29188 & 435   & 60  & 35143\\
 \cmidrule(){3-14}
& &\multicolumn{12}{c}{Tolerance on RMSE of 5.0E-5}\\
 \cmidrule(){3-14}
&0 & 112006 & 9220 & -     & 112908 & 9220 & -      & 24257 & 5330  & -   & 45849 & 19140 & -\\
&1 & 2326   & 30   & -     & 1380   & 40   & -      & 1585  & 270   & -   & 5621 & 3080 & -\\
&2 & 40    & 20   & -      & 174    & 20   & -      & 971  & 120    & -    & 1110 & 220 & -\\
&3 & 2     & 20   & 112487  & 27     & 20    & /     & 75   & 20     & 16965  & 263 & 60 & 20699\\
 \cmidrule(){3-14}
& &\multicolumn{12}{c}{Tolerance on RMSE of 6.5E-5}\\
 \cmidrule(){3-14}
&0 & 65771  & 7680 & -      & 66340  & 7680 & -    & 14568 & 2560   & -    & 27373 & 11070 & -\\
&1 & 1374    & 20   & -     & 823    & 30   & -    & 937   & 150    & -    & 3352 & 1770 & -\\
&2 & 40     & 20   & -      & 104    & 20   & -     & 599  & 120    & -    & 676 & 150 & -\\
&3 & 2      & 20   & 64569  & 13     & 20   & / & 49  & 20   & 2355 & 10170 & 30 & 12146\\
 \cmidrule(){3-14}
& &\multicolumn{12}{c}{Tolerance on RMSE of 8.4E-5}\\
 \cmidrule(){3-14}
&0 & 38107 & 2560 & -    & 38819 & 6400 & -    & 8650 & 2130 & -   & 15871 & 3700 & -\\
&1 & 803  & 20   & -     & 466   & 20   & -    &  563  & 120 & -    & 1939 & 840 & -\\
&2 & 40   & 20   & -     & 62    & 20   & -    & 362   & 80 & -     & 390 & 80 & -\\
&3 & 2    & 20   & 38479  & 8     & 20   & / & 28    & 20   & 5850 & 83 & 20 & 7282\\
 \cmidrule(){3-14}
& &\multicolumn{12}{c}{Tolerance on RMSE of 1.1E-4}\\
 \cmidrule(){3-14}
&0 & 22602 & 1470 & -    & 22895 & 3080 & -    & 5054  & 2130 & -  & 9310 & 3700 & -\\
&1 & 477  & 20   & -     & 285   & 20   & -    &  330  & 80 & -    & 1133 & 580 & -\\
&2 & 40   & 20   & -     & 40    & 20   & -    & 202    & 30 & -     & 235 & 60 & -\\
&3 & 2    & 20   & 23227  & 3     & 20   & 24078 & 15    & 20 & 3567 & 43 & 20 & 4255\\
 \cmidrule(){3-14}
& &\multicolumn{12}{c}{Tolerance on RMSE of 1.4E-4}\\
 \cmidrule(){3-14}
&0 & 13103 & 1220 & -    & 13427 & 2560 & -    & 2982 & 1470 & -  & 5172 & 2130 & -\\
&1 & 290  & 20   & -     & 156   & 20   & -    &  191  & 50 & -    & 646 & 220 & -\\
&2 & 40   & 20   & -     & 40    & 20   & -    & 115   & 20 & -     & 137 & 40 & -\\
&3 & 2    & 20   & 13565  & 2     & 20   & 13711 & 10    & 20 & 2125 & 9 & 20 & 2510\\
 \cmidrule(){3-14}
& &\multicolumn{12}{c}{Tolerance on RMSE of 1.8E-4}\\
 \cmidrule(){3-14}
&0 & 7606 & 480 & -     & 7695 & 1770 & -    & 1741 & 330 & -  & 3038 & 1470 & -\\
&1 & 173  & 20   & -    & 80   & 20   & -    &  105  & 20 & -    & 374 & 150 & -\\
&2 & 40   & 20   & -    & 40   & 20   & -    & 68   & 20 & -     & 76 & 30 & -\\
&3 & 2    & 20   & 7813 & 2    & 20   & 8171 & 7    & 20 & 1287 & 6 & 20 & 1491\\
 \cmidrule(){3-14}
& &\multicolumn{12}{c}{Tolerance on RMSE of 2.4E-4}\\
 \cmidrule(){3-14}
&0 & 4425 & 330 & -     & 4672 & 1220 & -    & 1016 & 330 & -  & 1792 & 840 & -\\
&1 & 109  & 20   & -    & 51   & 20   & -    &  63  & 20 & -    & 222 & 80 & -\\
&2 & 40   & 20   & -    & 40   & 20   & -    & 40   & 20 & -     & 47 & 20 & -\\
&3 & 2    & 20   & 4647 & 2    & 20   & 4440 & 7    & 20 & 818 & 4 & 20 & 941\\
\bottomrule
\end{tabular}}
\caption{Number of samples for MLMC, MLQMC and MC for the homogeneous and heterogeneous elastic cases.}
\label{tab:Samples_elast}
\end{table}

\begin{table}[H]
\vspace{-3cm}
\scalebox{0.8}{
 \begin{tabular}{c c c c c c c c c c c c c c c}
 \toprule
\parbox[t]{10mm}{\multirow{43}{*}{\rotatebox[origin=c]{90}{\scalebox{1.4}{--------------------------------------------------------------------------------Elastoplastic cases--------------------------------------------------------------------------------}}}} 
&\multirow{3}{*}{Level}& \multicolumn{6}{c}{Homogeneous Young's modulus} & \multicolumn{6}{c}{Heterogeneous Young's modulus}\\
\cmidrule(rl{4pt}){3-8} \cmidrule(rl{4pt}){9-14}
&&\multicolumn{3}{c}{p-ref.} & \multicolumn{3}{c}{h-ref.} & \multicolumn{3}{c}{p-ref.} & \multicolumn{3}{c}{h-ref.}\\
\cmidrule(rl{4pt}){3-5} \cmidrule(rl{4pt}){6-8} \cmidrule(rl{4pt}){9-11} \cmidrule(rl{4pt}){12-14}
&&\multicolumn{1}{c}{MLMC} & \multicolumn{1}{c}{MLQMC} & \multicolumn{1}{c}{MC} &\multicolumn{1}{c}{MLMC} & \multicolumn{1}{c}{MLQMC} & \multicolumn{1}{c}{MC}&\multicolumn{1}{c}{MLMC} & \multicolumn{1}{c}{MLQMC} & \multicolumn{1}{c}{MC}&\multicolumn{1}{c}{MLMC} & \multicolumn{1}{c}{MLQMC} & \multicolumn{1}{c}{MC}\\
 \cmidrule(){3-14}
& &\multicolumn{12}{c}{Tolerance on RMSE of 2.5E-6}\\
 \cmidrule(){3-14}
&0 & 15457 & 6400 & - & 11592 & 2130 & - & 4613 & 1010 & - & 6504 & 840 & -\\
&1 & 856 & 220    & - & 272 & 30 & -    & 558 & 80 & - & 859 & 180 & -\\
&2 & 425 & 50     & - & 40 & 20 & -     & 175 & 30 & - & 202 & 40 & -\\
&3 & 127 & 40     & / & 6 & 20 & /      & 81 & 20 & / & 62 & 20 & /\\
 \cmidrule(){3-14}
& &\multicolumn{12}{c}{Tolerance on RMSE of 3.2E-6}\\
 \cmidrule(){3-14}
&0 & 9075 & 2560 & - & 6668 & 1470 & - & 2768 & 700 & - & 3546 & 840 & -\\
&1 & 508 & 100    & - & 167 & 30 & -   & 342 & 50 & - & 471 & 180 & -\\
&2 & 244 & 40    & - & 40 & 20 & -    & 103 & 30 & - & 103 & 30 & -\\
&3 & 68 & 20     & / & 1 & 20 & /     & 52 & 20 & / & 23 & 20 & /\\
 \cmidrule(){3-14}
& &\multicolumn{12}{c}{Tolerance on RMSE of 4.2E-6}\\
 \cmidrule(){3-14}
&0 & 5457 & 2130 & - & 4070 & 840 & - & 1764 & 580 & - & 2091 & 580 & -\\
&1 & 310 & 80    & - & 100 & 20 & -   & 199 & 40 & - & 296 & 30 & -\\
&2 & 144 & 30    & - & 40 & 20 & -    & 66 & 20 & - & 63 & 20 & -\\
&3 & 44 & 20      & / & 4 & 20 & /     & 37 & 20 & 953 & 13 & 20 & /\\
 \cmidrule(){3-14}
& &\multicolumn{12}{c}{Tolerance on RMSE of 5.4E-6}\\
 \cmidrule(){3-14}
&0 & 3180 & 1010 & - & 2465 & 400 & - & 1002 & 330 & - & 1281 & 330 & -\\
&1 & 170  & 30   & - & 60 & 20 & -   & 113 & 30 & - & 175 & 40 & -\\
&2 & 47  & 20    & - & 40 & 20 & -   & 40 & 20 & - & 43 & 20 & -\\
&3 & 24  & 20    & / & 2 & 20 & /    & 18 & 20 & 513 & 9 & 20 & /\\
 \cmidrule(){3-14}
& &\multicolumn{12}{c}{Tolerance on RMSE of 7.1E-6}\\
 \cmidrule(){3-14}
&0 & 1803 & 480 & -    & 1475 & 270 & - & 607 & 180 & - & 776 & 220 & -\\
&1 & 97 & 20 & -      & 40 & 20 & -     & 71 & 20 & - & 100 & 30 & -\\
&2 & 40 & 20 & -      & 40 & 20 & -     & 40 & 20 & - & 40 & 20 & -\\
&3 & 11 & 20 & 1324   & 2 & 20 & /      & 12 & 20 & 287 & 6 & 20 & /\\
 \cmidrule(){3-14}
& &\multicolumn{12}{c}{Tolerance on RMSE of 9.2E-6}\\
 \cmidrule(){3-14}
&0 & 1077 & 180 & -   & 890 & 180 & - & 387 & 80 & - & 514 & 150 & -\\
&1 & 57 & 20 & -     & 40 & 20 & -    & 49 & 20 & - & 62 & 20 & -\\
&2 & 40 & 20 & -     & 40 & 20 & -    & 40 & 20 & - & 40 & 20 & -\\
&3 & 8 & 20 & 799    & 2 & 20 & /     & 9 & 20 & 171 & 5 & 20 & 168\\
 \cmidrule(){3-14}
& &\multicolumn{12}{c}{Tolerance on RMSE of 1.2E-5}\\
 \cmidrule(){3-14}
&0 & 653 & 100 & -   & 515 & 150 & - & 239 & 60 & - & 320 & 60 & -\\
&1 & 40 & 20 & -     & 40 & 20 & -   & 40 & 20 & - & 40 & 20 & -\\
&2 & 40 & 20 & -     & 40 & 20 & -   & 40 & 20 & - & 40 & 20 & -\\
&3 & 6 & 20 & 452    & 2 & 20 & /    & 6 & 20 & 90 & 3 & 20 & 103\\
 \cmidrule(){3-14}
& &\multicolumn{12}{c}{Tolerance on RMSE of 1.5E-5}\\
 \cmidrule(){3-14}
&0 & 389 & 100 & -   & 313 & 120 & -   & 154 & 40 & - & 211 & 50 & -\\
&1 & 40 & 20 & -    & 40 & 20 & -     & 40 & 20 & - & 40 & 20 & -\\
&2 & 40 & 20 & -    & 40 & 20 & -     & 40 & 20 & - & 40 & 20 & -\\
&3 & 4 & 20 & 285   & 2 & 20 & 288    & 5 & 20 & 60 & 3 & 20 & 64\\
\bottomrule
\end{tabular}}
\caption{Number of samples for MLMC, MLQMC and MC for the homogeneous and heterogeneous elastoplastic cases.}
\label{tab:Samples_elastplast}
\end{table}

\subsubsection{Runtime}
\label{Runtime}
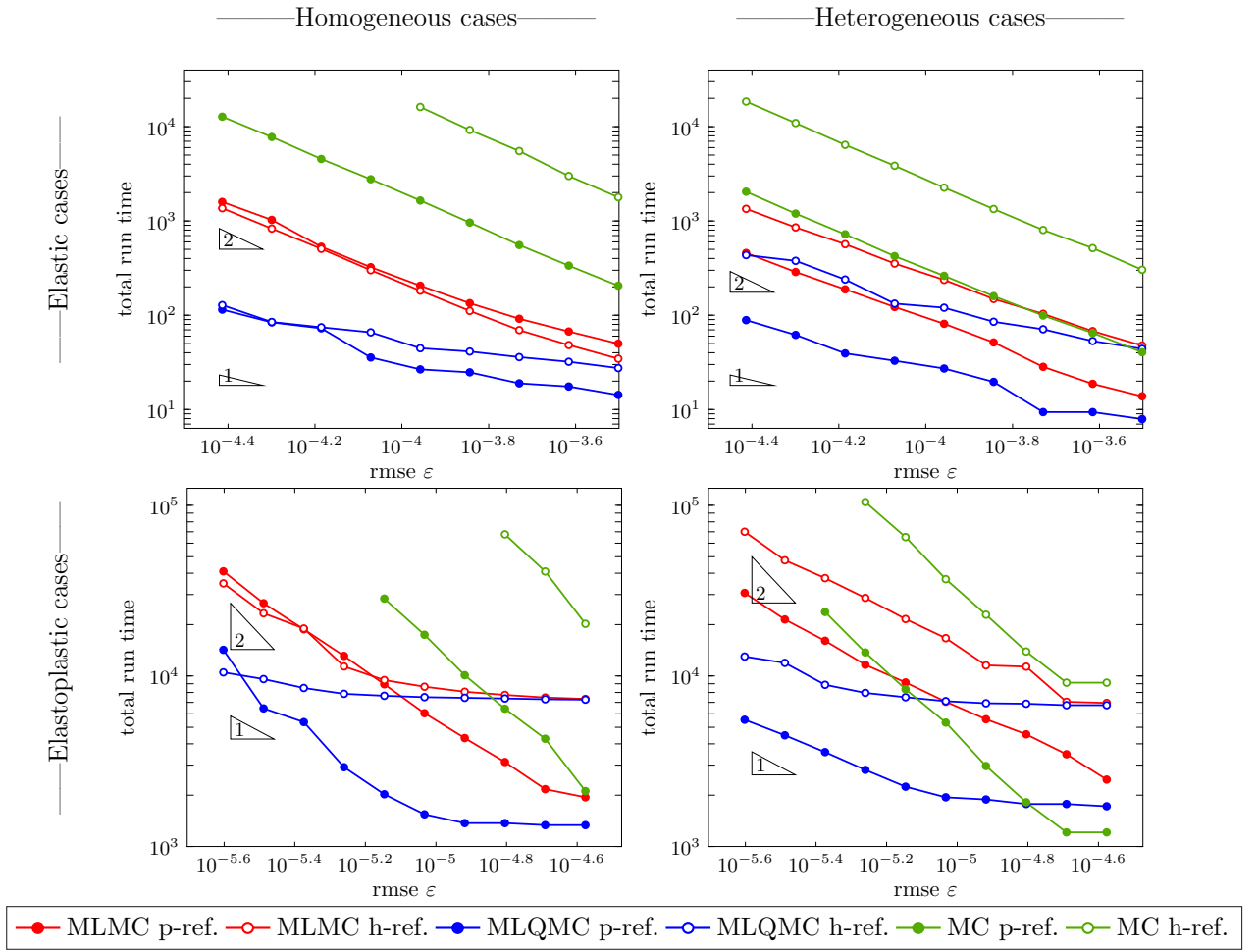
\begin{figure}[H]
\centering
        \begin{subfigure}[b]{0.475\textwidth}
        \scalebox{0.85}{
\newcommand{\littletriangle}[4]
{
    \pgfplotsextra
    {
        \pgfkeysgetvalue{/pgfplots/xmin}{\xmin}
        \pgfkeysgetvalue{/pgfplots/xmax}{\xmax}
        \pgfkeysgetvalue{/pgfplots/ymin}{\ymin}
        \pgfkeysgetvalue{/pgfplots/ymax}{\ymax}

        \pgfmathsetmacro{\xArel}{#1}
        \pgfmathsetmacro{\yArel}{#3}
        \pgfmathsetmacro{\xBrel}{#1-#2}
        \pgfmathsetmacro{\yBrel}{\yArel}
        \pgfmathsetmacro{\xCrel}{\xBrel}

        \pgfmathsetmacro{\lnxB}{\xmin*(1-(#1-#2))+\xmax*(#1-#2)} 
        \pgfmathsetmacro{\lnxA}{\xmin*(1-#1)+\xmax*#1} 
        \pgfmathsetmacro{\lnyA}{\ymin*(1-#3)+\ymax*#3} 
        \pgfmathsetmacro{\lnyC}{\lnyA+1.1*#4*(\lnxA-\lnxB)}
        \pgfmathsetmacro{\yCrel}{\lnyC-\ymin)/(\ymax-\ymin)}
        
        \coordinate (A) at (rel axis cs:\xArel,\yArel);
        \coordinate (B) at (rel axis cs:\xBrel,\yBrel);
        \coordinate (C) at (rel axis cs:\xCrel,\yCrel);

        \draw[black]   (A)--node[pos=0.9,yshift=1ex,xshift=0.5ex] {\small #4}
                    (B)--
                    (C)-- 
                    cycle;
    }
}


\begin{tikzpicture}
\begin{axis}[ticklabel style={{font=\small}}, major tick length={2pt}, every tick/.style={{black, line cap=round}}, axis on top, legend style={{draw=none, font=\small, at={(0.03,0.03)}, anchor=south west,ymin=10^0.8,ymax=10^4.6,xmin=10^-4.5,xmax=10^-3.5, fill=none, legend cell align=left}}, xmode={log}, ymode={log}, legend style={{draw=none, font=\small, at={(0.97,0.97)}, anchor=north east, fill=none, legend cell align=left}}, xlabel={rmse $\varepsilon$}, ylabel={total run time},legend style={font=\scriptsize},every axis plot/.append style={thick}]

    \addplot[mark={*}, mark size={1.5pt}, line cap={round}, mark options={solid}, color={red}]
        table[row sep={\\}]
        {
            \\
            0.0004094397207915301  43.355304058  \\
            0.00031495363137810006  49.950555596  \\
            0.00024227202413700007  67.210227667  \\
            0.00018636309549  91.741013446  \\
            0.0001433562273  134.225792282  \\
            0.00011027402100000002  205.76793256899998  \\
            8.482617e-5  323.256682813  \\
            6.525090000000001e-5  532.75565099  \\
            5.0193e-5  1026.147651703  \\
            3.861e-5  1592.054034528  \\
        }
        ;
        
                    \addplot[mark={*}, mark size={1.5pt}, line cap={round}, mark options={solid,fill = white}, color={red}]
        table[row sep={\\}]
        {
            \\
            0.0004094397207915301  26.709602477  \\
            0.00031495363137810006  34.595781675  \\
            0.00024227202413700007  48.280126593  \\
            0.00018636309549  69.529328338  \\
            0.0001433562273  111.194688304  \\
            0.00011027402100000002  182.241145812  \\
            8.482617e-5  300.24854480600004  \\
            6.525090000000001e-5  506.23101507900003  \\
            5.0193e-5  830.198222803  \\
            3.861e-5  1365.824669694  \\
        }
        ;

        \addplot[mark={*}, mark size={1.5pt}, line cap={round}, mark options={solid}, color={blue}]
        table[row sep={\\}]
        {
            \\
            0.0004094397207915301  14.273321802  \\
            0.00031495363137810006  14.279378484  \\
            0.00024227202413700007  17.528742559  \\
            0.00018636309549  18.941367595000003  \\
            0.0001433562273  24.775795671000004  \\
            0.00011027402100000002  26.625994252000005  \\
            8.482617e-5  35.62920648200001  \\
            6.525090000000001e-5  72.64560735200001  \\
            5.0193e-5  84.03720573700001  \\
            3.861e-5  114.69412520300001  \\
        }
        ;
         
                 \addplot[mark={*}, mark size={1.5pt}, line cap={round}, mark options={solid,fill = white}, color={blue}]
        table[row sep={\\}]
        {
            \\
            0.0004094397207915301  26.649936162  \\
            0.00031495363137810006  27.501116449  \\
            0.00024227202413700007  32.083091182000004  \\
            0.00018636309549  35.946382397  \\
            0.0001433562273  41.225794118  \\
            0.00011027402100000002  44.666603979  \\
            8.482617e-5  65.85655888400001  \\
            6.525090000000001e-5  74.487312903  \\
            5.0193e-5  84.718591231  \\
            3.861e-5  128.106845825  \\
        }
        ;

          \addplot[mark={*}, mark size={1.5pt}, line cap={round}, mark options={solid}, color={rgb:red,48;green,98;yellow,0}]
        table[row sep={\\}]
        {
            \\
            0.0004094397207915301  107.226865518  \\
            0.00031495363137810006  205.779184187  \\
            0.00024227202413700007  335.79354096  \\
            0.00018636309549  555.03278449  \\
            0.0001433562273  961.8892245  \\
            0.00011027402100000002  1647.888251212  \\
            8.482617e-5  2771.274196867  \\
            6.525090000000001e-5  4543.023608738  \\
            5.0193e-5  7769.837404939  \\
            3.861e-5  12755.958866338  \\
        }
        ;

          \addplot[mark={*}, mark size={1.5pt}, line cap={round}, mark options={solid,fill = white}, color={rgb:red,48;green,98;yellow,0}]
        table[row sep={\\}]
        {
            \\
            0.0004094397207915301  1719.273198663  \\
            0.00031495363137810006  1785.26724658  \\
            0.00024227202413700002  2998.9427496509998  \\
            0.00018636309549000004  5514.7133726249995  \\
            0.0001433562273  9226.226647021  \\
            0.000110274021  16177.072297863  \\
        }
        ;


    \littletriangle{0.18}{0.1}{0.5}{2};
    \littletriangle{0.18}{0.1}{0.12}{1};
\end{axis}
\node at (3.5,6.5) {\large---------Homogeneous cases---------};
\node at (-2,3) {\rotatebox{90}{\large------Elastic cases------}};
\end{tikzpicture}
}
 \end{subfigure}
         \begin{subfigure}[b]{0.475\textwidth}
                 \scalebox{0.85}{
\newcommand{\littletriangle}[4]
{
    \pgfplotsextra
    {
        \pgfkeysgetvalue{/pgfplots/xmin}{\xmin}
        \pgfkeysgetvalue{/pgfplots/xmax}{\xmax}
        \pgfkeysgetvalue{/pgfplots/ymin}{\ymin}
        \pgfkeysgetvalue{/pgfplots/ymax}{\ymax}

        \pgfmathsetmacro{\xArel}{#1}
        \pgfmathsetmacro{\yArel}{#3}
        \pgfmathsetmacro{\xBrel}{#1-#2}
        \pgfmathsetmacro{\yBrel}{\yArel}
        \pgfmathsetmacro{\xCrel}{\xBrel}

        \pgfmathsetmacro{\lnxB}{\xmin*(1-(#1-#2))+\xmax*(#1-#2)} 
        \pgfmathsetmacro{\lnxA}{\xmin*(1-#1)+\xmax*#1} 
        \pgfmathsetmacro{\lnyA}{\ymin*(1-#3)+\ymax*#3} 
        \pgfmathsetmacro{\lnyC}{\lnyA+1.1*#4*(\lnxA-\lnxB)}
        \pgfmathsetmacro{\yCrel}{\lnyC-\ymin)/(\ymax-\ymin)}
        
        \coordinate (A) at (rel axis cs:\xArel,\yArel);
        \coordinate (B) at (rel axis cs:\xBrel,\yBrel);
        \coordinate (C) at (rel axis cs:\xCrel,\yCrel);

        \draw[black]   (A)--node[pos=0.9,yshift=1ex,xshift=0.5ex] {\small #4}
                    (B)--
                    (C)-- 
                    cycle;
    }
}


\begin{tikzpicture}
\begin{axis}[ticklabel style={{font=\small}}, major tick length={2pt}, every tick/.style={{black, line cap=round}}, axis on top, legend style={{draw=none, font=\small, at={(0.03,0.03)}, anchor=south west,ymin=10^0.8,ymax=10^4.6,xmin=10^-4.5,xmax=10^-3.5, fill=none, legend cell align=left}}, xmode={log}, ymode={log}, legend style={{draw=none, font=\small, at={(0.97,0.97)}, anchor=north east, fill=none, legend cell align=left}}, xlabel={rmse $\varepsilon$}, ylabel={total run time},legend style={font=\scriptsize},every axis plot/.append style={thick}]

    \addplot[mark={*}, mark size={1.5pt}, line cap={round}, mark options={solid}, color={red}]
        table[row sep={\\}]
        {
            \\
            0.0004094397207915301  11.532561998  \\
            0.00031495363137810006  13.818767519  \\
            0.00024227202413700007  18.717934067999998  \\
            0.00018636309549  28.351934721  \\
            0.0001433562273  51.482960021  \\
            0.00011027402100000002  81.147985502  \\
            8.482617e-5  122.21362177099999  \\
            6.525090000000001e-5  188.313147468  \\
            5.0193e-5  287.628877686  \\
            3.861e-5  457.486740497  \\
        }
        ;
        
                    \addplot[mark={*}, mark size={1.5pt}, line cap={round}, mark options={solid,fill = white}, color={red}]
        table[row sep={\\}]
        {
            \\
            0.0004094397207915301  34.942879103  \\
            0.00031495363137810006  47.679816372000005  \\
            0.00024227202413700007  67.852014654  \\
            0.00018636309549  103.47732072  \\
            0.0001433562273  148.76519299699999  \\
            0.00011027402100000002  236.97102063  \\
            8.482617e-5  353.333339375  \\
            6.525090000000001e-5  567.180456519  \\
            5.0193e-5  854.345868537  \\
            3.861e-5  1344.902515102  \\
        }
        ;

        \addplot[mark={*}, mark size={1.5pt}, line cap={round}, mark options={solid}, color={blue}]
        table[row sep={\\}]
        {
            \\
            0.0004094397207915301  6.921052013  \\
            0.00031495363137810006  7.917233464  \\
            0.00024227202413700007  9.382647774  \\
            0.00018636309549  9.390247931000001  \\
            0.0001433562273  19.63172829  \\
            0.00011027402100000002  27.203467759000002  \\
            8.482617e-5  32.897927774  \\
            6.525090000000001e-5  39.369956623  \\
            5.0193e-5  61.83927878  \\
            3.861e-5  88.72247849600001  \\
        }
        ;
         
                 \addplot[mark={*}, mark size={1.5pt}, line cap={round}, mark options={solid,fill = white}, color={blue}]
        table[row sep={\\}]
        {
            \\
            0.0004094397207915301  41.229654415  \\
            0.00031495363137810006  44.013018226  \\
            0.00024227202413700007  53.249390417  \\
            0.00018636309549  70.97485190500001  \\
            0.0001433562273  85.311869192  \\
            0.00011027402100000002  120.159537865  \\
            8.482617e-5  132.80717643  \\
            6.525090000000001e-5  239.32669018299998  \\
            5.0193e-5  377.749513621  \\
            3.861e-5  434.764527445  \\
        }
        ;

          \addplot[mark={*}, mark size={1.5pt}, line cap={round}, mark options={solid}, color={rgb:red,48;green,98;yellow,0}]
        table[row sep={\\}]
        {
            \\
            0.0004094397207915301  16.216661989  \\
            0.00031495363137810006  40.503550254000004  \\
            0.00024227202413700007  64.89601855500001  \\
            0.00018636309549  99.23015766  \\
            0.0001433562273  159.533303958  \\
            0.00011027402100000002  261.958664218  \\
            8.482617e-5  423.34681298500004  \\
            6.525090000000001e-5  722.1441430750001  \\
            5.0193e-5  1197.6943270840002  \\
            3.861e-5  2042.9428903470002  \\
        }
        ;
%

          \addplot[mark={*}, mark size={1.5pt}, line cap={round}, mark options={solid,fill = white }, color={rgb:red,48;green,98;yellow,0}]
        table[row sep={\\}]
        {
            \\
            0.0004094397207915301  290.423049888  \\
            0.00031495363137810006  303.53339434199995  \\
            0.00024227202413700007  515.098949544  \\
            0.00018636309549  802.354610247  \\
            0.0001433562273  1338.6467418329999  \\
            0.00011027402100000002  2255.280942162  \\
            8.482617e-5  3851.8477988719997  \\
            6.525090000000001e-5  6414.5472711209995  \\
            5.0193e-5  10927.770025339  \\
            3.861e-5  18478.301326692  \\
        }
        ;

   \littletriangle{0.15}{0.1}{0.38}{2};
    \littletriangle{0.15}{0.1}{0.12}{1};
\end{axis}
\node at (3.5,6.5) {\large---------Heterogeneous cases---------};

\end{tikzpicture}
}
 \end{subfigure}
         \begin{subfigure}[b]{0.475\textwidth}
                 \scalebox{0.85}{
\newcommand{\littletriangle}[4]
{
    \pgfplotsextra
    {
        \pgfkeysgetvalue{/pgfplots/xmin}{\xmin}
        \pgfkeysgetvalue{/pgfplots/xmax}{\xmax}
        \pgfkeysgetvalue{/pgfplots/ymin}{\ymin}
        \pgfkeysgetvalue{/pgfplots/ymax}{\ymax}

        \pgfmathsetmacro{\xArel}{#1}
        \pgfmathsetmacro{\yArel}{#3}
        \pgfmathsetmacro{\xBrel}{#1-#2}
        \pgfmathsetmacro{\yBrel}{\yArel}
        \pgfmathsetmacro{\xCrel}{\xBrel}

        \pgfmathsetmacro{\lnxB}{\xmin*(1-(#1-#2))+\xmax*(#1-#2)} 
        \pgfmathsetmacro{\lnxA}{\xmin*(1-#1)+\xmax*#1} 
        \pgfmathsetmacro{\lnyA}{\ymin*(1-#3)+\ymax*#3} 
        \pgfmathsetmacro{\lnyC}{\lnyA+1.1*#4*(\lnxA-\lnxB)}
        \pgfmathsetmacro{\yCrel}{\lnyC-\ymin)/(\ymax-\ymin)}
        
        \coordinate (A) at (rel axis cs:\xArel,\yArel);
        \coordinate (B) at (rel axis cs:\xBrel,\yBrel);
        \coordinate (C) at (rel axis cs:\xCrel,\yCrel);

        \draw[black]   (A)--node[pos=0.9,yshift=1ex,xshift=0.5ex] {\small #4}
                    (B)--
                    (C)-- 
                    cycle;
    }
}


\begin{tikzpicture}
\begin{axis}[ticklabel style={{font=\small}}, major tick length={2pt}, every tick/.style={{black, line cap=round}}, axis on top, legend style={{draw=none, font=\small, at={(0.03,0.03)}, anchor=south west,ymin=10^3,ymax=10^5.1, fill=none, legend cell align=left}}, xmode={log}, ymode={log}, legend style={{draw=none, font=\small, at={(0.97,0.97)}, anchor=north east, fill=none, legend cell align=left}}, xlabel={rmse $\varepsilon$}, ylabel={total run time},legend style={font=\scriptsize},every axis plot/.append style={thick}]

    \addplot[mark={*}, mark size={1.5pt}, line cap={round}, mark options={solid}, color={red}]
        table[row sep={\\}]
        {
            \\
            2.651124843250001e-5  1943.631877268  \\
            2.0393268025000007e-5  2169.31571652  \\
            1.5687129250000006e-5  3128.259755313  \\
            1.2067022500000003e-5  4324.482956985  \\
            9.282325000000003e-6  6047.207132671  \\
            7.140250000000002e-6  8934.321417713  \\
            5.4925e-6  13103.70101797  \\
            4.225000000000001e-6  18774.564828884  \\
            3.2500000000000002e-6  26684.531535724  \\
            2.5e-6  41036.208305942  \\
        }
        ;
        
                    \addplot[mark={*}, mark size={1.5pt}, line cap={round}, mark options={solid,fill = white}, color={red}]
        table[row sep={\\}]
        {
            \\
            2.651124843250001e-5  7334.74122278  \\
            2.0393268025000007e-5  7474.588877085  \\
            1.5687129250000006e-5  7739.693737068  \\
            1.2067022500000003e-5  8081.740949156  \\
            9.282325000000003e-6  8648.46252427  \\
            7.140250000000002e-6  9455.155019479  \\
            5.4925e-6  11365.067162897  \\
            4.225000000000001e-6  18967.417960145  \\
            3.2500000000000002e-6  23300.715955608  \\
            2.5e-6  34879.071192204  \\
        }
        ;
        
        \addplot[mark={*}, mark size={1.5pt}, line cap={round}, mark options={solid}, color={blue}]
        table[row sep={\\}]
        {
            \\
            2.651124843250001e-5  1335.167383777  \\
            2.0393268025000007e-5  1335.172880332  \\
            1.5687129250000006e-5  1371.478678357  \\
            1.2067022500000003e-5  1371.485115278  \\
            9.282325000000003e-6  1544.585986409  \\
            7.140250000000002e-6  2024.414110664  \\
            5.4925e-6  2918.725456226  \\
            4.225000000000001e-6  5363.440396631  \\
            3.2500000000000002e-6  6452.106674758  \\
            2.5e-6  14219.104188671  \\
        }
        ;
         
                 \addplot[mark={*}, mark size={1.5pt}, line cap={round}, mark options={solid,fill = white}, color={blue}]
        table[row sep={\\}]
        {
            \\
            2.651124843250001e-5  7262.610505348  \\
            2.0393268025000007e-5  7299.965104717  \\
            1.5687129250000006e-5  7376.980824528  \\
            1.2067022500000003e-5  7442.946811694  \\
            9.282325000000003e-6  7509.0222330139995  \\
            7.140250000000002e-6  7651.506017086999  \\
            5.4925e-6  7855.527498924999  \\
            4.225000000000001e-6  8498.957997497  \\
            3.2500000000000002e-6  9582.383691111001  \\
            2.5e-6  10493.335854293002  \\
        }
        ;
         
          \addplot[mark={*}, mark size={1.5pt}, line cap={round}, mark options={solid}, color={rgb:red,48;green,98;yellow,0}]
        table[row sep={\\}]
        {
            \\
            2.6511248432500003e-5  2113.977968086  \\
            2.0393268025000004e-5  4285.15635275  \\
            1.568712925e-5  6424.318940502  \\
            1.2067022500000001e-5  10101.977542473  \\
            9.282325000000001e-6  17431.3633331  \\
            7.14025e-6  28402.381391176998  \\
        }
        ;
%

          \addplot[mark={*}, mark size={1.5pt}, line cap={round}, mark options={solid,fill = white}, color={rgb:red,48;green,98;yellow,0}]
        table[row sep={\\}]
        {
            2.6511248432500003e-5  20218.080351813  \\
            2.0393268025e-5  40953.239021055  \\
            1.568712925e-5  67602.96021550099  \\
        }
        ;

    \littletriangle{0.2}{0.1}{0.55}{2};
    \littletriangle{0.2}{0.1}{0.3}{1};
\end{axis}
\node at (-2,3) {\rotatebox{90}{\large------Elastoplastic cases------}};

\end{tikzpicture}
}
 \end{subfigure}
         \begin{subfigure}[b]{0.475\textwidth}
                 \scalebox{0.85}{
\newcommand{\littletriangle}[4]
{
    \pgfplotsextra
    {
        \pgfkeysgetvalue{/pgfplots/xmin}{\xmin}
        \pgfkeysgetvalue{/pgfplots/xmax}{\xmax}
        \pgfkeysgetvalue{/pgfplots/ymin}{\ymin}
        \pgfkeysgetvalue{/pgfplots/ymax}{\ymax}

        \pgfmathsetmacro{\xArel}{#1}
        \pgfmathsetmacro{\yArel}{#3}
        \pgfmathsetmacro{\xBrel}{#1-#2}
        \pgfmathsetmacro{\yBrel}{\yArel}
        \pgfmathsetmacro{\xCrel}{\xBrel}

        \pgfmathsetmacro{\lnxB}{\xmin*(1-(#1-#2))+\xmax*(#1-#2)} 
        \pgfmathsetmacro{\lnxA}{\xmin*(1-#1)+\xmax*#1} 
        \pgfmathsetmacro{\lnyA}{\ymin*(1-#3)+\ymax*#3} 
        \pgfmathsetmacro{\lnyC}{\lnyA+1.1*#4*(\lnxA-\lnxB)}
        \pgfmathsetmacro{\yCrel}{\lnyC-\ymin)/(\ymax-\ymin)}
        
        \coordinate (A) at (rel axis cs:\xArel,\yArel);
        \coordinate (B) at (rel axis cs:\xBrel,\yBrel);
        \coordinate (C) at (rel axis cs:\xCrel,\yCrel);

        \draw[black]   (A)--node[pos=0.9,yshift=1ex,xshift=0.5ex] {\small #4}
                    (B)--
                    (C)-- 
                    cycle;
    }
}


\begin{tikzpicture}
\begin{axis}[ticklabel style={{font=\small}}, major tick length={2pt}, every tick/.style={{black, line cap=round}}, axis on top, legend style={{draw=none, font=\small, at={(0.03,0.03)}, anchor=south west,ymin=10^3,ymax=10^5.1, fill=none, legend cell align=left}}, xmode={log}, ymode={log}, legend style={{draw=none, font=\small, at={(0.97,0.97)}, anchor=north east, fill=none, legend cell align=left}}, xlabel={rmse $\varepsilon$}, ylabel={total run time},legend style={font=\scriptsize},every axis plot/.append style={thick}]

    \addplot[mark={*}, mark size={1.5pt}, line cap={round}, mark options={solid}, color={red}]
        table[row sep={\\}]
        {
            \\
            2.651124843250001e-5  2468.171075039  \\
            2.0393268025000007e-5  3474.997050773  \\
            1.5687129250000006e-5  4546.979521533  \\
            1.2067022500000003e-5  5574.052782629  \\
            9.282325000000003e-6  7061.708707953  \\
            7.140250000000002e-6  9159.814649564  \\
            5.4925e-6  11622.122450318999  \\
            4.225000000000001e-6  16075.355716476999  \\
            3.2500000000000002e-6  21432.801040808998  \\
            2.5e-6  30666.775947962  \\
        }
        ;
        
                    \addplot[mark={*}, mark size={1.5pt}, line cap={round}, mark options={solid,fill = white}, color={red}]
        table[row sep={\\}]
        {
            \\
            2.651124843250001e-5  6947.156881139  \\
            2.0393268025000007e-5  7054.744941458  \\
            1.5687129250000006e-5  11333.878882742  \\
            1.2067022500000003e-5  11539.091617012999  \\
            9.282325000000003e-6  16647.075742087  \\
            7.140250000000002e-6  21555.283635924  \\
            5.4925e-6  28661.340946827  \\
            4.225000000000001e-6  37433.917901126  \\
            3.2500000000000002e-6  47623.086639069  \\
            2.5e-6  70039.25770138399  \\
        }
        ;
        
        \addplot[mark={*}, mark size={1.5pt}, line cap={round}, mark options={solid}, color={blue}]
        table[row sep={\\}]
        {
            \\
            2.651124843250001e-5  1721.714564519  \\
            2.0393268025000007e-5  1774.8758556340001  \\
            1.5687129250000006e-5  1774.8828537440002  \\
            1.2067022500000003e-5  1887.1269893540002  \\
            9.282325000000003e-6  1942.6437934820003  \\
            7.140250000000002e-6  2242.0538322780003  \\
            5.4925e-6  2813.753713544  \\
            4.225000000000001e-6  3576.6599857250003  \\
            3.2500000000000002e-6  4487.19511055  \\
            2.5e-6  5534.852567791  \\
        }
        ;
         
                 \addplot[mark={*}, mark size={1.5pt}, line cap={round}, mark options={solid,fill = white}, color={blue}]
        table[row sep={\\}]
        {
            \\
            2.651124843250001e-5  6739.798515872  \\
            2.0393268025000007e-5  6739.80514581  \\
            1.5687129250000006e-5  6872.78342488  \\
            1.2067022500000003e-5  6916.1945967579995  \\
            9.282325000000003e-6  7113.164792803999  \\
            7.140250000000002e-6  7502.960319489  \\
            5.4925e-6  7956.37446804  \\
            4.225000000000001e-6  8872.369791083  \\
            3.2500000000000002e-6  11923.374043865999  \\
            2.5e-6  12990.459557795999  \\
        }
        ;
         
          \addplot[mark={*}, mark size={1.5pt}, line cap={round}, mark options={solid}, color={rgb:red,48;green,98;yellow,0}]
        table[row sep={\\}]
        {
            \\
            2.651124843250001e-5  1210.778605844  \\
            2.0393268025000004e-5  1210.7792045469998  \\
            1.5687129250000003e-5  1816.7900774549998  \\
            1.2067022500000003e-5  2964.9727783219996  \\
            9.282325000000001e-6  5326.888134027  \\
            7.140250000000001e-6  8358.069108357  \\
            5.4925e-6  13733.332339828  \\
            4.225e-6  23708.885815717003  \\
        }
        ;

          \addplot[mark={*}, mark size={1.5pt}, line cap={round}, mark options={solid,fill = white}, color={rgb:red,48;green,98;yellow,0}]
        table[row sep={\\}]
        {
            \\
            2.6511248432500006e-5  9134.056101813  \\
            2.0393268025000004e-5  9134.056500873  \\
            1.5687129250000003e-5  13905.391338717  \\
            1.2067022500000001e-5  22861.965287123  \\
            9.282325000000001e-6  36925.281963841  \\
            7.140250000000001e-6  65153.201246959  \\
            5.4925e-6  104519.58169390599  \\
        }
        ;

    \littletriangle{0.2}{0.1}{0.68}{2};
    \littletriangle{0.2}{0.1}{0.2}{1};
\end{axis}
\end{tikzpicture}
}
 \end{subfigure}
 \begin{subfigure}[b]{0.99\textwidth}
\centering
\begin{tikzpicture}
\begin{axis}[%
    hide axis,
    xmin=10,
    xmax=50,
    ymin=0,
    ymax=0.4,
    legend style={draw=white!15!black,legend cell align=left,legend columns=-1},
    every axis plot/.append style={thick}
    ]
    \addlegendimage{red,mark=*}
    \addlegendentry{MLMC p-ref.};
    
     \addlegendimage{red,mark=*,mark options = {solid,fill=white}}
    \addlegendentry{MLMC h-ref.};
    
    \addlegendimage{blue,mark options = {solid}, mark=*}
    \addlegendentry{MLQMC p-ref.};
    
    \addlegendimage{blue,mark options = {solid,fill=white}, mark=*}
    \addlegendentry{MLQMC h-ref.};
    
    \addlegendimage{color={rgb:red,48;green,98;yellow,0},mark options = {solid}, mark=*}
    \addlegendentry{MC p-ref.};
    
    \addlegendimage{color={rgb:red,48;green,98;yellow,0},mark options = {solid,fill=white}, mark=*}
    \addlegendentry{MC h-ref.};    
    
    \end{axis}
\end{tikzpicture}
\end{subfigure}
\caption{Total simulation runtime of MC, MLMC and MLQMC with h- and p-refinement for the homogeneous elastic and elastoplastic cases (top and bottom left) and the heterogeneous elastic and elastoplastic cases (top and bottom right).}
\label{fig:Times}
\end{figure}

We plot the runtimes  for the different Monte Carlo methods combined with p-refinement and with h-refinement for a homogeneous and a heterogeneous modulus in Fig.\,\ref{fig:Times}. Here, the actual simulation time needed to reach a certain tolerance $\epsilon$ on the root-mean-square error (RMSE) for standard MC, MLMC and MLQMC is compared. In Tab.\,\ref{table:times}, we summarize these results. 
Note that the MC simulation is run at the highest level $L$ of the corresponding MLMC/MLQMC simulation, where $L$ is chosen according to Eq.\,\eqref{eq:bias_constraint}.   As can be seen, not all tolerances are simulated for the MC simulations. These simulations have not been done due to the long computation time that would be necessary, i.e., several days.

The first observation to be made is that the MLMC and MLQMC  simulations consistently outperform the MC simulations in terms of computational speed, except for one case in the pre-asymptotic phase (low tolerances), Fig.\ref{fig:Times} (bottom right). Speedups up to a factor 100 are observed. MLQMC outperforms MLMC by a factor  5 to 10. We observe that MLQMC tends to work well for cases where few uncertainties are considered, i.e., the homogeneous cases, Fig.\,\ref{fig:Times} (bottom and top left). While performing well, the cost for MLQMC for these cases tends to be proportional to $\epsilon^{-2}$ instead of the optimal value of $\epsilon^{-1}$. For the heterogeneous cases, Fig.\,\ref{fig:Times} (bottom and top right), we observe indeed that the cost for MLQMC is proportional to  $\epsilon^{-1}$. All MLMC costs are proportional to $\epsilon^{-2}$. This can be seen by investigating the starting points of the lines in combination with the triangles indicating the slopes.
\\
An interesting observation from Fig.\,\ref{fig:Times}, is that for half of the considered cases, it is advantageous to apply MLMC or MLQMC with a p-hierarchy of mesh refinements instead of h-refinements. All heterogeneous cases perform faster with p-refinement, Fig.\,\ref{fig:Times} (bottom and top right). However, for the homogeneous cases, MLMC with p-refinement yields no lower simulation time, Fig.\,\ref{fig:Times} (bottom and top left). While the simulation time for low tolerances is indeed lower for MLMC combined with p-refinement. For higher tolerances no gain is to be found for using one refinement method over the other. This seems only to be a problem when a very low number of uncertainties (homogeneous Young's modulus) is considered. MLQMC performs well but as already elaborated upon above, its cost is more likely to be proportional to $\epsilon^{-2}$ instead of $\epsilon^{-1}$. This is especially visible in Fig.\,\ref{fig:Times} (bottom left).
\\
When using the MC method, all homogeneous cases  have a higher simulation time than the heterogeneous cases. When using the MLMC/MLQMC method, the simulation times of the homogeneous cases are roughly equal to the ones of the heterogeneous cases. In general the number of samples of the homogeneous cases on the lowest level are an order of magnitude larger than the number of samples of the heterogeneous cases. For higher levels, the number of samples is higher for the heterogeneous cases. This creates a balancing effect resulting in both cases having roughly the same simulation time.

We thus conclude that MLMC and MLQMC achieve speedups up to factor of 100 with respect to standard MC. We have  empirically demonstrated that it is possible for MLQMC, applied to a structural engineering problem, to achieve an optimal cost of $\epsilon^{-1}$, under certain conditions.  Also, a p-refinement scheme is highly advantageous for problems where a high number of uncertainties are present.

\begin{table}[H]
\centering
  \scalebox{0.8}{
 \begin{tabular}{cccccccccccccc}
 \toprule
   & & \multicolumn{12}{c}{Time [sec]}  \\
   \cmidrule(rl{4pt}){3-14}  
& \multirow{3}{*}{RMSE [/]} & \multicolumn{6}{c}{Homogeneous Young's modulus}   & \multicolumn{6}{c}{Heterogeneous Young's modulus}  \\ [0.5ex]
 \cmidrule(rl{4pt}){3-8} \cmidrule(rl{4pt}){9-14} 
 & & \multicolumn{3}{c}{p-ref.} & \multicolumn{3}{c}{h-ref.} & \multicolumn{3}{c}{p-ref.} & \multicolumn{3}{c}{h-ref.}\\
  \cmidrule(rl{4pt}){3-5} \cmidrule(rl{4pt}){6-8} \cmidrule(rl{4pt}){9-11} \cmidrule(rl{4pt}){12-14}
 &  & MLMC & MLQMC& MC & MLMC &MLQMC& MC & MLMC &MLQMC& MC & MLMC &MLQMC& MC \\
\cmidrule{2-14}
\parbox[t]{2mm}{\multirow{10}{*}{\rotatebox[origin=c]{90}{Elastic cases}}} 
& 4.0E-4 & 43 & 14 & 107 & 26        & 26 & 1719  & 11   & 6  & 16  & 34   & 41 & 290\\
& 3.1E-4 & 49 & 14 & 205 & 34        & 27 & 1785  & 13   & 7  & 40  & 47   & 44 & 303\\
& 2.4E-4 & 67 & 17 & 335 & 48         & 32 & 2998 & 18   & 9  & 64  & 67   & 53 & 515\\
& 1.8E-4 & 91 & 18 & 555 & 69         & 35 & 5514 & 28   & 9  & 99  & 103  & 70 & 802\\
& 1.4E-4 & 134 & 24 & 961 & 111      & 41 & 9226  & 51   & 19 & 159  & 148  & 85 & 1338\\
& 1.1E-4 & 205 & 26 & 1647 & 182     & 44 & 16177 & 81   & 27 & 261  & 236  & 120 & 2255\\
& 8.4E-5 & 323 & 35 & 2771 & 300     & 65 & /     & 122  & 32 & 423 & 353  & 132 & 3851\\
& 6.5E-5 & 532 & 72 & 4543 & 506     & 74 & /     & 188  & 39 & 722 & 567  & 239 & 6414\\
& 5.0E-5 & 1026 & 84 & 7769 & 830    & 84 & /     & 287  & 61 & 1197 & 854  & 377 & 10927\\
& 3.8E-5 & 1592 & 114 & 12755 & 1365 & 128 & /     & 457 & 88 & 2042 & 1344 & 434 & 18478\\
 \cmidrule{2-14}
\parbox[t]{2mm}{\multirow{10}{*}{\rotatebox[origin=c]{90}{Elastoplastic cases}}} 
& 2.6E-5 & 1943 & 1335 & 2113 & 7334 & 7262 & 20218 & 2468 & 1721 & 1210 & 6947 & 6739 & 9134\\
& 2.0E-5 & 2169 & 1335 & 4285 & 7474 & 7299 & 40953 & 3474 & 1774 & 1210 & 7054 & 6739 & 9134\\
& 1.5E-5 & 3128 & 1371 & 6424 & 7739 & 7376 & 67602 & 4546 & 1774 & 1816 & 11333 & 6872 & 13905\\
& 1.2E-5 & 4324 & 1371 & 10101 & 8081 & 7442 & /    & 5574 & 1887 & 2964 & 11539 & 6916 & 22861\\
& 9.2E-6 & 6047 & 1544 & 17431 & 8648 & 7509 & /    & 7061 & 1942 & 5326 & 16647 & 7113 & 36925\\
& 7.1E-6 & 8934 & 2024 & 28402 & 9455 & 7651 & /    & 9159 & 2242 & 8358 & 21555 & 7502 & 65153\\
& 5.4E-6 & 13103 & 2918 & /    & 11365 & 7855 & /   & 11622 & 2813 & 13733 & 28661 & 7956 & 104519\\
& 4.2E-6 & 18774 & 5363 & /    & 18967 & 8498 & /   & 16075 & 3576 & 23708 & 37433 & 8872 & /\\
& 3.2E-6 & 26684 & 6452 & /    & 23300 & 9582 & /   & 21432 & 4487 & /     & 47623 & 11923 & /\\
& 2.5E-6 & 41036 & 14219 & /    & 34879 & 10493 & /   & 30666 & 5534 & /    & 70039 & 12990 & /\\
  \bottomrule
\end{tabular}}
\caption{Actual simulation time in seconds for MC, MLMC and MLQMC for the homogeneous elastic and elastoplastic cases and the heterogeneous elastic and elastoplastic cases.}
\label{table:times}
\end{table}

\section{Conclusion}
In this work, we considered a structural engineering problem where the uncertainty resides in the Young's modulus. To model this uncertainty, we considered both a homogeneous model, represented by a single random variable, and a heterogeneous model, represented by means of a random field. The stochastic responses were computed by means of the MC, the MLMC and the MLQMC method. We considered a mesh hierarchy based on h-refinement and on p-refinement for each method. In a first step we illustrated that the nature of the uncertainty in the Young's modulus has a major impact on the uncertainty characteristics of the simulation results. For realistic applications, the appropriate choice of the uncertainty model is of great importance. Then, we demonstrated that the MLMC method provides a significant computational cost reduction and speedup compared to the standard MC method up to a factor 100 regardless of the mesh hierarchy used. This has been shown by means of actual computing times. We further compared the speedup of MLQMC with respect to MLMC, and found speedups ranging from  5 to 10. When dealing with many uncertainty parameters, i.e., a heterogeneous Young's modulus, the multilevel methods can be accelerated even further by using a p-refinement mesh hierarchy. 
For MLQMC, we empirically showed that for many uncertainty parameters, the optimal cost in function of a desired tolerance is proportional to $\epsilon^{-1}$  with either an h- or a p-refinement mesh hierarchy. While for one uncertainty, i.e., the homogeneous Young's modulus, we showed that MLQMC combined with a p-refinement mesh hierarchy has a cost proportional to $\epsilon^{-2}$, the same as for all cases computed with MLMC.

Further paths of research will focus on ways to exploit the similarities between the responses for neighboring frequencies and to combine the advantage of both the h- and p-refinement mesh hierarchies in a Multi-Index setting \cite{PJ,PJ2}.

\section*{Acknowledgements}

This research was funded by project IWT/SBO EUFORIA: ``Efficient Uncertainty quantification For Optimization in Robust design of Industrial Applications" (IWT-140068) of the Agency for Innovation by Science and Technology, Flanders, Belgium. The authors gratefully acknowledge the support from the research council of KU Leuven though the funding of project C16/17/008 ``Efficient methods for large-scale PDE-constrained optimization in the presence of uncertainty and complex technological constraints". The authors also would like to thank the Structural Mechanics Section of the KU Leuven and Jef Wambacq for supplying their StaBIL code and providing support.

\section*{References}

\bibliography{Bib_ref.bib}
%
%
%
%

\end{document}